\input amstex
\loadbold
\loadeurm
\loadeusm
\openup.8\jot \magnification=1200

\def\sgn{{\text{sgn}}}
\def\g{{\bold g}}
\def\G{{\Cal G}}

\def\cls{{\text{ cls}}}
\def\Z{\Bbb Z}
\def\A{\Bbb A}
\def\F{\Bbb F}
\def\E{\Bbb E}
\def\disc{{\text{ disc}}}
\def\rank{{\text {rank}}}

\def\sym{{\text{sym}}}
\def\spn{{\text{span}}}

\def\H{{\Bbb H}}
\def\rad{{\text {rad}}}
\def\e{{\text{e}}}

\def\h{{\Cal H}}
\def\R{{\Bbb R}}
\def\C{{\Bbb C}}
\def\Q{{\Bbb Q}}

\def\gen{{\text{gen}}}
\def\cls{{\text{cls}}}

\def\lcm{{\text{lcm}}}

\def\B{{\Cal{B}}}
\def\Eis{{\Cal{E}}}

\def\stufe{{\Cal N}}

\documentstyle{amsppt}
\pageheight{7.7in}
\vcorrection{-0.05in}
\topmatter
\pageno 1
\title 
Toward explicit formulas for higher representation
numbers of quadratic forms
\endtitle
\rightheadtext{Higher representation numbers of quadratic forms}
\author Lynne H. Walling\endauthor
\address L.H. Walling, Department of Mathematics, University Walk,
University of Bristol, Bristol BS8 1TW England\endaddress

\abstract It is known that average Siegel theta series lie in the space of
Siegel Eisenstein series.  Also, every lattice equipped with an
even integral quadratic form
lies in a maximal lattice.  Here we consider average Siegel theta series
of degree 2 attached to maximal lattices; we
construct maps for which the average 
theta series is an eigenform.  We evaluate the
action of these maps on Siegel Eisenstein series of degree 2 (without knowing
their Fourier coefficients), and then realise the average theta series
as an explicit linear combination of the Eisenstein series.
\endabstract

\endtopmatter
\document

\head{\bf \S1. Introduction}\endhead
\smallskip
Quadratic forms are ubiquitous in mathematics and physics, as they capture geometry on
vector spaces.
In number theory we focus on questions about the integers; thus
given a vector space $V$ equipped with a quadratic form $Q$, we are interested in lattices
$K\approx\Z^m$ contained in $V$.  Siegel asked, with
$T$ another quadratic form of dimension $n\le m$, on how many sublattices of $K$
does $Q$ restrict to $T$?  In his Hauptsatz, Siegel used analytic number theory's
circle method to show that, ``on average'', the number of times $K$ represents $T$
is given by a product of $p$-adic densities.  
Siegel also introduced generalised theta
series whose Fourier coefficients are these representation numbers; he proved that 
in the case that $Q$ is positive definite and $m\ge 2n+3$, the
average theta series (averaging over the genus of $K$) lies in the space of Siegel
Eisenstein series.  (This latter result is elegently reproved by Freitag in [3]
in the case $m$ is even.)

In this paper we take these results one step further,
showing that with $K$ a maximal lattice of odd level $\stufe$
and rank $2k\ge 8$,
$$\theta(gen K)=\sum_{\stufe_0\stufe_1\stufe_2=\stufe}c_1(\stufe_1)c(\stufe_2)\E_{(\stufe_0,\stufe_1,\stufe_2)}$$
where the $\E_{(\stufe_0,\stufe_1,\stufe_2)}$ form a natural basis for the space of degree
2 Siegel Eisenstein series containing $\theta(\gen K)$, and 
$c_1(\stufe_1), c_2(\stufe_2)$ are multiplicative functions such that, for $q$ prime and $\chi$ the
character associated to $\theta(K)$,
$$c_1(q)=\cases -q^{-1}&\text{if $\chi_q=1$,}\\
\chi_{\stufe/q}(q)\left(\frac{\eta_q}{q}\right)q^{-1}\g(q)&\text{if $\chi_q\not=1$,}
\endcases$$
where $\Z_qK\simeq2\big<1,-1,\ldots,1,-1,\eta_q\big>\perp2q\big<\eta'_q\big>$ when $\chi_q\not=1$
and $\g(q)$ is the classical Gauss sum,
and 
$$c_2(q)=\cases q^{-2}&\text{if $\chi_q=1$,}\\
\left(\frac{-1}{q}\right)q^{-1}&\text{if $\chi_q\not=1$.}
\endcases$$

Our strategy here is to construct from $K$ lattices in the genus of $K$, and
then to use maps on modular forms to
mimic these lattice constructions.  Thus for each prime $q|\stufe$ we build a map $T_K(q)$
for which the average theta series $\theta(\gen K)$
is an eigenform (note that this was also our strategy in [16] where we considered
the case of Siegel degree 1).  The map $T_K(q)$ is an algebraic combination of Hecke operators, the standard
shift operator, as well as
several ``twist'' operators (defined in \S3).  
As the shift  operator associated to $q$ maps modular forms
of square-free level $\stufe$ to modular forms of level $q\stufe$,
we need to compute the action of the Hecke and twist operators on
Eisenstein series of level $q\stufe$; note that this is done without
having any knowledge of
the Fourier coefficients of these Eisenstein series.
(These computations comprise the bulk of the work in this paper.)  
The fact that $T_K(q)$ raises the level
gives us the leverage we need to show that there is a unique Eisenstein
series $\E$ of level $\stufe$ and with 0th Fourier coefficient 1 so that
for each prime $q|\stufe$, $\E|T_K(q)$ has level $\stufe$.

Together with closed-form formulas for the Fourier coefficients of the Siegel Eisenstein
series, our result herein will produce closed-form formulas for the average representation numbers
of maximal lattices of odd level.  Many authors have produced closed-form formulas
for the Fourier coefficients of Eisenstein series under various conditions
(for instance, 
see [8], [10], [11] for level 1, degree 2; [6] for level 1, degree 3;
[1], [2], [7], [9] for level 1, arbitrary degree; [12] for degree 2, odd square-free level $\stufe$,
primitive character, and the Eisenstein series associated to
$\Gamma_{\infty}\backslash\Gamma_0(\stufe)$; [15] for degree 2, arbitrary level $\stufe$, primitive character,
and the Eisenstein series associated to
$\Gamma_{\infty}\backslash\Gamma_0(\stufe)$).
So as yet, we do not have such formulas for the Fourier coefficients of a basis
for the space of Siegel Eisenstein series containing $\theta(\gen K)$,
but this is an active area of research.

We should note that Siegel's Hauptsatz is not just restricted to positive definite
quadratic forms.
When $Q$ is indefinite, Siegel considered the ``measure'' of the
representation, measuring the density of subspaces of $L$ on which $Q$ restricts to
$T$.  Siegel also introduced (non-holomorphic) theta series whose Fourier coefficients
reflect these densities, but the situation here is noteably more complicated (see
[4] for a derivation of explicit formulas for these measures in the case of 
Siegel degree 1).

\bigskip
\head{\bf \S2. Preliminaries}\endhead
\smallskip

Given a dimension $2k$
vector space $V$ with quadratic form $Q$, we say a lattice $K\approx\Z^{2k}$
in $V$ is even integral if $Q(K)\subseteq 2\Z$; note that we can always scale
$Q$ to accomplish this.  The symmetric bilinear form $B$ that we associate to $Q$
is given by the relation $Q(x+y)=Q(x)+Q(y)+2B(x,y)$ (so $Q(x)=B(x,x)$).
Let $(v_1,\ldots,v_{2k})$ be a $\Z$-basis for $K$; slightly abusing notation,
we write $Q$ for the matrix $\big(B(v_i,v_j)\big)$, which represents the quadratic
form relative to the given basis, and we write $K\simeq (B(v_i,v_j))$.
Note that since $Q(K)\subseteq 2\Z$, $Q$ is an even integral matrix
(meaning integral with even diagonal entries).

We assume from now on that $Q$ is positive definite with $k\ge 4$, and that $K$ is a maximal
even integral lattice that is unimodular at 2
(meaning $Q$ is invertible over $\Z_2$).  The level of $K$ is the smallest
positive integer $\stufe$ so that $\stufe Q^{-1}$ is also even integral.  
This means that for $p$ prime, $p\nmid \stufe$
if and only if $\Z_pK$ is unimodular; since $K$ is maximal, for a prime $q|\stufe$ we have
$$\Z_qK\simeq
\cases 2\big<1,-1,\ldots,1,-1,\eta\big>\perp 2q\big<\eta'\big>
&\text{if $q\parallel\stufe$,}\\
2\big<1,-1,\ldots,1,-1,1,-\omega\big>\perp 2q\big<1,-\omega\big>
&\text{if $q^2|\stufe$}\endcases$$
where $\eta,\eta'\in\Z_q^{\times}$ and
$\omega\in\Z$ is chosen so that $\left(\frac{\omega}{q}\right)=-1$ and
$\big<*,\ldots,*\big>$ denotes a diagonal matrix.

The degree 2 Siegel theta series is defined as follows.
For 
$$\tau\in\h_{(2)}\{X+iY:\ X,Y\in\R^{2,2}_{\sym},\ Y>0\ \},$$
 set
$$\theta(L;\tau)=\sum_{C}\e\{\,^tCQC\tau\}$$
where $C$ varies over $\Z^{2k,2}$ and $\e\{*\}=\exp(\pi i Tr(*))$.
(Here $Y>0$ means that, as a quadratic form, $Y$ is positive definite.)
Since $Q$ is positive definite, this series is absolutely uniformly convergent
on compact regions,
and hence is an analytic function (in all variables of $\tau$).  As a Fourier series,
we have
$$\theta(L;\tau)=\sum_T r(Q,T)\e\{T\tau\}$$
where $T$ varies over all positive, semi-definite, even integral $2\times 2$ matrices,
and $$r(Q,T)=\{C\in\Z^{2k,2}:\ ^tCQC=T\ \},$$
which is the number of times $Q$ represents $T$.  For $C\in\Z^{2k,2}$, let
$(y_1\ y_2)=(v_1\ \ldots\ v_{2k})C$
and set $\Lambda=\Z y_1+\Z y_2$; thus $\Lambda$ is a sublattice
of $L$ with (formal) rank 2 and quadratic form
$^tCQC$.  Hence we can rewrite the theta series as
$$\theta(L;\tau)=\sum_{\Lambda}\e\{\Lambda\tau\}$$
where $\Lambda$ varies over all sublattices of $L$ of (formal) rank 2, and with
$T_{\Lambda}$ a matrix for $Q$ restricted to $\Lambda$,
$$\e\{\Lambda\tau\}=\sum_{G}\e\{\,^tGT_{\Lambda}G\tau\};$$
when $k$ is even, $G$ varies over $GL_2(\Z)$, and when $k$ is odd, we equip
$\Lambda$ with an orientation and $G$ varies over $SL_2(\Z)$.

This generalised theta series is a prototypical example of a 
degree 2 Siegel modular form of weight $k$ and level $\stufe$,
meaning that for $\pmatrix A&B\\C&D\endpmatrix\in Sp_2(\Z)$ with $\stufe|C$,
$$\theta(K;(A\tau+B)(C\tau+D)^{-1})=\chi(\det D)\,\det(C\tau+D)^k\,\theta(K;\tau).$$
Here $\chi$ is the character modulo $\stufe$ given by
$$\chi(d)=(\sgn d)^k\left(\frac{(-1)^k\disc K}{|d|}\right)$$
where $\disc K=\det Q=\det \big(B(v_i,v_j)\big)$ and $\left(\frac{*}{*}\right)$
is the Jacobi symbol extended so that
$\left(\frac{*}{2}\right)$ denotes the Kronecker symbol.

Since $K$ is even integral and unimodular at 2, we have
$$\Z_2K\simeq\pmatrix 0&1\\1&0\endpmatrix\perp\ldots\perp\pmatrix 0&1\\1&0\endpmatrix
\text{ or } 
\pmatrix 0&1\\1&0\endpmatrix\perp\ldots\perp\pmatrix 0&1\\1&0\endpmatrix\perp\pmatrix 2&1\\1&2\endpmatrix.$$
Thus $(-1)^k\disc K\equiv1\ (4)$, and given
the structure of $\Z_qK$ for $q|\stufe$,
we have $\disc K=\Delta_0\Delta_1^2$ where $\Delta_0,\Delta_1\in\Z$; so $(-1)^k\Delta_0\equiv1\ (4)$, and
since $Q$ is
positive definite, $\Delta_0>0$.
Note that by
Quadratic Reciprocity and the fact that $\left(\frac{2}{\Delta_0}\right)=(-1)^{(\Delta_0^2-1)/2}$,
for $p$ prime we have
$\chi(p)=\left(\frac{(-1)^k\disc K}{p}\right)
=\left(\frac{p}{\Delta_0}\right).$ 
Hence $\chi(d)=\left(\frac{d}{\disc K}\right)$, and for  a prime $q|\stufe$,
$\chi_q=1$ if and only if
$$\Z_qK\simeq2\big<1,-1,\ldots,1,-1,1,-\omega\big>\perp2q\big<1,-\omega\big>.$$

Suppose $K,M$ are lattices on $V$ (meaning $K,M$ are lattices in $V$ with $\Q K=\Q M=V$).  
Then it is well-known that there is a $\Z$-basis $v_1,\dots,v_{2k}$ for $K$
and rational numbers $r_1,\ldots,r_{2k}$ so that
$r_1v_1,\ldots,r_{2k}v_{2k}$ is a $\Z$-basis for $M$ with $r_{i+1}\in\Z r_i$;
the numbers $r_i$ are called the invariant factors (or elementary divisors) of $M$ in $K$,
and we write $\{K:M\}$ to denote $(r_1,\ldots,r_{2k})$.
We say $K$ and $M$ are isometric,
denoted $K\simeq M$, if there is an isomorphism $\sigma$ from $K$ onto $M$
so that $B(\sigma x,\sigma y)=B(x,y)$ for every $x,y\in K$.
The orthogonal group of $K$, denoted
$O(K)$, is the set of isometries from $K$ to itself.  
We say $M$ is in the genus of $K$, denoted $\gen K$, if for 
every prime $p$, $\Z_pM\simeq\Z_pK$.

Since $Q$ is positive definite, it is known that there are finitely many
isometry classes within $\gen K$, and that the orthogonal group for any
lattice in $V$ is finite.
We define the 
average theta series to be
$$\theta(\gen K)=\frac{1}{\text{mass}\,K}\sum_{\cls M\in\gen K}\frac{1}{o(M)}\theta(M)$$
where $\cls M$ varies over all the isometry classes in $\gen K$,
$o(M)$ is the order of $O(M)$, and
$$\text{mass}\,K=\sum_{\cls M\in\gen K}\frac{1}{o(M)}.$$
Hence the 0th Fourier coefficient of $\theta(\gen K)$ is 1, and the $T$th coefficient
is the average number of times $T$ is represented by the lattices in the genus of $K$.

We know that $\theta(\gen K)$ lies in the space of degree 2 Siegel Eisenstein
series of weight $k$, level $\stufe$, and character $\chi$ (see [3]).
The elements in a natural basis for this space correspond to elements of
$\Gamma_{\infty}\backslash Sp_2(\Z)\slash\Gamma_0(\stufe)$ where
$$\Gamma_{\infty}=\left\{\gamma\in Sp_2(\Z):\ \gamma=\pmatrix *&*\\ 0&*\endpmatrix\ \right\}$$
and 
$$\Gamma_0(\stufe)=\left\{\gamma\in Sp_2(\Z):\ \gamma\equiv \pmatrix *&*\\ 0&*\endpmatrix\ (\stufe)\ \right\}.$$
For $\gamma_0\in Sp_2(\Z)$, the Eisenstein series associated to 
$\Gamma_{\infty}\gamma_0\Gamma_0(\stufe)$ is 
$$\E_{\gamma_0}(\tau)=\sum\overline\chi(\det D_{\gamma})\,1|\gamma(\tau)$$
where $\Gamma_{\infty}\gamma$ varies through the $\Gamma_0(\stufe)$-orbit of
$\Gamma_{\infty}\gamma_0$, 
and
$$1|\pmatrix A&B\\C&D\endpmatrix(\tau)=\det(C\tau+D)^{-k}.$$
Since $k\ge 4$,
$\E_{\gamma_0}$ is absolutely convergent and analytic (in all variables of $\tau$);
for $\gamma'\in\Gamma_0(\stufe)$, $\Gamma_{\infty}\gamma\gamma'$ varies
over the $\Gamma_0(\stufe)$-orbit of $\Gamma_{\infty}\gamma_0$ as
$\Gamma_{\infty}\gamma$ does, and hence
$\E_{\gamma_0}|\gamma'=\chi(\det D_{\gamma'})\,\E_{\gamma_0}.$
As noted in [3], these Eisenstein series are linearly independent,
and the 0th Fourier coefficient of $\E_{\gamma_0}$ is 0 unless $\gamma_0\in\Gamma_0(\stufe)$,
in which case it is 1.

With $\gamma_0=\pmatrix *&*\\ M_0&N_0\endpmatrix$ and $\gamma=\pmatrix *&*\\ M&N\endpmatrix,$
we have $\Gamma_{\infty}\gamma_0\Gamma_0(\stufe)=\Gamma_{\infty}\gamma\Gamma_0(\stufe)$ if and
only if, for all primes $q|\stufe$, $\rank_qM_0=\rank_qM$, where
$\rank_q$ denotes rank modulo $q$.
Thus for each (multiplicative) 
partition $\rho=(\stufe_0,\stufe_1,\stufe_2)$  of $\stufe$ (so
$\stufe_0\stufe_1\stufe_2=\stufe$), take $M_{\rho}$ to be a diagonal matrix so that for each
prime $q|\stufe$, modulo $q$ we have
$$M_{\rho}\equiv\cases 0&\text{if $q|\stufe_0$,}\\
\pmatrix 1\\&0\endpmatrix&\text{if $q|\stufe_1$,}\\
I&\text{if $q|\stufe_2$.}
\endcases$$
Then set $\E_{\rho}=\E_{\gamma_{\rho}}$ where $\gamma_{\rho}=\pmatrix I&0\\ M_{\rho}&I\endpmatrix.$
Hence a natural basis for this space of Eisenstein series is
$$\{\E_{\rho}:\ \rho=(\stufe_0,\stufe_1,\stufe_2),\ \stufe_0\stufe_1\stufe_2=\stufe\ \}.$$

When $\gamma=\pmatrix *&*\\ M&N\endpmatrix\in Sp_2(\Z)$, $(M\ N)$ is a coprime symmetric
pair, meaning that $M,N$ are integral matrices with $M\ ^tN$ symmetric, and $(GM\ GN)$
is integral only if $G$ is.  (Here $^tN$ denotes the transpose of $N$.)
Conversely, given a coprime symmetric pair $(M\ N)$, there is some
$\gamma=\pmatrix *&*\\ M&N\endpmatrix\in Sp_2(\Z)$.  
Note that $(M\ N)$ is a coprime pair if and only if, for every prime $p$,
$\rank_p(M\ N)=2$.
It will be more convenient for us
to describe $\E_{\rho}$ in terms of coprime symmetric pairs.
We say $(M\ N)$ is $\stufe$-type $\rho$ if $(M\ N)$ is a coprime symmetric pair so that
$GL_2(\Z)(M\ N)$ is in the $\Gamma_0(\stufe)$-orbit of $GL_2(\Z)(M_{\rho}\ I)$, or equivalently,
$(M\ N)$ is a coprime symmetric pair so that for each prime $q|\stufe$,
$\rank_qM=\rank_qM_{\rho}$;
sometimes we simply say $M$ is $\stufe$-type $\rho$, and for $q$ prime, we simply say
$(M\ N)$ is $q$-type $i$ when $(M\ N)$ is a coprime symmetric pair with
$\rank_qM=i$.  So $\E_{\rho}$ can be defined as a sum over a set of 
representatives for the $GL_2(\Z)$-equivalence
classes of $\stufe$-type $\rho$.
Note that
$$GL_2(\Z)(M_{\rho}\ I)
=SL_2(\Z)(M_{\rho}\ I)\cup SL_2(\Z)(M_{\rho}\ I)\pmatrix -1\\&1\\&&-1\\&&&1\endpmatrix,$$
so we can consider $2\E_{\rho}$ to be supported on
a set of representatives for the
$SL_2(\Z)$-equivalence classes of $\stufe$-type $\rho$, and for $\gamma\in\Gamma_0(\stufe)$
so that $SL_2(\Z)(M\ N)=SL_2(\Z)(M_{\rho}\ I)\gamma$, we can describe $\chi(\det D_{\gamma})$
in terms of $M,N,\rho$ as follows.

Suppose $(M\ N)=(M_{\rho}\ I)\gamma$ where $\gamma=\pmatrix A&B\\C&D\endpmatrix\in\Gamma_0(\stufe)$
and $\rho=(\stufe_0,\stufe_1,\stufe_2)$ is a (multiplicative) partition of $\stufe$.
For each prime $q|\stufe_0$, we have $N\equiv D\ (q)$, so $\chi_q(\det D)=\chi_q(\det N)$.
For each prime $q|\stufe_2$, we have $M\equiv A\equiv\,^t\overline D\ (q)$,
so $\chi_q(\det D)=\overline\chi_q(\det M)$.  Now take a prime $q|\stufe_1$; write
$D=\pmatrix d_1&d_2\\d_3&d_4\endpmatrix$.  
Since $A\,^tD\equiv I\ (q)$, 
$$A\equiv\overline{\det D}
\pmatrix d_4&-d_3\\-d_2&d_1\endpmatrix,\ M\equiv\overline{\det D}
\pmatrix d_4&-d_3\\0&0\endpmatrix,\ N\equiv\pmatrix *&*\\d_3&d_4\endpmatrix\ (q).$$
We know $q\nmid(d_3\ d_4)$ and $\chi^2=1$, so $\chi_q(\det D)=\chi_q(m_1)\chi_q(n_4)$ or
$\chi_q(-m_2)\chi_q(n_3)$, whichever is non-zero.  
Take $E\in SL_2(\Z)$
so that $q$ divides row 2 of $EM$;
thus $E\equiv\pmatrix\alpha&\beta\\0&\overline\alpha\endpmatrix\ (q)$,
and with 
$M=\pmatrix m_1&m_2\\qm_3&qm_4\endpmatrix$,
$N=\pmatrix n_1&n_2\\n_3&n_4\endpmatrix$,
$EM=\pmatrix m_1'&m_2'\\qm_3'&qm_4'\endpmatrix$,
$EN=\pmatrix n_1'&n_2'\\n_3'&n_4'\endpmatrix$, we have
$\chi_q(m_1')\chi_q(n_4')=\chi_q(m_1)\chi_q(n_4)$ and
$\chi_q(-m_2')\chi_q(n_3')=\chi_q(-m_2)\chi_q(n_3)$.
So when $SL_2(\Z)(M\ N)=SL_2(\Z)(M_{\rho}\ I)\gamma$ and $q|\stufe_1$,
we can choose $E\in SL_2(\Z)$ so that $EM=\pmatrix m_1&m_2\\qm_3&qm_4\endpmatrix$,
$EN=\pmatrix n_1&n_2\\n_3&n_4\endpmatrix$; thus, setting
$\chi_{(1,q,1)}(M,N)=\chi_q(m_1)\chi_q(n_4)$ or
$\chi_q(-m_2)\chi_q(n_3)$ (whichever is non-zero), we have $\chi_q(\det D)=\chi_{(1,q,1)}(M,N)$.
We set
$$\chi_{\rho}(M,N)
=\prod_{q|\stufe_0}\chi_q(\det N)\prod_{q|\stufe_1}\chi_{(1,q,1)}(M,N)\prod_{q|\stufe_2}\overline\chi_q(\det M).$$
Hence $2\E_{\rho}(\tau)=\sum\overline\chi_{\rho}(M,N)\det(M\tau+N)^{-k}$ where 
$(M\ N)$ varies over a set of $SL_2(\Z)$-equivalence class representatives for pairs
of $\stufe$-type $\rho$.
Also note that for $G\in SL_2(\Z)$, $\chi_{\rho}(GM,GN)=\chi_{\rho}(M,N)$, and since
$\pmatrix G\\&^tG^{-1}\endpmatrix\in Sp_2(\Z)$, 
we have $\chi_{\rho}(MG,N\,^tG^{-1})=\chi_{\rho}(M,N)$.

We will often use the theory of quadratic forms over a finite field $\F=\Z/p\Z$ where 
$p$ is an odd prime.  For $V$ a vector space over $\F$ equipped with a quadratic
form $Q$ and associated bilinear form $B$, 
we say a non-zero vector $v$ is isotropic if $Q(v)=0$, and
we define the radical of $V$ to be
$$\rad V=\{v\in V: B(v,w)=0\text{ for all }w\in V\ \};$$
we say $V$ is regular if $\rad V$ is trivial.  When $V$ is regular with dimension
2, either $V\simeq\H$ or $\A$ where $\H\simeq\big<1,-1\big>$ is called a hyperbolic plane,
and with $\omega$ a non-square in $\F$, $\A\simeq\big<1,-\omega\big>$ is called an 
anisotropic plane.   Note we also have $\H\simeq\pmatrix 0&1\\1&0\endpmatrix.$
For $T,T'\in \F^{2,2}_{\sym}$, we write $T\simeq T'$ if there is some
$G\in GL_2(\F)$ so that $^tGTG=T'$.

\bigskip
\head{\bf \S3.  Defining Hecke, shift, and twist operators}\endhead
\smallskip

For each prime $p$, we have Hecke operators $T(p)$ and $T_1(p^2)$ that
act on degree 2 Siegel modular forms, and
$\{T(p),T_1(p^2):\ p\text{ prime }\}$ generates the Hecke algebra.
For $f$ a degree 2 Siegel modular form of weight $k$, level $\stufe'$, and
character $\chi'$, and for $\gamma'=\pmatrix A&B\\C&D\endpmatrix\in Sp_2(\Q)$, we set
$$f(\tau)|\gamma'=(\det\gamma')^{k/2}\det(C\tau+D)^{-k}f(\gamma'\circ\tau).$$  Then
$$f|T(p)=p^{k-3}\sum_{\gamma}\overline\chi'(\det D_{\gamma})f|\delta^{-1}\gamma$$
where $\delta=\pmatrix pI_2\\&I_2\endpmatrix$ and $\gamma$ varies over a set
of coset representatives for
$$\left(\delta\Gamma_0(\stufe')\delta^{-1}\cap\Gamma_0(\stufe')\right)\backslash\Gamma_0(\stufe').$$
Somewhat similarly,
$$f|T_1(p^2)=p^{k-3}\sum_{\gamma}\overline\chi'(\det D_{\gamma})f|\delta_1^{-1}\gamma$$
where $\delta_1=\pmatrix X\\&X^{-1}\endpmatrix$, $X=\pmatrix p\\&1\endpmatrix$, and $\gamma$
varies over a set of coset representatives for
$$\left(\delta_1\Gamma_0(\stufe')\delta_1^{-1}\cap\Gamma_0(\stufe')\right)
\backslash\Gamma_0(\stufe').$$
By Propositions 2.1, 3.1 and Theorem 6.1 of [5], for a prime $q|\stufe$ we know
$$f|T(q)=q^{k-3}\sum_Y f|\delta^{-1}\pmatrix I&Y\\&I\endpmatrix$$
where $Y\in\Z^{2,2}_{\sym}$ varies modulo $q$, and
$$f|T_1(q^2)=q^{k-3}\sum_{G,Y} f|\delta_1^{-1}\pmatrix G^{-1}&Y\,^tG\\&^tG\endpmatrix$$
where $G$ varies over
$$\G_1(q)=\left\{E\in SL_2(\Z):\ E\equiv\pmatrix *&0\\ *&*\endpmatrix\ (q)\ \right\}
\backslash SL_2(\Z)$$
and $Y=\pmatrix y_1&y_2\\y_2&0\endpmatrix$ with $y_1$ varying modulo $q^2$, $y_2$ varying
modulo $q$.

We define $B(q)$ by
$f|B(q)(\tau)=f(q\tau)$.  

We define $A_0^*(q)$ by
$$f|A^*_0(q)=\sum_{G,Y}f|\pmatrix ^tG^{-1}&YG/q\\&G\endpmatrix$$
where $Y=\pmatrix y\\&0\endpmatrix$ with $y$ varying modulo $q$,
and $G$ varies over $\G_1(q)$.
Normalising,
we define $A_0(q)=\frac{1}{q}A_0^*(q)$.  

We define $A^*_1(q)$ by
$$f|A^*_1(q)=\sum_{G,Y}\left(\frac{y_1}{q}\right)\,f|\pmatrix G^{-1}&Y\,^tG/q\\&^tG\endpmatrix$$
where $G$ varies as above and $Y=\pmatrix y_1&y_2\\y_2&y_4\endpmatrix$ varies over $\Z^{2,2}_{\sym}$ modulo $q$.
Normalising, we define $A_1(q)=\frac{1}{q^2\g(q)}A_1^*(q)$ where $\g(q)$ is the
standard Gauss sum, defined by
$$\g(q)=\sum_{a\,(q)}\e^{2\pi i a^2/q}.$$

We define $A^*_2(q)$ by
$$f|A^*_2(q)=\sum_{Y}\left(\frac{\det Y}{q}\right)\,f|\pmatrix I&Y/q\\&I\endpmatrix$$
where $Y$ varies over $\Z^{2,2}_{\sym}$ modulo $q$.  Normalising, we define
$A_2(q)=\frac{1}{q}\left(\frac{-1}{q}\right)A_2^*(q)$.

Say $E\in SL_2(\Z)$ so that $E\equiv\pmatrix *&0\\*&*\endpmatrix\ (q)$.  Then
$E^{-1}Y\,^tE^{-1}$ varies over $\Z^{2,2}_{\sym}$ modulo $q$ as $Y$ does, and 
$^tE\pmatrix \alpha&0\\0&0\endpmatrix E=\pmatrix \alpha'&0\\0&0\endpmatrix$
with $\alpha'$ varying modulo $q$ as $\alpha$ does.  From this
it is straight-forward to verify that the maps 
$A_i(q)$ are well-defined.  
Note that we may choose the matrices
$G$ so that, with $q^r||\stufe'$, we may choose
$$G\equiv\pmatrix 1&\alpha\\0&1\endpmatrix\text{ or }\pmatrix0&-1\\1&0\endpmatrix\ (q),$$
and $G\equiv I\ (\stufe'/q^r)$.  Similarly, we may choose $Y\equiv0\ (\stufe'/q^r)$.

Then we have the following.

\proclaim{Proposition 3.1}  Suppose $f$ is a degree 2 Siegel modular form of weight $k$,
level $\stufe'$ and character $\chi'$.  Then $f|B(q)$ is a degree 2 Siegel modular
form of weight $k$, level $q\stufe'$ and character $\chi'$; for $i=0,1,2$,
$f|A_i(q)$ is a degree 2 Siegel modular form of weight $k$, level
$\lcm(\stufe',q^2)$ and character $\chi'$.  With $c(T)$ denoting the $T$th coefficient
of $f$, the $T$th coefficient of $f|A_0(q)$ is
$$\cases 
2c(T)&\text{if $T\simeq\H\ (q)$,}\\
0&\text{if $T\simeq\A\ (q)$,}\\
qc(T)&\text{if $T\simeq\big<0,2\nu\big>\ (q)$, $\nu\not\equiv0\ (q)$,}\\
(q+1)c(T)&\text{if $T\simeq\big<0,0\big>\ (q).$}
\endcases$$
The $T$th coefficient of $f|A_1(q)$ is
$\left(\frac{\nu}{q}\right)c(T)$ if $T\simeq\big<2\nu,0\big>\ (q)$,
and 0 otherwise.
The $T$th coefficient of $f|A_2(q)$ is 
$$\cases qc(T)&\text{if $T\simeq\big<0,*\big>\ (q)$,}\\
-c(T)&\text{if $T\simeq\H$ or $\A\ (q)$.}\endcases$$
\endproclaim

\demo{Proof}  
It is easy to see $B(q)$ raises the level by $q$.  To see
that the $A_i(q)$ map forms to level $\lcm(\stufe',q^2)$, we first consider $i=1$.
Take $\gamma=\pmatrix A&B\\C&D\endpmatrix
\in \Gamma_0(\stufe')$ with $q^2|C$, $G\in\G_1(q)$, and $Y\in\Z^{2,2}_{\sym}$;
set $\alpha=\det A$.  Then
$\det \,^tDG\pmatrix 1\\&\alpha\endpmatrix\equiv1\ (q)$, so there is some
$E\in SL_2(\Z)$ so that $E\equiv\,^tDG\pmatrix 1\\&\alpha\endpmatrix\ (q).$
Set $U=\pmatrix 1\\&\overline\alpha\endpmatrix Y\pmatrix 1\\&\overline\alpha\endpmatrix$
where $\overline\alpha\in\Z$ so that $\alpha\overline\alpha\equiv1\ (q).$
Then $E$ varies over $SL_2(\Z)$ modulo $q$ as $G$ does; for fixed $G$ and $E$,
$U$ varies over $\Z^{2,2}_{\sym}$ modulo $q$ as $Y$ does, and
$\left(\frac{y_1}{q}\right)=\left(\frac{u_1}{q}\right)$.  Also, with
$$\gamma'=\pmatrix A'&B'\\C'&D'\endpmatrix
=\pmatrix G^{-1}&\frac{1}{q}Y\,^tG\\&^tG\endpmatrix\gamma
\pmatrix E&-\frac{1}{q}EU\\&^tE^{-1}\endpmatrix,$$
we have $\gamma'\in\Gamma_0(\stufe')$ with $q^2|C'$ and $\det D'\equiv\det D\ (q)$.
Therefore
$$\align
f|A'_1(q)|\gamma
&=\sum_{G,Y}\left(\frac{y_1}{q}\right)\,f|\pmatrix G^{-1}&\frac{1}{q}Y\,^tG\\&^tG\endpmatrix|\gamma\\
&=\sum_{E,U} \left(\frac{u_1}{q}\right)\,f|\gamma'|\pmatrix E^{-1}&\frac{1}{q}U\,^tE\\&^tE\endpmatrix\\
&=\chi'(\det D)\,f|A'_1(q).
\endalign$$
Similar arguments show that $f|A_i(q)|\gamma=\chi'(\det D)\,f|A_i(q)$ for $i=0,2$.

Let us now consider the effect of the $A_i(q)$ on Fourier coefficients.
Since we have $Tr(AB)=Tr(BA)$, we have
$$\align
f(\tau)|A^*_0(q)
&=\sum_{T,G,Y} c(T)\e\{T\,^tG^{-1}\tau G^{-1}+TY/q\}\\
&=\sum_{T,G,Y} c(T)\e\{G^{-1}T\,^tG^{-1}\tau\}\sum_{Y}\e\{TY/q\}.
\endalign$$
Since $G\in SL_2(\Z)$, we have $c(T)=c(GT\,^tG)$, so changing variables we get
$$f(\tau)|A^*_0(q)=\sum_T c(T)\e\{T\tau\}\sum_{G,Y}\e\{GT\,^tGY/q\}.$$
Let $V=\F x_1\oplus\F x_2$ be equipped with the quadratic form $Q$ given by $T$
relative to the basis $(x_1\ x_2)$.  Then
with $(x_1'\ x_2')=(x_1\ x_2)\,^tG$, $Q$ is given by $GT\,^tG$ relative to
$(x_1'\ x_2')$, and
$\F x_1'$ varies over all lines in $V$ as $G$ varies;
thus the sum on $Y$ tests whether $\F x_1'$ is isotropic (over $\F$).
So the sum on $G$ and $Y$ is $q$ times the number of isotropic lines in $V$; thus the sum is
$$\cases 2q&\text{if $T\simeq\H$,}\\
0&\text{if $T\simeq\A$,}\\
q^2&\text{if $T\simeq\big<0,2\nu\big>$, $\nu\not=0$,}\\
q(q+1)&\text{if $T\simeq\big<0,0\big>$.}\endcases$$

Similarly, we have
$$\align
f(\tau)|A^*_1(q)
&=\sum_{G,Y,T}\left(\frac{y_1}{q}\right)\,c(T)\,\e\left\{TG^{-1}\tau\,^tG^{-1}+TY\right\}\\
&=\sum_{G,T}c(T)\,\e\{T\tau\}\sum_Y\left(\frac{y_1}{q}\right)\,\e\{\,^tGTGY/q\}.
\endalign$$

Let $V=\F x_1\oplus\F x_2\simeq T$
(relative to $(x_1\ x_2)$).  For $G\in \G_1(q)$ and $(x'_1\ x'_2)=(x_1\ x_2)G$,
$Q$ is represented by $^tGTG\ (q)$ relative to the basis $(x'_1\ x'_2)$, and 
$\F x_2'$ varies over all lines in $V$ as $G$ varies.  Also,
$$\sum_Y\left(\frac{y_1}{q}\right)\e\left\{^tGTGY/q\right\}=0$$
unless $\F x_1'$ is an anisotropic line and
$\F x'_2=\rad V$ (which is only possible in the case that $\dim\rad V=1$).
So suppose $\dim\rad V=1$.  Then for 1 choice of $G\in\G_1(q)$ we have
we have $\F x_2'=\rad V$, so $^tGTG\equiv\big<2\nu,0\big>\ (q)$ with $q\nmid \nu$, and
the sum over $Y$ gives us $\left(\frac{\nu}{q}\right)\g(q)$.

Evaluating the action of $A_2(q)$ on Fourier coefficients is similar; we need to use the fact
that, summing over all $Y\in\F^{2,2}_{\sym}$,
$$\sum_Y\left(\frac{\det Y}{q}\right)\e\{TY/q\}
=\cases -q\g(q)^{-1}&\text{if $T\simeq\H$ or $\A\ (q)$,}\\
q\g(q)^{-1}(q-1)&\text{otherwise}\endcases$$
(see Theorem 1.3 of [14]).
$\square$
\enddemo

\bigskip
\head{\bf\S4. Building $T_K(q)$ so that
$\theta(\gen K)$ as an eigenform, $q$ a prime dividing $\stufe$}\endhead
\smallskip

Recall that $K$ is a maximal even integral lattice of 
rank $2k\ge 8$, square-free odd level $\stufe$
and associated character $\chi$ where $\chi^2=1$.
Throughout this section, we fix a prime $q$ 
dividing $\stufe$ and fix $\omega\in\Z$ so that $\left(\frac{\omega}{q}\right)=-1$.  
Set $\epsilon=\left(\frac{-1}{q}\right)$ and let $\F$ denote $\Z/q\Z$.
Take $R$ to be the preimage in $K$ of
$\rad K/qK$.
At the end of this section we will show that $\theta(\gen R)=\theta(\gen K)|T_{K/R}(q)$ and that
$\theta(\gen K)$ is an eigenform for
$T_K(q)$, where these maps are defined as follows.

\noindent{\bf Definitions.}  Set
$$T_{K/R}(q)=\left[T_1(q^2)T(q)-\gamma T(q)+q^{k-1}\frac{\lambda(q)}{\lambda(p)}T(p)\right]B(q)$$
where $\gamma=(q^{2k-3}+1)(q+1)+q^{k-1}$, $\lambda(t)=(t^{k-1}+1)(t^{k-2}+1)$, and the prime $p$ is
chosen so that $\chi(p)=1$ and for all primes $q'|\stufe/q$, $\chi_{q'}(p)=\chi_{q'}(q).$
(In Lemma 4.1 below we show such a prime $p$ exists.)
When $\chi_q=1$, set
$$\align
&T_K(q)=\\
&[T(q)^2+aT_1(q^2)-q^{k-2}A_0(q)T_1(q^2)+cA_0(q)+q^{2k-1}T(q)B(q)]\\
&\quad+T_{K/R}(q)[q^{2k-4}T_1(q^2)T_1(q^2)+q^{2k-3}aT_1(q^2)-q^{3k-5}T_1(q^2)A_0(q)+q^{4k-6}]
\endalign$$
where $a=\frac{q^2(q^{2k-6}-1)}{(q-1)}+(q^{k-2}-1)$,
$c=-\frac{q^{k+1}(q^{2k-6}-1)}{(q-1)}+q^{k-1}$.
When $\chi_q\not=1$, set
$$\align
&T_K(q)=\\
&\left[T(q)^2+aT_1(q^2)+\left(\frac{\eta}{q}\right)q^{k-2}A_1(q)T_1(q^2)
+\left(\frac{\eta}{q}\right)q^{k-1}a A_1(q)+q^{2k-4}A_2(q)\right]\\
&+ T_{K/R}(q)\left[q^{2k-4}T_1(q^2)T_1(q^2)+q^{2k-3}aT_1(q^2)+\left(\frac{\eta}{q}\right)q^{3k-5}T_1(q^2)A_1(q)
+q^{4k-6}\right]
\endalign$$
where $a=\frac{q(q^{2k-5}-1)}{(q-1)}-1$ and 
$\Z_qK\simeq2\big<1,-1,\ldots,1,-1,\eta\big>\perp 2q\big<\eta'\big>.$
\smallskip

Here we prove that there exists a prime $p$ as claimed in the above definitions.

\proclaim{Lemma 4.1}  Let $q$ be a prime dividing $\stufe$.
(1)  Say $\chi_q\not=1$; write
$$\Z_qK\simeq2\big<1,-1,\ldots,1,-1,\eta\big>\perp2q\big<\eta'\big>.$$
There exists a prime $p$ so that $\chi(p)=1$, $\chi_{q'}(p)=\chi_{q'}(q)$ for all primes
$q'|\Delta_0/q$, and
$\theta(\gen K)|T(p)=\lambda(p)\theta(\gen K^*)$
where $\lambda(p)=(p^{k-1}+1)(p^{k-2}+1)$,
$$\align
\Z_qK^*&\simeq\Z_qK^{\eta\eta'},\\
\Z_{q'}K^*&\simeq\Z_{q'}K^q \text{ for all primes }q'|\Delta_0/q,\\
\Z_{\ell}K^*&\simeq\Z_{\ell}K \text{ for all primes }\ell\nmid \Delta_0.
\endalign$$
(2) For $\chi_q=1$, there exists a prime $p$ so that
$\chi(p)=1$, $\chi_{q'}(p)=\chi_{q'}(q)$ for all primes $q'|\Delta_0$, and
 $\theta(\gen K)|T(p)=\lambda(p)\theta(\gen K^*)$
where $\lambda(p)=(p^{k-1}+1)(p^{k-2}+1)$,
$$\align
\Z_{q'}K^*&\simeq\Z_{q'}K^q \text{ for all primes }q'|\Delta_0,\\
\Z_{\ell}K^*&\simeq\Z_{\ell}K \text{ for all primes }\ell\nmid \Delta_0.
\endalign$$
\endproclaim

\demo{Proof}   Take $\Delta_0,\Delta_1\in\Z$ so that $\disc K=\Delta_0\Delta_1^2$
(and hence $\stufe=\Delta_0\Delta_1$).  So for $q|\stufe$, we have $\chi_q\not=1$ if and
only if $q|\Delta_0$.
By the argument used to prove Theorem 1.2 (1) of [17], when $\chi(p)=1$
we get $\theta(\gen K)|T(p)=\lambda(p)\theta(\gen K^*)$ where
$\Z_pK^*\simeq\Z_pK$ and $\Z_{\ell}K^*\simeq\Z_{\ell}K^p$ for all other primes $\ell$.
($K^p$ denotes the lattice $K$ whose quadratic form has been scaled by $p$.)
For any prime $\ell$, the Jordan components of $\Z_{\ell}K$ have even rank unless
$\ell|\Delta_0$, so for $\ell\nmid \Delta_0p$, we have $\Z_{\ell}K^p\simeq\Z_{\ell}K$.
Thus to prove the lemma, we need to show there exists a prime $p$ so that
$\chi(p)=1$, $\left(\frac{p}{q'}\right)=\left(\frac{q}{q'}\right)$
for all primes $q'\not=q$ with $q'|\Delta_0$, and $\left(\frac{p}{q}\right)
=\left(\frac{\eta\eta'}{q}\right)$ if $q|\Delta_0$.

(1) Suppose $q|\Delta_0$.  Using the Chinese Remainder Theorem and Dirichlet's Theorem,
we can choose a prime $p\nmid \stufe$ so that $\left(\frac{p}{q}\right)=\left(\frac{\eta\eta'}{q}\right)$
and $\left(\frac{p}{q'}\right)=\left(\frac{q}{q'}\right)$ for all primes $q'|\Delta_0/q$.
Note that since $\disc\Z_qK=(-1)^{k-1}\eta\eta' q$, we have
$$\eta\eta'\equiv(-1)^{k-1}\Delta_0/q\cdot u^2\ (q)$$ for some $u\not\equiv0\ (q).$
As discussed in \S2, $(-1)^k\disc K\equiv 1\ (4)$, so
$(-1)^{k-1}\Delta_0/q\equiv-q\ (4)$.
If $q\equiv1\ (4)$ then $\left(\frac{-1}{q}\right)=1$ and by Quadratic Reciprocity
$\left(\frac{\eta\eta'}{q}\right)=\left(\frac{q}{\Delta_0/q}\right)$.  If $q\equiv-1\ (4)$ then
$\Delta_0/q\equiv(-1)^{k-1}\ (4)$ and hence
$$\left(\frac{\eta\eta'}{q}\right)=(-1)^{k-1}\cdot(-1)^{(\Delta_0/q-1)/2}\left(\frac{q}{\Delta_0/q}\right)
=\left(\frac{q}{\Delta_0/q}\right).$$
In either case, 
$$\chi(p)=\left(\frac{p}{\Delta_0}\right)=\left(\frac{q}{\Delta_0/q}\right)\left(\frac{\eta\eta'}{q}\right)=1.$$

(2) Now suppose $q|\Delta_1$; so
$$\Z_qK\simeq2\big<1,-1,\ldots,1,-1,1,-\omega\big>\perp2q\big<1,-\omega\big>.$$
Choose a prime $p\nmid \stufe$ so that
$\left(\frac{p}{q'}\right)=\left(\frac{q}{q'}\right)$ for all primes $q'|\Delta_0$.  
We know $(-1)^kq^2\equiv \Delta_0u^2q^2\ (q^3)$ for some $u\not\equiv0\ (q)$, so
again using Quadratic Reciprocity and the fact that $\left(\frac{-1}{q}\right)=(-1)^{(q-1)/2}$, we have
$\left(\frac{(-1)^k\Delta_0}{q}\right)=1=\left(\frac{q}{\Delta_0}\right).$  Hence
$\chi(p)=\left(\frac{p}{\Delta_0}\right)=\left(\frac{q}{\Delta_0}\right)=1.$
$\square$
\enddemo

To prove $T_{K/R}(q)$ and $T_K(q)$ have the action on $\theta(\gen K)$ as claimed, we first
evaluate the action of the various maps used to construct $T_K(q)$. 
Note that $\Z_pR=\Z_pK$ for each prime $p\not=q$, and
$$\Z_qR\simeq 
\cases 2q\big<\eta'\big>\perp 2q^2\big<1,-1,\ldots,1,-1,\eta'\big>
&\text{if $\chi_q\not=1$,}\\
2q\big<1,-\omega\big>\perp2q^2\big<1,-1,\ldots,1,-1,1,-\omega\big>
&\text{if $\chi_q=1$.}\endcases$$
By Proposition 1.4 of [16], we know
that as $K'$ varies over the isometry classes in $\gen K$, $R'$ varies over the
isometry classes in $\gen R$ (where $R'$ denotes the preimage in $K'$ of $\rad K'/qK'$),
and $O(R')=O(K')$.

Let $\Omega$ be a sublattice of $\frac{1}{q}K$ with (formal) rank 2; 
throughout this section, we decompose $\Omega$ as
$(\frac{1}{q}\Omega_0\oplus\Omega_1)+\Omega_2$ where $\Omega_2\subset R$,
$\Omega_0\oplus\Omega_1$ is primitive in $K$ modulo $R$, meaning that
$\Omega_0,\Omega_1\subset K$,
and with $d_i(\Omega)=\rank\Omega_i$,
$\dim\overline{\Omega_0+\Omega_1}$ in $K/R$ is
$d_0(\Omega)+d_1(\Omega)$.  When it is clear we are referring to $\Omega,$ we simply
write $d_i$ rather than $d_i(\Omega)$.  Also, we will often write
$\overline\Omega$ to refer to $\Omega/q\Omega$.
Recall that 
$\G_1=\G_1(q)$
is a set of representatives for
$$\left\{G\in SL_2(\Z):\ G\equiv\pmatrix *&0\\ *&*\endpmatrix\ (q)\ \right\}\backslash SL_2(\Z),$$
and we can choose the elements of $\G_1$ to be congruent modulo $q$ to
$\pmatrix 1&\alpha\\0&1\endpmatrix$ or $\pmatrix 0&-1\\1&0\endpmatrix$, and congruent
modulo $\stufe/q$ to $I$.

\proclaim{Proposition 4.2}  We have
$$\theta(K;\tau)|T(q)=\sum_{{\Lambda\subset K}\atop{\overline\Lambda\simeq<0,0>}} \e\{\Lambda\tau/q\},$$
and  $\theta(R;\tau)|T(q)=\theta(R;\tau/q).$
\endproclaim

\demo{Proof} This follows immediately from Theorem 6.1 of [5]. $\square$
\enddemo

\proclaim{Proposition 4.3} For $\Omega\subset\frac{1}{q}K$,
write $d_i$ for $d_i(\Omega)$.  We have
$$\theta(K;\tau)|T_1(q^2)=\sum_{{\Omega\subset\frac{1}{q}K}\atop{d_0=1}}\e\{\Omega\tau\}
+(q+1)\theta(K;\tau),$$
and
$$\theta(R;\tau)|T_1(q^2)=\sum_{{\Omega\subset K}\atop{d_1=1}}\e\{\Omega\tau\}
+(q+1)\theta(R;\tau).$$
\endproclaim

\demo{Proof}  We know
$$f|T_1(q^2)=q^{k-3}\sum f|\pmatrix X^{-1}\\&X\endpmatrix \pmatrix G^{-1}&Y\,^tG\\&^tG\endpmatrix$$
where $X=\pmatrix q\\&1\endpmatrix$, $G$ varies over $\G_1$, and
$Y=\pmatrix y_1&y_2\\y_2\endpmatrix$ with $y_1$ varying modulo $q^2$, $y_1$ varying modulo $q$.
For a sublattice $\Lambda$ of $K$, let $T_{\Lambda}$ be a matrix for $Q$ restricted to $\Lambda$.
Thus with $G,Y$ varying as above and $E$ varying over $GL_2(\Z)$, we have
$$\align
\theta(K;\tau)|T_1(q^2)
&=q^{-3}\sum_{\Lambda\subset K}\sum_{E,G,Y}\e\{\,^tET_{\Lambda}EX^{-1}(G^{-1}\tau+Y\,^tG)\,^tG^{-1}X^{-1}\}\\
&=p^{-3}\sum_{{\Lambda\subset K}\atop{E,G}}\e\{\,^tG^{-1}X^{-1}\,^tET_{\Lambda}EX^{-1}G^{-1}\tau\}\\
&\qquad\cdot \sum_Y\e\{(\,^tET_{\Lambda}E)(X^{-1}YX^{-1})\}.
\endalign$$
The sum on $Y$ is $q^3$ if $X^{-1}\,^tET_{\Lambda}EX^{-1}$ is integral, and 0 otherwise.
When this sum is $q^3$, let $\Omega=\Lambda EX^{-1}.$  Notice that with
$\Omega'=\Lambda E'X^{-1}$, $E'\in GL_2(\Z)$, we have $\Omega=\Omega'$ if and only if
$E^{-1}E'\equiv\pmatrix *&0\\*&*\endpmatrix\ (q)$.  Hence
$$\align
\theta(K;\tau)|T_1(q^2)&=\sum_{{\Lambda\subset K}\atop{E\in GL_2(\Z)}}
\sum_{{\Omega\text{ integral}}\atop{\{\Omega:\Lambda\}=(1,q)}} \e\{\,^tET_{\Omega}E\tau\} \\
&=\sum_{{\Omega\subset\frac{1}{q}K}\atop{\Omega\text{ integral}}}
\left(\sum_{{\Lambda\subset K}\atop{\{\Omega:\Lambda\}=(1,q)}} 1 \right) \e\{\Omega\tau\}.
\endalign$$
Using our standard decomposition of $\Omega$, we see $d_0\le 1$ else there are no
$\Lambda\subset K$ so that $\{\Omega:\Lambda\}=(1,q)$, and when $d_0=1$
the only choice for $\Lambda$ is $\Omega\cap K$.  If 
$d_0=0$ then there are $q+1$ choices for $\Lambda$ (corresponding to
the lines in $\Omega/q\Omega$). 

Evaluating $\theta(R;\tau)|T_1(q^2)$ is similar, noting that for $x,y\in R$,
we have $\Z\frac{1}{q}x\oplus\Z y$ integral if and only if $x\in qK$.
$\square$
\enddemo

\proclaim{Proposition 4.4} Writing $d_i=d_i(\Omega)$ for $\Omega\subset\frac{1}{q} K$, we have
$$\theta(R;\tau)|T_1(q^2)T(q)
=\sum_{{\Omega\subset K}\atop{\overline\Omega_1\simeq<0>}} \e\{\Omega\tau/q\}+(q+1)\theta(R;\tau/q)$$
and
$$\align
\theta(R;\tau)|T_1(q^2)T_1(q^2)
&=\sum_{{\Omega\subset K}\atop{d_1=0}}\e\{\Omega\tau\}+(q+1)\sum_{{\Omega\subset K}\atop{d_1=2}}\e\{\Omega\tau\}\\
&+(2q+1)\sum_{{\Omega\subset K}\atop{d_1=1}}\e\{\Omega\tau\}+(q+1)^2\theta(R;\tau).
\endalign$$
\endproclaim

\demo{Proof} 
To see $\theta(R;\tau)|T_1(q^2)T(q)$ is as claimed, first apply Proposition 4.3 and then use the fact 
that $T(q)$ annihilates terms with 
$\overline\Omega\not\simeq\big<0,0,\big>$ and replaces $\tau$ by $\tau/q$.

Again using Proposition 4.3, we have
$$\align
\theta(R;\tau)|T_1(q^2)T_1(q^2)
&=\sum_{{\Omega\subset K}\atop{d_1=1}}\e\{\Omega\tau\}|T_1(q^2)
+(q+1)\sum_{{\Omega\subset K}\atop{d_1=1}}\e\{\Omega\tau\}\\
&+(q+1)^2\theta(R;\tau).
\endalign$$
With $m(\Omega)$ is the number of lattices $\Lambda\subset K$ so that
$d_1(\Lambda)=1$ and $\{\Omega:\Lambda\}=(1,q)$, we have
$$\sum_{d_0=0,\,d_1=1}\e\{\Omega\tau\}|T_1(q^2)=\sum_{\Omega\subset\frac{1}{q}K}m(\Omega)\e\{\Omega\tau\}.$$
Suppose $\Omega\not\subset K$; then to have $\Lambda\subset K$ so that $\{\Omega:\Lambda\}=(1,q)$, we must
have $\Lambda=\Omega\cap K$ and $d_0(\Omega)+d_1(\Omega)=1$.
Hence $m(\Omega)=1$ if $d_0(\Omega)=1$ and $d_1(\Omega)=0$, and otherwise $m(\Omega)=0$.

Now suppose $\Omega\subset K$.  If
$d_1(\Omega)=2$ then every line in $\overline\Omega$ corresponds to a
sublattice $\Lambda$ meeting the criteria.  If $d_1(\Omega)=1$ then
$q$ lines in $\overline \Omega$ correspond to sublattices $\Lambda$
meeting the criteria.  If $d_1(\Omega)=0$ then no sublattices
$\Lambda$ meet the criteria. $\square$
\enddemo

\proclaim{Proposition 4.5} With $\nu$ denoting any non-zero value modulo $q$, we have
$$\align
&\theta(K;\tau)|A_0(q)\\
&\quad = 
2\sum_{{\Omega\subset K}\atop{\overline\Omega\simeq\H}}\e\{\Omega\tau\}
+\sum_{{\Omega\subset K}\atop{\overline\Omega\simeq<0,2\nu>}}\e\{\Omega\tau\}
+(q+1)\sum_{{\Omega\subset L}\atop{\overline\Omega\simeq<0,0>}}\e\{\Omega\tau\}
\endalign$$
and
$$\align
&\theta(K;\tau)|A_0(q)T_1(q^2)\\
&\quad =
(2q+1)\theta(K;\tau)+\sum_{{d_0=1}\atop{\overline\Omega_1\simeq<2\nu>}}\e\{\Omega\tau\}
+(q+1)\sum_{{d_0=1}\atop{\overline{\Omega\cap K}\simeq<0,0>}}\e\{\Omega\tau\}\\
&\quad\quad +
q\sum_{{d_0=0}\atop{\overline\Omega\simeq\H}}\e\{\Omega\tau\}-q\sum_{{d_0=0}\atop{\overline\Omega\simeq\A}}
\e\{\Omega\tau\}+q^2\sum_{{d_0=0}\atop{\overline\Omega\simeq<0,0>}}\e\{\Omega\tau\}.
\endalign$$
\endproclaim

\demo{Proof} The first claim follows immediately from the definition of $A_0(q)$, since
$A_0(q)$ counts the number of isotropic lines in $\overline\Omega$.
So using the definition of $T_1(q^2)$,
$$\align
&\theta(K;\tau)|A_0(q)T_1(q^2)\\
&\quad =
2\sum_{{\Omega\subset K}\atop{\overline\Omega\simeq\H}}
\sum_{{\Lambda\text{ integral}}\atop{\{\Lambda:\Omega\}=(1,q)}} \e\{\Lambda\tau\}
+\sum_{{\Omega\subset K}\atop{\overline\Omega\simeq<0,2\nu>}}
\sum_{{\Lambda\text{ integral}}\atop{\{\Lambda:\Omega\}=(1,q)}} \e\{\Lambda\tau\}\\
&\qquad
+(q+1)\sum_{{\Omega\subset K}\atop{\overline\Omega\simeq<0,0>}}
\sum_{{\Lambda\text{ integral}}\atop{\{\Lambda:\Omega\}=(1,q)}} \e\{\Lambda\tau\}.
\endalign$$
Changing the order of summation, we first sum over integral $\Lambda\subset\frac{1}{q}K$,
then over $\Omega\subset\Lambda\cap K$ with $\{\Lambda:\Omega\}=(1,q)$.  Now take integral
$\Lambda\subset\frac{1}{q}K$; then for any $\Omega$ so that $\{\Lambda:\Omega\}=(1,q)$,
$\overline\Omega$ will necessarily have a radical of dimension at least 1.  Hence 
$\overline\Omega\simeq\big<0,*\big>$.

Suppose $\Lambda\not\subset K$; to have any $\Omega\subset\Lambda\cap K$ so that
$\{\Lambda:\Omega\}=(1,q)$, we must have $d_0(\Lambda)=1$, $\Omega=\Lambda\cap L$, and
$\overline\Omega\simeq\big<0,*\big>$.
With $\Omega=\Lambda\cap L$, we have
 $\overline\Omega\simeq\big<0,2\nu\big>$ if and
only if $\overline\Lambda_1\simeq\big<2\nu\big>$; similarly, $\overline\Omega\simeq\big<0,0\big>$
if and only if $\overline\Lambda\simeq\big<0\big>$ or $d_1(\Lambda)=0$.
(Recall that here $\nu$ represents any element of $\F^{\times}$.)

Now suppose $\Lambda\subset K$. Then each $\Omega\subset\Lambda\cap K$ with $\{\Lambda:\Omega\}=(1,q)$
corresponds to a line in $\Lambda/q\Lambda$.  Each isotropic line in $\Lambda/q\Lambda$ 
corresponds to $\Omega$ where
$\overline\Omega\simeq\big<0,0\big>$, and each anisotropic line in
$\Lambda/q\Lambda$ corresponds to $\Omega$ where
$\Omega\simeq\big<0,2\nu\big>$.
We know $\Lambda/q\Lambda$ has $q+1$ lines; also, $\H$ has 2 isotropic lines, $\A$ has no
isotropic lines, $\big<0,2\nu\big>$ has 1, and $\big<0,0\big>$ has $q+1$.
From this the second claim  follows. $\square$
\enddemo

\proclaim{Proposition 4.6}  We have
$$\theta(K;\tau)|A_1(q)=\sum_{{\Omega\subset K}\atop{\overline\Omega\simeq<0,2\nu>}}\left(\frac{\nu}{q}\right)
\e\{\Omega\tau\},\ \theta(R;\tau)|A_1(q)=0,$$
and
$$\align
&\theta(K;\tau)|A_1(q)T_1(q^2)\\
&\quad = \sum_{{d_0=1}\atop{\overline\Omega_1\simeq<2\nu>}}\left(\frac{\nu}{q}\right)\e\{\Omega\tau\}
+ q\sum_{{d_0=0}\atop{\overline\Omega\simeq<0,2\nu>}}\left(\frac{\nu}{q}\right)\e\{\Omega\tau\}.
\endalign$$
\endproclaim

\demo{Proof}  From the definition, we have
$$\align
\theta(K;\tau)|A_1(q)
&=\frac{1}{q^2\g(q)}\sum_{x_1,x_2\in K}\sum_G\e\{(B(x_i,x_j))\tau\}\\
&\qquad\cdot\sum_Y\left(\frac{y_1}{q}\right)\e\left\{\frac{1}{q}\,^tG(B(x_i,x_j))GY\right\} 
\endalign$$
where $G$ varies over $\G_1$, 
$Y=\pmatrix y_1&y_2\\y_2&y_4\endpmatrix\in\Z^{2,2}_{\sym}$ varying modulo $q$.
Thus the sum on $Y$ is $\left(\frac{\nu}{q}\right)q^2\g(q)$ if 
$^tG(B(x_i,x_j))G\equiv\pmatrix \nu&0\\0&0\endpmatrix\ (q)$, and 0 otherwise.
Let $V=\F x_1\oplus\F x_2$, $(x_1'\ x_2')=(x_1\ x_2)G$.  Then as $G$ varies,
$\F x_2'$ varies over all lines in $V$.  If $\rad V=\{0\}$ or $V$, then the sum
on $Y$ is 0 for each choice of $G$; otherwise, there is a unique $G$ so that
$\F x_2'=\rad V$, and in this case the sum on $Y$ is $\left(\frac{\nu}{q}\right)q^2\g(q)$
where $V\simeq\big<\nu,0\big>.$  

Since $A_1(q)$ annihilates terms $c(\Lambda)\e\{\Lambda\tau\}$ where 
$\overline\Lambda\simeq\big<0,0\big>$, $\theta(R;\tau)|A_1(q)=0$.

We also have
$$
\theta(K;\tau)|A_1(q)T_1(q^2)
=\sum_{{\Omega\subset K}\atop{\overline\Omega\simeq<0,2\nu>}}\left(\frac{\nu}{q}\right)
\sum_{{\Lambda\text{ integral}}\atop{\{\Lambda:\Omega\}=(1,q)}}\e\{\Lambda\tau\}.
$$
Interchanging the order of summation, we sum first over integral $\Lambda\subset\frac{1}{q}K$,
then over $\Omega\subset\Lambda\cap K$ so that $\{\Lambda:\Omega\}=(1,q)$ and
$\overline\Omega\simeq\big<0,2\nu\big>$.  

Suppose first $\Lambda\not\subset K$.
Then $\Lambda=\frac{1}{q}\Lambda_0\oplus\Lambda_1\oplus\Lambda_2$
with $d_0(\Lambda)>0$, and there are no $\Omega$ 
meeting the criteria if $d_0(\Lambda)>1$.  So suppose $d_0(\Lambda)=1$; then we need to take
$\Omega=\Lambda\cap K$ to have $\{\Lambda:\Omega\}=(1,q)$, and then $\overline\Omega\simeq\big<0,0\big>$
unless $\overline\Lambda_1\simeq\big<2\nu\big>$ (and hence $d_1(\Lambda)=1$).

Now suppose $\Lambda\subset K$.  Then each $\Omega$ where $\{\Lambda:\Omega\}=(1,q)$
corresponds to a line in $\Lambda/q\Lambda$, and any line that represents a quadratic
residue cannot represent a quadratic non-residue, and vice-versa.  If $\Lambda/q\Lambda\simeq\H$ or $\A$ or
$\big<0,0\big>$ then the number of lines of $\Lambda/q\Lambda$ that represent quadratic
residues is equal to the number of lines that represent quadratic non-residues.  If
$\Lambda/q\Lambda\simeq\big<0,2\nu'\big>$ with $q\nmid\nu'$, then $\Lambda/q\Lambda$ has
1 istotropic line, and $q$ lines $\F x\simeq\big<2\nu'\big>$.  From this the proposition follows.
$\square$
\enddemo

\proclaim{Proposition 4.7}  We have
$$\theta(K;\tau)|A_2(q)=-\theta(K;\tau)+q\sum_{{\Omega\subset K}\atop{\overline\Omega\simeq<0,*>}}
\e\{\Omega\tau\},$$
and
$$\theta(R;\tau)|A_2(q)=q\theta(R;\tau).$$
\endproclaim

\demo{Proof}  We have
$$\theta(K;\tau)|A_2(q)=\frac{1}{\epsilon q}\sum_{x_1,x_2\in K}
\sum_Y\left(\frac{\det Y}{q}\right) \e\{(B(x_i,x_j))(\tau+Y/q)\},$$
$Y$ varying over symmetric matrices modulo $q$.  For $x_1,x_2\in K$,
$$\sum_Y\left(\frac{\det Y}{q}\right)\e\{(B(x_i,x_j))Y/q\}
=\cases -\epsilon q&\text{if $(B(x_i,x_j))\simeq\H$ or $\A\ (q)$,}\\
\epsilon(q-1)q&\text{otherwise,}\endcases$$
from which the proposition follows. $\square$
\enddemo

\proclaim{Proposition 4.8}  With $q$ a prime dividing $\stufe$, there is a prime $p$ so that
$\chi(p)=1$, $\chi_{\stufe/q}(p)=\chi_{\stufe/q}(q)$, and 
$$\theta(\gen R)|[q^{k-2}T_1(q^2)T(q)B(q)+q^{2k-3}]
=\theta(\gen K)|\left[\frac{\lambda(q)}{\lambda(p)}\cdot T(p)B(q)-T(q)B(q)\right]$$
where $\lambda(t)=(t^{k-1}+1)(t^{k-2}+1).$
\endproclaim

\demo{Proof} First suppose $\chi_q\not=1$.  By \S91C and 92:2 [13],
there is a basis $(x_1,\ldots,x_{2k})$ for $\Z_qK$ so that, relative to this basis,
$\Z_qK\simeq2\big<1,-1,\ldots,1,-1,\eta\big>\perp2q\big<\eta'\big>.$
(So $R=\Z_q x_{2k}+\Z_q qK.$)
Let $\overline C$ be a maximal totally isotropic subspace of $K/R\simeq\Z_qK/\Z_qR$;
let $M$ be the preimage in $K$ of $\overline C$.  
Since $SL_{2k}(\Z)$ projects onto $SL_{2k}(\F)$, there is a basis
$(x_1',\ldots,x_{2k-1}',x_{2k})$ of $\Z_qK$ so that 
$(\overline x_1',\ldots,\overline x_{k-1}')$ is a basis for $\overline C$.
(So $(x_1',\ldots,x_{k-1}',qx_k',\ldots,qx_{2k-1}',x_{2k})$ is a basis for $\Z_qM$.)
By 92:2 [13], we can 
adjust $x_k',\ldots,x_{2k-1}'$ so that
$$\Z_q x_1'\oplus\cdots\oplus\Z_q x_{2k-2}'\simeq 2\pmatrix 0&I\\I&0\endpmatrix\ (q),$$
and hence $\Z_q x_1'\oplus\cdots\oplus\Z_q x_{2k-2}'$
 is an orthogonal sum of hyperbolic planes.
So $$\Z_qx_1'\oplus\cdots\oplus\Z_qx_{k-1}'\oplus\Z_q qx_k'\oplus\cdots\oplus\Z_q qx_{2k-2}'
\simeq 2q\pmatrix *&I\\I&0\endpmatrix\ (q^2),$$
which means $\Z_qx_1'\oplus\cdots\oplus\Z_qx_{k-1}'\oplus\Z_q qx_k'\oplus\cdots\oplus\Z_q qx_{2k-2}'$
is an orthogonal sum of hyperbolic planes scaled by $q$.
Again appealing to 92:2 [13], we have
$$\Z_qM\simeq 2q\big<1,-1,\ldots,1,-1\big>\perp 2q^2\big<\eta\big>\perp 2q\big<\eta'\big>.$$

Take $\Omega\subset K$ with
formal rank 2; decompose $\Omega$ as $\Omega_1+\Omega_2$ where
$\Omega_2\subset R$ and $d_1=\rank\Omega_1=\dim\overline\Omega_1$ in $K/R$.
So $\Omega\subset M$ if and only if $\overline\Omega_1\subseteq\overline C$.
Also, $\Omega$ can only be in $M$ if $\overline\Omega_1$ is totally isotropic
in $K/R$.  So suppose $\overline\Omega_1$ is totally isotropic; then the number
of $M$ containing $\Omega$ is the number of ways to extend $\overline\Omega_1$
to a maximal totally isotropic subspace $\overline C$ of $K/R$; hence the number
of $M$ containing $\Omega$ is
$$\cases \delta=\prod_{i=3}^{k-1}(q^{k-i}+1)&\text{if $d_1=2$,}\\
(q^{k-2}+1)\delta&\text{if $d_1=1$,}\\
(q^{k-1}+1)(q^{k-2}+1)\delta&\text{if $d_1=0$.}
\endcases$$

We know that for any prime $\ell\not=q$, $\Z_{\ell}M=\Z_{\ell}K.$
Using Propositions 4.2 and 4.3, we see
$$\delta^{-1}\sum_M\theta(M;\tau/q)
=\theta(K;\tau)|T(q)+\theta(R;\tau)|[q^{k-2}T_1(q^2)T(q)+q^{2k-3}T(q)].$$
Further, with $M^{1/q}$ denoting the lattice $M$
with its quadratic form scaled by $1/q$, we have $\Z_qM^{1/q}\simeq\Z_qK$ and
$\Z_{\ell}M^{1/q}\simeq\Z_{\ell}K^{1/q}\simeq\Z_{\ell}K^q$ for all primes $\ell\not=q$.
Sorting the $M$ into isometry classes (and fixing a choice of $M$ to represent the genus),
we have
$$\sum_M\theta(M;\tau/q)=\sum_{\cls M'\in\gen M}\frac{\#\{\sigma\in O(\Q K):\ qK\subset\sigma M'\subset K\ \}}
{o(M')} \theta(M';\tau/q).$$
Note that for $M'\in\gen M$, $\Q K=\Q M=\Q M'$, $o(qK)=o(K)$, and
$$\#\{\sigma\in O(\Q K):\ qK\subset\sigma M'\subset K\ \}=
\#\{\sigma\in O(\Q M'):\ qM'\subset\sigma qK\subset M'\ \}.$$
Thus averaging over the genus of $K$ and using that $O(K)=O(R)$
(see Proposition 1.4 [16]), we get
$$\align
&\theta(\gen K;\tau)|T(q)+\theta(\gen R;\tau)|[q^{k-2}T_1(q^2)T(q)+q^{2k-3}T(q)]\\
&\quad = \delta^{-1}\sum_{{K'\in\gen K}\atop{M'\in\gen M}}
\frac{\#\{\sigma\in O(\Q M'):\ qM'\subset\sigma qK'\subset M'\ \}}{o(qK')o(M')}\theta(M';\tau/q)\\
&\quad = (q^{k-1}+1)(q^{k-2}+1)\theta(\gen M;\tau/q)
\endalign$$
since $$\sum_{K'\in \gen K}\frac{\#\{\sigma\in O(\Q K):\ qK\subset\sigma M'\subset K\ \}}{o(qK')}$$
is the number of maximal totally isotropic subspaces of $M'/qM'$ scaled by $1/q$,
which is $(q^{k-1}+1)(q^{k-2}+1)\delta$.

Let $K^*=M^{1/q}$; thus $\theta(\gen M;\tau/q)=\theta(\gen K^*;\tau)$
and for $\ell$ a prime,
$$\Z_{\ell}K^*\simeq
\cases \Z_{\ell}K \text{ if $\ell\nmid \Delta_0/q$,}\\
\Z_{\ell}K^q \text{if $\ell|\Delta_0/q$.}\endcases$$  
By Lemma 4.1, we can choose a prime $p$ so that
$\theta(\gen K;\tau)|T(p)=\lambda(p)\theta(\gen K^*;\tau).$  Substituting for 
$\theta(\gen M;\tau/q)$, applying $B(q)$, and noting that
$\theta(R;\tau)|T(q)B(q)=\theta(R;\tau)$ (see Proposition 4.2), we get the result.

In the case $\chi_q=1$, $\Z_qK\simeq 2\big<1,-1,\ldots,1,-1,1,-\omega\big>\perp2q\big<1,-\omega\big>$;
otherwise, the argument is virtually identical to the argument above. $\square$
\enddemo

\proclaim{Proposition 4.9}
Let $q$ be a prime dividing $\stufe$.
\item{(1)} Suppose $\chi_q\not=1$; 
with $\alpha=(q+1)(q^{2k-3}+1)$, we have
$$\theta(\gen R)|q^{2k-3}T_1(q^2)
=\theta(\gen K)|[\alpha-T_1(q^2)-q^{k-1}\left(\frac{\eta}{q}\right)A_1(q)]$$
where $\Z_qK\simeq\big<1,-1,\ldots,1,-1,\eta\big>\perp q\big<\eta'\big>.$
\item{(2)} Suppose $\chi_q=1$.  With $\alpha=(q+1)(q^{2k-3}-q^{k-1}+1)$, we have
$$\theta(\gen R)|q^{2k-3}T_1(q^2)
=\theta(\gen K)|[\alpha-T_1(q^2)+q^{k-1}A_0(q)].$$
\endproclaim

\demo{Proof}  Suppose $\chi_q\not=1$.  
Take a basis $(x_1,\ldots,x_{2k})$ for $\Z_qK$ so that, relative to this basis,
$$\Z_qK\simeq2\big<1,-1,\ldots,1,-1,\eta\big>\perp 2q\big<\eta'\big>.$$
Let $\overline C$ be an isotropic line in $K/R\simeq\Z_qK/\Z_qR$.  
Since $K/R$ is regular, there is a line
$\overline D$ in $K/R$ so that $\overline C\oplus\overline D\simeq\H$.  Since $SL_{2k}(\Z)$
projects onto $SL_{2k}(\Z/q\Z)$ , there is a basis
$(x_1',\ldots,x_{2k-1}',x_{2k})$ for $\Z_qK$ so that, in $K/R$, $\overline C=\F\overline x_1'$
and $\overline D=\F\overline x_2'$.
By 92:2 [13], $\Z_qx_1'\oplus\Z_qx_2'\simeq\H$ and
by 82:16 [13], $\Z_qx_1'\oplus\Z_qx_2'$ splits $\Z_qK$.  Again using 92:2 [13],
we can assume that, relative to the basis $(x_1',\ldots,x_{2k-1}',x_{2k})$,
$$\Z_qK\simeq 2\pmatrix qa&b\\b&c\endpmatrix\perp2\big<1,-1,\ldots,1,-1,\eta\big>\perp 2q\big<\eta'\big>$$
where $q\nmid b$.
So $(x_1',qx_2',\ldots,qx_{2k-1}',x_{2k})$ is a basis for $\Z_qK'$, and relative to this basis,
$$\Z_qK'\simeq 2q\pmatrix a&b\\b&qc\endpmatrix\perp2q^2\big<1,-1,\ldots,1,-1,\eta\big>\perp 2q\big<\eta'\big>.$$

Now scale $K/qK'$ by $1/q$ and take $R'$ to be the preimage in $K'$ of $\rad K'/qK'$
(so $(qx_1',q^2x_2',qx_3',\ldots,qx_{2k-1}',qx_{2k})$ is a basis for $\Z_qR'$).
In $K'/R'$ scaled by $1/q$, let $\overline C'$ be an isotropic line so that
$\overline C'\not=\F\overline{qx_2'}$.  
So $\overline C'$ is generated by some $\overline x_1''$ where $x_1''=x_1'+\alpha qx_2',$
$\alpha\in\Z$.
Take $qK_1$ to be the preimage in $K'$ of
$\overline C'$; thus $(x_1''/q,qx_2',x_3',\ldots,x_{2k-1}',x_{2k})$ is a basis for $\Z_qK_1$
and relative to this basis,
$$\Z_qK_1\simeq 2\pmatrix a'&b'\\b'&q^2c\endpmatrix\perp2\big<1,-1,\ldots,1,-1,\eta\big> \perp 2q\big<\eta'\big>,$$
$q\nmid b'$.
Since $2\pmatrix a'&b'\\b'&q^2c\endpmatrix\simeq\H$
by 92:2 [13], we have $\Z_qK_1\simeq\Z_qK$.
We also have $\Z_{\ell}K_1=\Z_{\ell}K$ for every prime $\ell\not=q$,
so $K_1\in\gen K$.

Now take $\Omega$ in $\frac{1}{q}K$ with formal rank 2; we count how many of these lattices
$K_1$ contain $\Omega$.  Decompose
$\Omega$ as $\left(\frac{1}{q}\Omega_0\oplus\Omega_1\right)+\Omega_2$ so that $\Omega_2\subset R$
and $\Omega_0\oplus\Omega_1$ is primitive in $K$ modulo $R$;
let $d_i=d_i(\Omega)$.
Then $\Omega\subset K_1$ if and only if $\overline\Omega_0\subseteq\overline C\subseteq
\overline\Omega_1^{\perp}$ in $K/R$, and $\overline\Omega_0\subseteq\overline C'$
in $K'/R'$.  So when $d_0=1$, there is 1 $K_1$ containing $\Omega$.  When $d_0=0$, the
number of $\overline C\subseteq\overline\Omega_1^{\perp}$ is
$$\cases
(q^{2k-4}-1)/(q-1)&\text{if $\overline\Omega_1\simeq\H$ or $\A$ or $\big<0,0\big>$,}\\
(q^{2k-4}-1)/(q-1)+\left(\frac{\eta\nu}{q}\right)q^{k-2}&\text{if $\overline\Omega_1\simeq\big<0,2\nu\big>$, $\nu\in\F^{\times}$,}\\
(q^{2k-3}-1)/(q-1)&\text{if $\overline\Omega_1\simeq\big<0\big>$,}\\
(q^{2k-3}-1)/(q-1)+\left(\frac{\eta\nu}{q}\right)q^{k-2}
&\text{if $\overline\Omega_1\simeq\big<2\nu\big>$, $\nu\in\F$,}\\
(q^{2k-2}-1)/(q-1)&\text{if $d_1=0$.}\endcases$$
Then the number of choices for $\overline C'$ is $q$ (which is the
number of isotropic lines in $\F\overline x\oplus\F\overline{qy}$ scaled by $1/q$
that are independent of $\F\overline{qy}$).  Thus, using Propositions 4.3 and 4.6, with
$c=\frac{q^2(q^{2k-5}-1)}{(q-1)}-1$ 
$$\theta(K;\tau)|[T_1(q^2)+\left(\frac{\eta}{q}\right)q^{k-1}A_1(q)+c]
+\theta(R;\tau)|q^{2k-3}T_1(q^2)
=\sum_{K_1}\theta(K_1;\tau).$$
So averaging over $\gen K$ gives us
$$\align
&\theta(\gen K;\tau)|[T_1(q^2)+\left(\frac{\eta}{q}\right)q^{k-1}A_1(q)+c]
+\theta(\gen R;\tau)|q^{2k-3}T_1(q^2)\\
&\quad = 
\sum_{M,M_1\in\gen K}\frac{1}{o(M)}\frac{\#\{\sigma\in O(\Q K):\ \{M:\sigma M_1\}
=(1/q,1,\ldots,1,q)\ \}}{o(M_1)} \theta(M_1;\tau)\\
&\quad =
\sum_{M,M_1\in\gen K}\frac{1}{o(M_1)}\frac{\#\{\sigma\in O(\Q K):\ \{M_1:\sigma M\}
=(1/q,1,\ldots,1,q)\ \}}{o(M)}\theta(M_1;\tau)\\
&\quad =
\frac{q(q^{2k-2}-1)}{(q-1)}\theta(\gen K;\tau).
\endalign$$
Rearranging this equation yields the result for $\chi_q\not=1$.

Now suppose $\chi_q=1$.  Let $\overline C$ be an isotropic line in $K/R$,
$K'$ the preimage in $K$ of $\overline C$.  So there is
a basis $(x,y,w_1,w_2,z_1,\ldots,z_{2k-4})$ for $\Z_qK$ so that
$(x,qy,w_1,w_2,qz_1,\ldots,qz_{2k-4})$ is a basis for $\Z_qK'$, and relative to this basis
the quadratic form on $K'$ is 
$$2q\pmatrix a&b\\b&qc\endpmatrix\perp 2q\big<1,-\omega\big>
\perp2q^2\big<1,-1,\ldots,1,-1,1,-\omega\big>$$
with $b\not\equiv 0\ (q)$.  Now in $K'/qK'$ scaled by $1/q$, let 
$R'$ be the preimage in $K'$ of $\rad K'/qK'$.  Then in $K'/R'$ scaled by $1/q$,
let $\overline C'$ be an
isotropic line independent of $\overline {qK}$; let $qK_1$
be the preimage in $K'$ of $\overline C'\oplus\rad K'/qK'$.  
As when $\chi_q\not=1$, $K_1\in\gen K$ with
$\{K:K_1\}=(1/q,1,\ldots,1,q)$.

For $\Omega\subset\frac{1}{q}K$ (decomposed as before), there is 1 $K_1$ containing
$\Omega$ if $d_0=1$.  If $d_0=0$ the number of $\overline C\subset\overline\Omega_1^{\perp}$ is
$$\cases
(q^{2k-5}-1)/(q-1)-q^{k-3}&\text{if $\overline\Omega_1\simeq\H$,}\\
(q^{2k-5}-1)/(q-1)+q^{k-3}&\text{if $\overline\Omega_1\simeq\A$,}\\
(q^{2k-5}-1)/(q-1)&\text{if $\overline\Omega_1\simeq\big<0,2\nu\big>$, $\nu\in\F^{\times}$,}\\
(q^{2k-5}-1)/(q-1)-q^{k-2}&\text{if $\overline\Omega_1\simeq\big<0,0\big>$,}\\
(q^{2k-4}-1)/(q-1)&\text{if $\overline\Omega_1\simeq\big<2\nu\big>$, $\nu\in\F^{\times}$,}\\
(q^{2k-4}-1)/(q-1)-q^{k-2}&\text{if $\overline\Omega_1\simeq\big<0\big>$,}\\
(q^{2k-3}-1)/(q-1)-q^{k-2}&\text{if $d_1=0$.}\endcases$$
The number of choices for $\overline C'$ is
$$\frac{(q^2+1)(q-1)-(q-1)}{(q-1)}=q^2.$$
Using Propositions 4.3 and 4.5, then averaging over $\gen K$ yields the result for $\chi_q=1$.
$\square$
\enddemo

\proclaim{Theorem 4.10} With $\gamma=\frac{q^{2k-3}+1)}{(q-1)}+q^{k-1}$,
$\lambda(t)=(t^{k-1}+1)(t^{k-2}+1)$ and a prime $p$ chosen as in Lemma 4.1,
we have
$$q^{3k-4}\theta(\gen R)
=\theta(\gen K)\Big|\left[T_1(q^2)T(q)-\gamma T(q)
+q^{k-1}\frac{\lambda(q)}{\lambda(p)}T(p)\right]B(q).$$
\endproclaim

\demo{Proof}  We apply $T(q)B(q)$ to the identity of the preceeding proposition;
note that $A_1(q)T(q)=0$ since $A_1(q)$ annihilates all coefficients attached
to lattices $\Lambda\simeq\big<0,0\big>\ (q)$, and $T(q)$ annihilates
all coefficients attached to lattices $\Lambda\not\simeq\big<0,0\big>\ (q)$.
Similarly, $A_0(q)T(q)=(q+1)T(q)$.  Using this together with Proposition 4.8
yields the result. $\square$
\enddemo

\proclaim{Theorem 4.11} With $K$ maximal of odd level $\stufe$, $q$ a prime dividing
$\stufe$, and $T_K(q)$ as defined at the beginning of this section,
we have $$\theta(\gen K)|T_K(q)=\kappa(q)\theta(\gen K)$$ for
some $\kappa(q)\in\Q$.
\endproclaim

\demo{Proof}  The argument here is much the same as that of Proposition 4.9:  We 
first count how often $\Omega\subset\frac{1}{q}K$ lies in lattices $K_2$ where
$K_2\in\gen K$ with $\{K:K_2\}=(1/q,1/q,1,\ldots,1,q,q)$.  Then we realise
$\sum_{K_2}\theta(K_2)$ in terms of images of $\theta(K)$ and $\theta(R)$; averaging
over $\gen K$ yields our desired result.

Whether $\chi_q=1$ or not, we begin by taking $\overline C$ to be a dimension 2
totally isotropic subspace of $K/R$ so that $\overline C$;
let $K'$ be the preimage in $K$ of $\overline C$.
Let $R'$ be the preimage in $K'$ of $\rad K'/qK'$, where $K'/qK'$ is scaled by $1/q$.
Then in $K'/R'$ scaled by $1/q$,
take $\overline C'$ to be a dimension 2 totally isotropic subspace
that is independent of $\overline{qK}$; let $qK_2$ be the preimage
in $K'$ of $\overline C'$.  So $\{K:K_2\}=(1/q,1/q,1,\ldots,1,q,q)$, and an
argument virtually identical to that used in the proof of Proposition 4.9 shows that
$K_2\in\gen K$.  Take $\Omega\subset\frac{1}{q}K$ with $\Omega
=(\frac{1}{q}\Omega_0\oplus\Omega_1)+\Omega_2$, $\Omega_2\subset R$, $\Omega_0\oplus\Omega_1$
primitive in $K$ modulo $R$; let $d_i=\dim \Omega_i$.  As in the proof of Proposition 8.3, we have $\Omega\in K_2$
if and only if $\overline\Omega_0\subseteq\overline C\subseteq\overline\Omega_1^{\perp}$ in
$K/R$ and $\overline\Omega_0\subseteq\overline C'$ in $K'/R'$.

To count the number of $K_2$ containing $\Omega\subset\frac{1}{q}K$, first suppose
$\chi_q\not=1$.
If $d_0=2$ then the number of $K_2$ containing $\Omega$ is 1.  Say $d_0=1$; then the number
of $\overline C$ so that $\overline\Omega_0\subseteq\overline C\subseteq\overline\Omega_1^{\perp}$ is
$$\cases (q^{2k-5}-1)/(q-1)+\left(\frac{\eta\nu}{q}\right)q^{k-2}
&\text{if $\overline\Omega\simeq\big<2\nu\big>$, $\nu\in\F$,}\\
(q^{2k-4}-1)/(q-1)&\text{if $d_1=0$.}\endcases$$
Then the number of ways to choose $\overline C'$ so that $\Omega\in K_2$ is $q$, which is
the number of ways to extend $\overline\Omega_0$ to a dimension 2 totally isotropic subspace
of $K'/qK'$, independent $\overline{qK}$.  Now say $d_0=0$; then
the number 
of $\overline C$ so that $\overline\Omega_0\subseteq\overline C\subseteq\overline\Omega_1^{\perp}$ is
$$\cases
\frac{(q^{2k-4}-1)}{(q-1)}&\text{if $\overline\Omega_1\simeq\H$ or $\A$,}\\
\frac{(q^{2k-4}-1)}{(q-1)}+q^{2k-6}+\left(\frac{\eta\nu}{q}\right)q^{k-3}\frac{(q^{2k-5}-1)}{(q-1)}
&\text{if $\overline\Omega_1\simeq\big<2\nu,0\big>$, $\nu\in\F$,}\\
\frac{(q^{2k-4}-1)^2}{(q^2-1)(q-1)} + \left(\frac{\eta\nu}{q}\right)\frac{q^{k-3}(q^{2k-4}-1)}{(q-1)}
&\text{if $\overline\Omega_1\simeq\big<2\nu\big>$,}\\
\frac{(q^{2k-2}-1)(q^{2k-4}-1)}{(q^2-1)(q-1)}&\text{if $d_1=0.$}\endcases$$
The number of ways to choose $\overline C'$ so that $\Omega\subset K_2$ is $q^3$.

This means that 
$$\align
&\sum_{K_2}\theta(K_2)\\
&= \theta(K)|\left[T(q)^2+aT_1(q^2)+\left(\frac{\eta}{q}\right)q^{k-2}A_1(q)T_1(q^2)\right]\\
&\quad + \theta(K)|\left[\left(\frac{\eta}{q}\right)q^{k-1}aA_1(q)+q^{2k-4}A_2(q)+e\right]\\
&\quad + \theta(R)|\left[ q^{2k-4}T_1(q^2)T_1(q^2)+q^{2k-3}aT_1(q^2)+\left(\frac{\eta}{q}\right)
q^{3k-5}T_1(q^2)A_1(q)+q^{4k-6}\right]
\endalign$$
where $e=\frac{q^3(q^{2k-4}-1)(q^{2k-6}-1)}{(q^2-1)(q-1)}-a(q+1)-(q^{2k-3}+1).$  Averaging over
$\gen K$ gives us the result in the case $\chi_q\not=1$.

Now suppose $\chi_q=1$; let $\Omega$ be as above.  
If $d_0=2$ then the number of $K_2$ containing $\Omega$ is 1.  Say $d_0=1$; then the number
of $\overline C$ so that $\overline\Omega_0\subseteq\overline C\subseteq\overline\Omega_1^{\perp}$ is
$$\cases
(q^{2k-6}-1)/(q-1)&\text{if $\overline\Omega_1\simeq\big<2\nu\big>$, $\nu\in\F^{\times}$,}\\
(q^{2k-6}-1)/(q-1)-q^{k-3}&\text{if $\overline\Omega_1\simeq\big<0\big>$,}\\
(q^{k-2}+1)(q^{k-3}-1)/(q-1)&\text{if $d_1=0$.}\endcases$$
Then the number of ways to choose $\overline C'$ 
so that $\Omega\subset K_2$ is $q^2$.  Say $d_0=0$; then the number 
of $\overline C$ so that $\overline\Omega_0\subseteq\overline C\subseteq\overline\Omega_1^{\perp}$ is
$$\cases
\frac{(q^{2k-6}-1)}{(q-1)}\left(\frac{(q^{2k-1}-1)}{(q-1)}-q^{k-4}\right)
&\text{if $\overline\Omega_1\simeq\H$,}\\
\frac{(q^{2k-6}-1)}{(q-1)}\left(\frac{(q^{2k-1}-1)}{(q-1)}+q^{k-4}\right)
&\text{if $\overline\Omega_1\simeq\A$,}\\
\frac{(q^{2k-6}-1)^2}{(q^2-1)(q-1)}
&\text{if $\overline\Omega_1\simeq\big<2\nu,0\big>$, $\nu\in\F^{\times}$,}\\
\frac{(q^{2k-6}-1)^2}{(q^2-1)(q-1)}-\frac{q^{k-3}(q^{2k-6}-1)}{(q-1)}+q^{2k-6}
&\text{if $\overline\Omega_1\simeq\big<0,0\big>$,}\\
\frac{(q^{2k-4}-1)(q^{2k-6}-1)}{(q^2-1)(q-1)}
&\text{if $\overline\Omega_1\simeq\big<2\nu\big>$, $\nu\in\F^{\times}$,}\\
\frac{(q^{2k-4}-1)(q^{2k-6}-1)}{(q^2-1)(q-1)}+\frac{q^{k-3}(q^{k-2}+1)(q^{k-3}-1)}{(q-1)}
&\text{if $\overline\Omega_1\simeq\big<0\big>$,}\\
\frac{(q^{k-1}+1)(q^{2k-4}-1)(q^{k-3}-1)}{(q^2-1)(q-1)}&\text{if $d_1=0$.}\endcases$$
Then the number of ways to choose $\overline C'$ 
so that $\Omega\subset K_2$ is $q^5$.

Thus
$$\align
&\sum_{K_2}\theta(K_2)\\
&=\theta(K)|[T(q)^2+aT_1(q^2)-q^{k-2}A_0(q)T_1(q^2)+cA_0(q)+q^{2k-1}T(q)B(q)+e]\\
&\quad+\theta(R)|[q^{2k-4}T_1(q^2)T_1(q^2)+q^{2k-3}aT_1(q^2)-q^{3k-5}T_1(q^2)A_0(q)+q^{4k-6}]
\endalign$$
where $e=\frac{q^5(q^{2k-6}-1)^2}{(q^2-1)(q-1)}-1-a(q+1)-c-q^{2k-4}(q+1).$  Averaging over
$\gen K$ yields the result in the case $\chi_q=1$.

The proof in the case that $\chi_q\not=1$ is virtually identical, using the following
(where $\Omega\subset\frac{1}{q}K$, as above).
If $d_0=2$ then the number of $K_2$ containing $\Omega$ is 1.  Say $d_0=1$; then the
number of $\overline C$ so that $\overline\Omega_0\subseteq\overline C\subseteq\Omega_1^{\perp}$
is
$$\cases
(q^{2k-5}-1)/(q-1)+\left(\frac{\eta\nu}{q}\right)q^{k-3}
&\text{if $\overline\Omega_1\simeq\big<2\nu\big>$, $\nu\in\F^{\times}$,}\\
(q^{2k-5}-1)/(q-1)&\text{if $\overline\Omega_1\simeq\big<0\big>$,}\\
(q^{2k-4}-1)/(q-1)&\text{if $d_1=0$,}\endcases$$
and then the number of ways to choose $\overline C'\oplus\rad K'/qK'$ so that
$\Omega\subset K_2$ is $q$.  If $d_0=0$, then the number of $\overline C$ so that
$\overline C\subset\overline\Omega_1^{\perp}$ is
$$\cases
\frac{(q^{2k-4}-1)(q^{2k-6}-1)}{(q^2-1)(q-1)}&\text{if $\overline\Omega_1\simeq\H$ or $\A$,}\\
\frac{(q^{2k-4}-1)(q^{2k-6}-1)}{(q^2-1)(q-1)}+q^{2k-6}
+\left(\frac{\eta\nu}{q}\right)\frac{q^{k-3}(q^{2k-5}-1)}{(q-1)}
&\text{if $\overline\Omega_1\simeq\big<0,2\nu\big>$, $\nu\in\F^{\times}$,}\\
\frac{(q^{2k-4}-1)(q^{2k-6}-1)}{(q^2-1)(q-1)}+q^{2k-6}
&\text{if $\overline\Omega_1\simeq\big<0,0\big>$,}\\
\frac{(q^{2k-4}-1)^2}{(q^2-1)(q-1)}+\left(\frac{\eta\nu}{q}\right)\frac{q^{k-3}(q^{2k-4}-1)}{(q-1)}
&\text{if $\overline\Omega_1\simeq\big<2\nu\big>$, $\nu\in\F^{\times}$,}\\
\frac{(q^{2k-4}-1)^2}{(q^2-1)(q-1)}&\text{if $\overline\Omega_1\simeq\big<0\big>$,}\\
\frac{(q^{2k-2}-1)(q^{2k-4}-1)}{(q^2-1)(q-1)}&\text{if $d_1=0$.}\endcases$$
Then the number of $\overline C'\oplus K'/qK'$ so that $\Omega\subset K_2$ is $q^3$.
From this the theorem follows for $\chi_q\not=1$.
$\square$
\enddemo

\bigskip
\head{\bf \S5. Action of shift, Hecke, and twist operators on Eisenstein series of level $q\stufe$,
where $\stufe$ is square-free and $q$ is a prime dividing $\stufe$}\endhead
\smallskip

Fix a prime $q$ dividing $\stufe$.
Our next task is to evaluate the action of $T_K(q)$ on the Eisenstein series
of level $\stufe$; since the shift and twist operators associated to $q$ map
Eisenstein series of level $\stufe$ to Eisenstein series of level $q\stufe$,
we need to evaluate the action of the Hecke and twist operators
on Eisenstein series of level $q\stufe$.

Fix a partition $\sigma=(\stufe_0,\stufe_1,\stufe_2)$
of $\stufe/q$; take $M_{\sigma}$ as defined in \S2.
Fix $\omega\in\Z$ so that
$\left(\frac{\omega}{q}\right)=-1$; let $\F$ denote $\Z/q\Z$ and
$\G_1$ denote $\G_1(q)$ (as defined in \S2).
Let $\psi(*)=\left(\frac{*}{q}\right)$;
using notation introduced in \S2, we abbreviate
$\psi_{(1,q,1)}(M,N)$ to $\psi_1(M,N)$.  Similarly, we abbreviate
$\chi_{(1,q,1)}(M,N)$ to $\chi_1(M,N)$.

We can decompose $\E_{\sigma}$ as $\E_{\sigma}=\E_0+\E_1+\E_2$ where the pairs
in the support of $\E_i$ have $q$-type $i$.  So $2\E_1$ is supported on a set of 
$SL_2(\Z)$-equivalence
class representatives for the $\Gamma_0(\stufe)$-orbit of $SL_2(\Z)(M_1\ I)$ where $M_1$
is a diagonal matrix with $M_1\equiv M_{\sigma}\ (\stufe/q)$, $M_1\equiv\pmatrix 1\\&0\endpmatrix\ (q)$;
this 
$\Gamma_0(\stufe)$-orbit splits futher into the $\Gamma_0(q\stufe)$-orbits corresponding to 
diagonal matrices $M_{1,1}$,
$M_{1,2,+}$ and $M_{1,2,-}$ where each is congruent modulo $\stufe/q$ to $M_{\sigma}$, and
$M_{1,1}\equiv\pmatrix 1\\&0\endpmatrix\ (q^2)$, $M_{1,2,+}\equiv\pmatrix 1\\&q\endpmatrix\ (q^2)$,
$M_{1,2,-}\equiv\pmatrix 1\\&q\omega\endpmatrix\ (q^2).$  Similarly, $2\E_0$ is supported on
a set of $SL_2(\Z)$-equivalence class representatives for 
the $\Gamma_0(\stufe)$-orbit of
$SL_2(\Z)(M_0\ I)$ where $M_0\equiv M_{\sigma}\ (\stufe/q)$, $M_0\equiv0\ (q)$, and this 
$\Gamma_0(\stufe)$-orbit splits into
$\Gamma_0(q\stufe)$-orbits corresponding to diagonal matrices
$M_{0,0}$, $M_{0,1,+}$, $M_{0,1,-}$, $M_{0,2,\alpha}$
where $\alpha$ varies over $\F^{\times}$; these matrices are each congruent modulo $\stufe/q$
to $M_{\sigma}$, and $M_{0,0}\equiv 0\ (q^2)$, $M_{0,1,+}\equiv\pmatrix q\\&0\endpmatrix\ (q^2)$,
$M_{0,1,-}\equiv\pmatrix q\omega\\&0\endpmatrix\ (q^2)$, $M_{0,2,\alpha}\equiv\pmatrix q\\&q\alpha\endpmatrix
\ (q^2)$.  (The $\Gamma_0(\stufe)$-orbit of $SL_2(I\ I)$ is that same as the $\Gamma_0(q\stufe)$-orbit
of $SL_2(\Z)(I\ I).$)  
With $M_*$ one of these matrices and $SL_2(\Z)(M\ N)$ in the $\Gamma_0(q\stufe)$-orbit
of $SL_2(\Z)(M_*\ I)$, we say $(M\ N)$ is $q^2$-type $*$.  When $(M\ N)$ is $q^2$-type
$i,j,\nu$ with $\nu=+$ or $-$,
we sometimes simply say $(M\ N)$ is $q^2$-type $i,j$.

We write $2\E_{1,1}$ to denote the level $q\stufe$,
character $\chi$ Eisenstein series supported on the $\Gamma_0(q\stufe)$-orbit of 
$SL_2(\Z)(M_{1,1}\ I)$,
and so forth.  We let $\E_{1,2}=\E_{1,2,+}+\E_{1,2,-}$, $\E_{0,1}=\E_{0,1,+}+\E_{0,1,-}$,
$\E_{0,2,+}=\sum\E_{0,2,\alpha}$ where the sum is over all $\alpha\in\F^{\times}$ with $\left(\frac{\alpha}{q}\right)=1$,
$\E_{0,2,-}=\sum\E_{0,2,\alpha}$ where the sum is over all $\alpha\in\F^{\times}$ with $\left(\frac{\alpha}{q}\right)=-1$,
and $\E_{0,2}=\E_{0,2,+}+\E_{0,2,-}$.  
With $SL_2(\Z)(M\ N)=SL_2(\Z)(M_*\ I)\gamma$ for $\gamma\in\Gamma_0(q\stufe)$,
as we saw in \S2 we have
$$\chi(\det D_{\gamma})=\chi_{\sigma}(M,N)\cdot
\cases \chi_q(\det N)&\text{if $*=0,0$ or $0,1,\nu$ or $0,2,\nu$,}\\
\chi_1(M,N)&\text{if $*=1,1$ or $1,2,\nu$,}\\
\chi_q(\det M)&\text{if $*=2$.}
\endcases$$
Recall that for $E\in SL_2(\Z)$, $\chi_{\sigma}(EM,EN)=\chi_{\sigma}(M,N)$,
and $\chi_1(EM,EN)=\chi_1(M,N)$.
For the rest of this section, we only consider pairs $(M\ N), (M'\ N')$ of
$\stufe/q$-type $\sigma$.

When evaluating the action of operators on these Eisenstein series,
we will frequently use the following.

\proclaim{Proposition 5.1} Suppose $\gamma=\pmatrix A&B\\C&D\endpmatrix\in\Gamma_0(\stufe^2)$ and
$(M\ I)=(M_{*}\ I)\gamma$ where $M_*$ is one of the above defined matrices.
\item{(a)}  If $M_*=M_{1,2,\nu}$ with $\nu=+$ or $-$, then $\nu\left(\frac{(\det M)/q}{q}\right)
=\psi_1(M,N).$
\item{(b)}  If $M_*=M_{0,1,\nu}$ with $\nu=+$ or $-$, then $\nu\left(\frac{\det N}{q}\right)
=\psi_1(M/q,N).$
\item{(c)} If $M_*=M_{0,2,\alpha}$ with $\alpha\in\F^{\times}$, then $\det\frac{1}{q}MN\equiv\alpha\ (q)$.
\endproclaim

\demo{Proof}  Recall that with $\gamma$ as above, $q^2|C$, so $\det AD\equiv1\ (q^2).$

(a) Say $M_*=M_{1,2,\nu}$.  Then $M\equiv\pmatrix 1\\&q\beta\endpmatrix A\ (q^2)$,
where $\beta=1$ or $\omega$ with $\nu\left(\frac{\beta}{q}\right)=1$.
So $\nu\left(\frac{\det M/q}{q}\right)=\left(\frac{\det A}{q}\right)=\left(\frac{\det D}{q}\right)$
and by the discussion in \S2, $\left(\frac{\det D}{q}\right)=\psi_1(M,N)$.

(b) Say $M_*=M_{0,1,\nu}$.  Then 
$M\equiv\pmatrix q\beta\\&0\endpmatrix A\ (q^2)$ where $\beta=1$ or $\omega$ so that
$\nu\left(\frac{\beta}{q}\right)=1$, and
$N\equiv D\ (q)$.  Thus, following the analysis used in \S2, we find
$\psi_1(M/q,N)=\nu\left(\frac{\det D}{q}\right)=\nu\left(\frac{\det N}{q}\right).$

(c) Say $M_*=M_{0,2,\alpha}$.  Thus 
$M\equiv q\pmatrix 1\\&\alpha\endpmatrix A\ (q^2)$ and $N\equiv D\ (q)$.
Thus $\det\frac{1}{q}MN\equiv\alpha\cdot\det AD\equiv\alpha\ (q).$
$\square$
\enddemo

The following proposition is easily verified, and thus the proof is omitted.

\proclaim{Proposition 5.2} For $\chi_q\not=1$, we have
$$\align
\E_0|B(q)&=\chi_{\stufe_1}(q)\E_{0,0},\\
\E_1|B(q)&=\chi_{\stufe_0\stufe_2}(q)q^{-k}\E_{1,1}+\chi_{\stufe_1}(q)(\E_{0,1,+}-\E_{0,1,-}),\\
\E_2|B(q)&=\chi_{\stufe_1}(q)q^{-2k}\E_2+\chi_{\stufe_0\stufe_2}(q)q^{-k}(\E_{1,2,+}-\E_{1,2,-})\\
&\quad +\chi_{\stufe_1}(q)(\E_{0,2,+}-\E_{0,2,-}).
\endalign$$
For $\chi_q=1$, we have
$$\align
\E_0|B(q)&=\chi_{\stufe_1}(q)\E_{0,0},\\
\E_1|B(q)&=\chi_{\stufe_0\stufe_2}(q)q^{-k}\E_{1,1}+\chi_{\stufe_1}(q)\E_{0,1},\\
\E_2|B(q)&=\chi_{\stufe_1}(q)q^{-2k}\E_2+\chi_{\stufe_0\stufe_2}(q)q^{-k}\E_{1,2}+\chi_{\stufe_1}(q)\E_{0,2}.
\endalign$$
\endproclaim

In [18] we evaluated the action of $T(q)$, $T_1(q^2)$ on Eisenstein series
of level $\stufe$.  Here we need to consider level $q\stufe$, but the techniques are
very similar.  

Proposition 3.5 of [18] states the following.

\proclaim{Proposition 5.3} $\E_2|T(q)=\chi_{\stufe_1}(q)q^{2k-3}\E_2$.
\endproclaim

\proclaim{Proposition 5.4}  $\E_{1,1}|T(q)=\chi_{\stufe_0\stufe_1}(q) q^{k-1}\E_1$.
\endproclaim

\demo{Proof}  
By Proposition 3.1 [5], we have
$$\align
&2\E_{1,1}(\tau)|T(q)\\
&\quad = q^{-3}\sum\chi_{\sigma}\chi_1(M,N)\det\left(\frac{1}{q}M\tau+N+\frac{1}{q}MY\right)^{-k}\\
&\quad = q^{2k-3}\sum\chi_{\sigma}\chi_1(M,N)\det(M\tau+qN+MY)^{-k}
\endalign$$
where the sums are over $SL_2(\Z)$-equivalence class representatives $(M\ N)$ for pairs
of $\stufe/q$-type $\sigma$ and $q^2$-type $1,1$.

Take $(M\ N)$ of $q^2$-type $1,1$.
We can assume $q^2$ divides row 2 of $M$.  Set
$$(M'\ N')=\pmatrix 1\\&1/q\endpmatrix(M\ qN+MY).$$
Thus $\rank_qM'=1$ and $(M',N')=1$ since $q$ does not divide row 1 of $M'$ or
row 2 of $N'$.  Also, $\chi_{\sigma}(M',N')=\chi_{\stufe_0\stufe_2}(q)\chi_{\sigma}(M,N)$.

Reversing, take $(M'\ N')$ of $q$-type 1; assume $q$ divides row 2 of $M'$.
We need to count the equivalence classes $SL_2(\Z)(M\ N)$ so that
$\pmatrix 1\\&1/q\endpmatrix(M\ qN+MY)\in SL_2(\Z)(M'\ N')$.  For any
$E\in SL_2(\Z)$, we have $\pmatrix 1\\&q\endpmatrix E\pmatrix 1\\&1/q\endpmatrix\in SL_2(\Z)$
if and only if $E\equiv \pmatrix *&0\\ *&*\endpmatrix \ (q)$; thus we need to count the integral,
coprime pairs 
$$(M\ N)=\pmatrix 1\\&q\endpmatrix E\left(M'\ \frac{1}{q}(N'-M'Y)\right)$$
where $E$ varies over $\G_1$.
For $q$ choices of $E$ we have $q$ dividing row 2 of $EM'$, and for
1 choice of $E$ we have $q$ dividing row 1 of $EM'$.  When $q$ divides row 1 of $EM'$,
we have $q$ dividing $M$ so $M$ is not $q^2$-type $1,1$.
When $q$ divides row 2 of $EM'$, $N$ is integral with $(M,N)=1$ for $q$ choices of $Y$.

Since $\det(M\tau+qN+MY)^{-k}=q^{-k}\cdot\det(M'\tau+N')^{-k}$, and each pair $(M'\ N')$
corresponds to $q^2$ pairs $(M\ N)$, this shows
$2\E_{1,1}|T(q)=\chi_{\stufe_0\stufe_2}(q)q^{k-1}\cdot 2\E_1$.
$\square$
\enddemo

The above proof demonstrates our method of argument throughout this section:
We begin with an Eisenstein series written in terms of $SL_2(\Z)$-equivalence
class representatives $(M\ N)$,
apply an operator which alters this pair, then construct a coprime pair $(M'\ N')$.
We then reverse this construction to count how many times the pairs $(M\ N)$ lead to
pairs $(M'\ N')$ in the same $SL_2(\Z)$-equivalence class.  Hence whenever our construction
of $(M'\ N')$ involves left-multiplication by something other than an element of $\Q$,
we need to consider whether left-multiplying our representatives for $SL_2(\Z)(M'\ N')$
by $E\in SL_2(\Z)$ changes the equivalence class of the resulting $(M\ N)$.
When our construction of $(M'\ N')$ involves left-multiplication by
$\pmatrix 1\\&1/q\endpmatrix$ or $\pmatrix q\\&1\endpmatrix$, our reverse construction
of $(M\ N)$ involves left-multiplication by $\pmatrix 1\\&q\endpmatrix E$ or
$\pmatrix 1/q\\&1\endpmatrix E$ where $E\in\G_1$; when our construction of $(M'\ N')$
involves left-multiplication by $\pmatrix 1\\&q\endpmatrix$ or $\pmatrix 1/q\\&1\endpmatrix$,
our reverse construction involves left-multiplication by $\pmatrix 1\\&1/q\endpmatrix E$
or $\pmatrix q\\&1\endpmatrix E$ where $E\in\,^t\G_1$.

\proclaim{Proposition 5.5} For $\nu=+$ or $-$, 
$$
\E_{1,2,\nu}|T(q)=
\cases \nu \chi_{\stufe_0\stufe_2}(q) q^{k-3}(q^2-1)/2\cdot\E_2 &\text{if $\chi_q\not=1$,}\\
\chi_{\stufe_0\stufe_2}(q) q^{k-3}(q^2-1)/2\cdot\E_2 &\text{if $\chi_q=1$.}
\endcases$$
\endproclaim

\demo{Proof}  Take $(M\ N)$ of $q^2$-type $1,2,\nu$.
We can assume $q$ divides row 2 of $M$.
Set $$(M'\ N')=\pmatrix 1\\&1/q\endpmatrix (M\ qN+MY).$$
Then $\rank_qM'=2$ so $(M',N')=1$.  Also,
$\chi_{\sigma}(M',N')=\chi_{\stufe_0\stufe_2}(q)\chi_{\sigma}(M,N)$.

Reversing, take $(M'\ N')$ of $q$-type 2.  Set
$$(M\ N)=\pmatrix 1\\&q\endpmatrix  E\left(M'\ \frac{1}{q}(N'-M'Y)\right).$$
For each $E$, there are $q-1$ choices for $Y$ so that $N$ is integral and
$(M,N)=1$.
To have $(M\ N)$ of $q^2$-type $1,2,\nu$, we need $\psi_1(M,N)=\nu\psi\left(\frac{1}{q}\det M\right).$
Write $EM'=\pmatrix m_1&m_2\\m_3&m_4\endpmatrix$, $EN'=\pmatrix n_1&n_2\\n_3&n_4\endpmatrix.$
We need to choose $Y$ so that $(n_1\ n_2)\equiv (m_1\ m_2)Y\ (q)$.  If $q\nmid m_1$ we can
choose any $y_4$, then solve for $y_2, y_1$.  If $q\nmid m_2$ we can choose any $y_1$, then
solve for $y_2, y_4$.  In either case, 1 solution is $Y\equiv \overline M'N'\ (q)$, which
gives $N\equiv 0\ (q)$ and hence $(M,N)\not=1$.
So when $q\nmid m_1m_4$ there is 1 choice of $y_4$ so that $n_4\equiv m_3y_2+m_4y_4\ (q);$
similarly, when $q\nmid m_2m_3$ there is 1 choice of $y_1$ so that $n_3\equiv m_3y_1+m_4y_2\ (q).$
Consequently $\psi_1(M,N)=\nu\psi\left(\frac{1}{q}\det M\right)$ for $(q-1)/2$ choices
of $Y$.  Also, $\psi(\det M')=\psi\left(\frac{1}{q}\det M\right)$, yielding the result.
$\square$
\enddemo

\proclaim{Proposition 5.6} $\E_{0,0}|T(q)=\chi_{\stufe_1}(q)\E_0$.
\endproclaim

\demo{Proof}  Take $(M\ N)$ of $q^2$-type $0,0$. Set
$(M'\ N')=\frac{1}{q}(M\ qN+MY).$  Thus $\rank_qM'=0$,
$\rank_qN'=\rank_qN=2$.  Also, $\chi_{\sigma}(M',N')=\chi_{\stufe_1}(q)\chi_{\sigma}(M,N)$.

Reversing, $\chi(\det N)=\chi(\det N')$ for all choices of $Y$, yielding the result.
$\square$
\enddemo

\proclaim{Proposition 5.7} Take $\nu=+$ or $-$. Then
$$\E_{0,1,\nu}|T(q)=
\cases \nu \chi_{\stufe_1}(q) q^{-1}(q-1)/2\cdot\E_1 &\text{if $\chi_q\not=1$,}\\
\chi_{\stufe_1}(q) q^{-1}(q-1)/2\cdot\E_1 &\text{if $\chi_q=1$.}
\endcases $$
\endproclaim

\demo{Proof} Take $(M\ N)$ of $q^2$-type $0,1,\nu$; can assume $q$ divides row 2 of $M$.
Set $$(M'\ N')=\frac{1}{q}(M\ qN+MY).$$
So $\rank_q(M',N')=2$, and $(M',N')=1$.
Also, $\chi_{\sigma}(M',N')=\chi_{\stufe_1}(q)\chi_{\sigma}(M,N)$.

Reversing, take $(M'\ N')$ of $q$-type 1; 
we find there are $q^2(q-1)/2$ choices for 
$Y$ so that $\nu\psi(\det N)=\psi_1\left(\frac{1}{q}M,N\right).$
$\square$
\enddemo

\proclaim{Proposition 5.8} For $\nu=+$ or $-$, we have
$$\E_{0,2,\nu}|T(q)=
\cases \nu \chi_{\stufe_1}(q) q^{-2}(q-1)(q+\nu\epsilon)/2\cdot\E_2 &\text{if $\chi_q\not=1$,}\\
\chi_{\stufe_1}(q) q^{-2}(q-1)(q+\nu\epsilon)/2\cdot\E_2 &\text{if $\chi_q=1$.}
\endcases$$
\endproclaim

\demo{Proof}  Take $(M\ N)$ of $q^2$-type $0,2,\nu$.  Set
$$(M'\ N')=\frac{1}{q}(M\ qN+MY).$$
Thus $\rank_qM'=2$ so $(M',N')=1$, and $\chi_{\sigma}(M',N')=\chi_{\stufe_1}(q)\chi_{\sigma}(M,N)$.

Reversing, take $(M'\ N')$ of $q$-type 2. 
We need to choose $Y$ so that $\nu\psi(\det N)=\psi(\det M'),$
or equivalently, we need 
$$\nu\psi(\det(N'\,^tM'-M'Y\,^tM'))=1.$$
As $Y$ varies over $\F^{2,2}_{\sym}$, so does $U=N'\,^tM'-M'Y\,^tM'$.  The number of
$U\in\F^{2,2}_{\sym}$ so that $\psi(\det U)=1$ is $q(q-1)(q+\epsilon)/2$, and the number
of $U$ so that $\psi(\det U)=-1$ is $q(q-1)(q-\epsilon)/2$.
$\square$
\enddemo

Proposition 3.8 of [18] tells us the following.

\proclaim{Proposition 5.9} We have
$\E_2|T_1(q^2)=(q+1)q^{2k-3}\E_2$.
\endproclaim

\proclaim{Proposition 5.10} We have
$$\E_{1,1}|T_1(q^2)=
\cases q^{2k-2}\E_{1,1}+q\E_1 &\text{if $\chi_q\not=1$,}\\
q^{-1}(q^2-1)\E_2+q^{2k-2}\E_{1,1}+q\E_1 &\text{if $\chi_q=1$.}
\endcases$$
\endproclaim

\demo{Proof}  
By Proposition 2.1 [5],
$$\align
&2\E_{1,1}(\tau)|T_1(q^2)\\
&\quad = q^{-3}\sum \chi_{\sigma}\chi_1(M,N)\\
&\qquad \cdot \det\left(M\pmatrix 1/q\\&1\endpmatrix(G^{-1}\tau+Y\,^tG)\,^tG^{-1}\pmatrix 1/q\\&1\endpmatrix+N\right)^{-k}\\
&\quad = 
q^{k-3}\sum \chi_{\sigma}\chi_1(M,N)\\
&\qquad \cdot \det\left(M\pmatrix 1/q\\&1\endpmatrix\tau
+N\pmatrix q\\&1\endpmatrix\,^tG+M\pmatrix 1/q\\&1\endpmatrix Y\right)^{-k},
\endalign$$
where the sum is over $SL_2(\Z)$-equivalence class representatives of
$\stufe/q$-type $\sigma$, $q^2$-type $1,1$.

Take $(M\ N)$ of $q^2$-type $1,1$; we can assume $q^2$ divides row 2 of $M$.

Suppose $q$ does not divide column 1 of $M$.  By symmetry, $q$ divides
row 2 of $N\pmatrix q\\&1\endpmatrix$ if and only if $q$ divides row 2 of $N$; hence
$q$ does not divide row 2 of $N\pmatrix q\\&1\endpmatrix$.  Thus $(M',N')=1$ where
$$(M'G\ N'\,^tG^{-1})=\pmatrix q\\&1\endpmatrix
\left(M\pmatrix 1/q\\&1\endpmatrix\ \ N\pmatrix q\\&1\endpmatrix+M\pmatrix 1/q\\&1\endpmatrix Y
\right).$$
Thus $M'$ is $q^2$-type $1,1$, since $q\nmid M'$, $q^2|\det M'$, and
$\chi_{\sigma}(M',N')=\chi_{\sigma}(M,N)$.
Reversing, take $(M'\ N')$ of $q^2$-type $1,1$ and set
$$(M\ N)=\pmatrix 1/q\\&1\endpmatrix E\left( M'G\pmatrix q\\&1\endpmatrix
\  (N'\,^tG^{-1}-M'GY)\pmatrix 1/q\\&1\endpmatrix\right).$$
We can assume $q^2$ divides row 2 of $M'$.
If $q^2$ divides row 1 of $EM'$, then we cannot solve $(n_1\ n_2)\equiv(m_1\ m_2)Y\ (q).$
So suppose $q^2$ divides row 2 of $EM'$ (this is the case for $q$ choices of $E$).
To have $M$ integral, choose the unique $G$ so that $q|m_2$.  Then $q\nmid m_1$,
and by symmetry, $q|n_3$ (and hence $q\nmid n_4$).  Now choose $y_2$ so that $n_2\equiv m_1y_2\ (q);$
then choose $y_1$ so that $n_1\equiv m_1y_1+m_2y_2\ (q^2).$  
So the contribution from this case is $q^{k-3}q^kq\E_{1,1}.$

Suppose $q$ divides column 1 of $M$.  So $\rank_qM\pmatrix 1/q\\&1\endpmatrix=1$
with $q$ dividing its 2nd row.  Thus $M\equiv\pmatrix 0&b\\0&0\endpmatrix\ (q)$, $q\nmid b$,
so by symmetry $N\equiv \pmatrix *&*\\c&0\endpmatrix \ (q)$ with $q\nmid c$ since $(M,N)=1$.
Thus with
$$(M'G\ N'\,^tG^{-1})=\pmatrix 1\\&1/q\endpmatrix\left(M\pmatrix 1/q\\&1\endpmatrix
\ \ N\pmatrix q\\&1\endpmatrix+M\pmatrix 1/q\\&1\endpmatrix Y\right),$$
we have $(M'\ N')$ integral and coprime.
Also, we know $M'$ is $q^2$-type 2 or $1,1$ since $q^2|\det M$, $q\nmid M'$, and
$\chi_{\sigma}(M',N')=\chi_{\sigma}(M,N)$.
Reversing, have
$$(M\ N)=\pmatrix 1\\&q\endpmatrix E\left(M'G\pmatrix q\\&1\endpmatrix
\ \ (N'\,^tG^{-1}-M'GY)\pmatrix 1/q\\&1\endpmatrix\right),$$
where $E\in\G_1$.

Suppose first that $M'$ is $q^2$-type $1,1$; we can assume $q$ divides 
row 2 of $M'$.  To have $\rank_qM=1$, we choose $G$ so that
$q\nmid m_2$.  There are $q$ choices of $E$ so that $q$ does not divide row 1 of $EM'$; then
have $q$ choices of $G$ so that $q\nmid m_2$.  Then by symmetry, $q\nmid n_3$.
Now, for any choice of $y_1$, choose the unique $y_2$ modulo $q$ so that
$n_1\equiv m_1y_1+m_2y_2\ (q)$.  Then $N$ is integral, and
the contribution from this situation is $q^{k-3}q^{-k}q^4\E_1$.

Now suppose $\rank_qM'=2$.  For each choice of $E$ we have $q$ choices of $G$
so that $M$ is $q^2$-type $1,1$.
Summing over those $Y$ so that $N$ is integral with $(M,N)=1$, we have
$$\sum_Y\chi_1(M,N)=\cases 
0&\text{if $\chi_q\not=1$,}\\
q(q-1)&\text{if $\chi_q=1$,}
\endcases$$
proving the proposition. $\square$
\enddemo

\proclaim{Proposition 5.11}  Let $\nu=+$ or $-$.  We have
$$
\E_{1,2,\nu}|T_1(q^2)
=\cases
q^{2k-2}\E_{1,2,\nu}+\nu \chi_{\stufe/q}(q) q^{k-3}(q^2-1)/2\cdot \E_2 &\text{if $\chi_q\not=1$,}\\
q^{2k-2}\E_{1,2,\nu}+\chi_{\stufe/q}(q) q^{k-3}(q^2-1)/2\cdot \E_2 &\text{if $\chi_q=1$.}\endcases$$
\endproclaim

\demo{Proof}  Take $(M\ N)$ of $q^2$-type $1,2,\nu$; we can assume $q$ divides row 2 of $M$.

Suppose $q$ does not divide column 1 of $M$.  Set
$$(M'G\ N'\,^tG^{-1})=\pmatrix q\\&1\endpmatrix\left(M\pmatrix 1/q\\&1\endpmatrix
\ \ N\pmatrix q\\&1\endpmatrix+M\pmatrix 1/q\\&1\endpmatrix Y\right).$$
So $(M\ N)$ is $q^2$-type $1,2$ and $\det M=\det M'$,
$\chi_{\sigma}(M',N')=\chi_{\sigma}(M,N)$.
Reversing, take $(M'\ N')$ $q^2$-type $1,2$; 
we find there are $q$ choices for $E$, and for each $E$ there are unique choices for $G$ and $Y$
so that $M,N$ are integral and coprime.
Thus the contribution from this case is
$$q^{k-3}q^kqE_{1,2,\nu}.$$

Suppose $q$ divides column 1 of $M$. 
Set
$$(M'G\ N'\,^tG^{-1})=\left(M\pmatrix 1/q\\&1\endpmatrix\ \ N\pmatrix q\\&1\endpmatrix
+M\pmatrix 1/q\\&1\endpmatrix Y\right).$$
We have $\rank_qM\pmatrix 1/q\\&1\endpmatrix=2$ and
$\chi_{\sigma}(M',N')=\chi_{\stufe/q}(q)\chi_{\sigma}(M,N)$.
Reversing, take $(M'\ N')$ of $q$-type 2; since $\rank_qM'=2$,
for each $G$, we can replace $(M'\ N')$ by $E(M'\ N')$ so that $q|m_4$ (and
hence $q\nmid m_2m_3$); as well, there are
unique $u,y_2$ modulo $q$ so that 
$$\pmatrix n_1\\n_3\endpmatrix\equiv M'G\pmatrix u\\y_2\endpmatrix\ (q).$$
Set $y_1=u+qu'$, $u'$ varying modulo $q$.  Then
$\psi_1(M,N)=\psi(a+m_2m_3u')$ where $a$ depends only on $M',N',u,y_2$.  So there
are $(q-1)/2$ ways to choose $u'$ so that $\nu\psi_1(M,N)=\psi\left(\frac{1}{q}\det M\right).$
Also, $\psi\left(\frac{1}{q}\det M\right)=\psi(\det M')$.  So, 
recalling that we have $q+1$ choices for $G$, the contribution from this
case is $\nu\chi_{\stufe/q}(q) q^{k-3}(q^2-1)/2\cdot\E_2$ when $\chi_q\not=1$, and 
$\chi_{\stufe/q}(q) q^{k-3}(q^2-1)/2\cdot\E_2$
when $\chi_q=1$. $\square$
\enddemo

\proclaim{Proposition 5.12} We have
$$
\E_{0,0}|T_1(q^2)
=\cases q\E_{0,0}+\E_0+\epsilon(\E_{0,2,+}-\E_{0,2,-}) &\text{if $\chi_q\not=1$,}\\
q^{-1}(q-1)\E_{1,1}+ q\E_{0,0}+\E_0+\epsilon(\E_{0,2,+}-\E_{0,2,-}) &\text{if $\chi_q=1$.}
\endcases$$
\endproclaim

\demo{Proof} With $(M\ N)$ of $q^2$-type $0,0$.  We have $q^2|M$ and $\rank_qN\pmatrix q\\&1\endpmatrix=1.$
We can adjust the representative so that
$$(M'G\ N'\,^tG^{-1})=\pmatrix 1\\&1/q\endpmatrix \left(M\pmatrix 1/q\\&1\endpmatrix
\ \ N\pmatrix q\\&1\endpmatrix+M\pmatrix 1/q\\&1\endpmatrix Y\right)$$
is a coprime pair.
We have $q^4|\det M$ so $q^2|\det M'$; thus $\rank_qM'=0$ or $M'$ is $q^2$-type $1,1$,
$\chi_{\sigma}(M',N')=\chi_{\sigma}(M,N)$.

Reversing,
first take $(M'\ N')$ of $q^2$-type $1,1$.  We find there are unique $E$ and $G$ so that
$N$ is integral and $M$ is $q^2$-type $0,0$; summing over $Y$ so that $\rank_qN=2$,
we have
$$\sum_Y\chi(\det N)=
\cases 0&\text{if $\chi_q\not=1$,}\\
q^2(q-1)&\text{if $\chi_q=1$.}\endcases$$
So the contribution from these terms is 0 if $\chi_q\not=1$, and 
$q^{k-3}\cdot q^{-k}\cdot q^2(q-1)\E_{1,1}$ otherwise.

Now suppose $(M'\ N')$ is $q^2$-type $0,0$; for each $E$ there is a unique $G$
so that $N$ is integral.  Then for all $Y$, $\chi_q(\det N)=\chi_q(\det N')$, so we
get a contribution of $(q+1)\E_{0,0}$.

Next suppose $(M'\ N')$ is $q^2$-type $0,1$.  Then for any $Y$, there are unique
$E$ and $G$ so that $(M\ N)$ is an integral pair of $q^2$-type $0,0$,
giving us a contribution of $\E_{0,1}$.

Last, suppose $(M'\ N')$ is $q^2$-type $0,2$.
For each $E$, there is a unique $G$ so that $q|n_1$.
Let $V=\F x_1\oplus\F x_2$ be equipped with the quadratic form given by $\frac{1}{q}M'\,^tN'$
relative to $(x_1\ x_2)$.  Thus relative to $(x_1'\ x_2')=(x_1\ x_2)E$, the quadratic form
is given by 
$$\frac{1}{q}EM'\,^tN'\,^tE=\frac{1}{q}\pmatrix m_2n_2&*\\*&*\endpmatrix\in\F^{2,2}$$
where $EM'G=\pmatrix m_1&m_2\\m_3&m_4\endpmatrix,$
$EN'\,^tG^{-1}=\pmatrix n_1&n_2\\n_3&n_4\endpmatrix$.
Also, $\F x_1'$ varies over all lines in $V$ as $E$ varies.
We need $q^2|m_2$ for $(M\ N)$ to be $q^2$-type $0,0$; there are 2 choices for $E$ that do this
if $V\simeq\H$, and 0 choices otherwise.
So the contribution $2\E_{0,2,\epsilon}.$
$\square$
\enddemo

\proclaim{Proposition 5.13} For $\nu=+$ or $-$, $\chi_q\not=1$,
$$
\E_{0,1,\nu}|T_1(q^2)=
\nu \chi_{\stufe/q}(q)q^{k-2}(q-1)/2\cdot\E_{1,1}+q\E_{0,1,\nu}+q/2\cdot\E_{0,2}-\epsilon/2\cdot(\E_{0,2,+}-\E_{0,2,-}).
$$
For $\nu=+$ or $-$, $\chi_q=1$, 
$$\align
&\E_{0,1,\nu}|T_1(q^2)\\
&\quad =q^{-1}(q-1)\E_{1,2,\nu}+\chi_{\stufe/q}(q)q^{k-2}(q-1)/2\cdot\E_{1,1}
+q\E_{0,1,\nu}\\
&\qquad
+q/2\cdot \E_{0,2}-1/2\cdot(\E_{0,2,\epsilon}-\E_{0,2,-\epsilon}).
\endalign$$
\endproclaim

\demo{Proof}  Take $(M\ N)$ $q^2$-type $0,1,\nu$.

Suppose $(M',N')=2$ where
$$(M'G\ N'\,^tG^{-1})=\left(M\pmatrix 1/q\\&1\endpmatrix
\ \ N\pmatrix q\\&1\endpmatrix+M\pmatrix 1/q\\&1\endpmatrix Y\right).$$
Since $\rank_qN\pmatrix q\\&1\endpmatrix=1$, we must have  
$\rank_qM'=1$.  We know $q^3|\det M$, so
$q^2|\det M'$, thus $M'$ must be $q^2$-type $1,1$ and
$\chi_{\sigma}(M',N')=\chi_{\stufe/q}(q)\chi_{\sigma}(M,N)$.
Reversing, with $(M'\ N')$ $q^2$-type $1,1$, 
we have a unique choice of $G$ and $q(q-1)/2$ choices of $Y$ so that $(M\ N)$
is an integral pair of $q^2$-type $0,1,\nu$.  So the contribution from these terms is
$\nu\chi_{\stufe/q}(q) q^{k-3}q(q-1)/2\cdot\E_{1,1}$ when $\chi_q\not=1$, 
and $\chi_{\stufe/q}(q)q^{k-3}q(q-1)/2\cdot\E_{1,1}$
when $\chi_q=1$.

Now suppose $\rank_q\left(M\pmatrix 1/q\\&1\endpmatrix
\ N\pmatrix q\\&1\endpmatrix\right)=1.$  Adjust the representative so that
$$(M'G\ N'\,^tG^{-1})=\pmatrix 1\\&\frac{1}{q}\endpmatrix E\left(M\pmatrix 1/q\\&1\endpmatrix
\ \ N\pmatrix q\\&1\endpmatrix+M\pmatrix \frac{1}{q}\\&1\endpmatrix Y\right)$$
is integral; note that $\rank_qN'=2$.  Also, $q^3|\det M$ so $q|\det M'$.
Thus $\rank_qM'=0$ or 1, and $\chi_{\sigma}(M',N')=\chi_{\sigma}(M,N)$.

Suppose 
first that $(M'\ N')$ is of $q$-type 1.  There are unique choices of $E$ and $G$, and
$q^2(q-1)$ choices of $Y$ so that $(M\ N)$ is $q^2$-type $0,1$, and
$(M\ N)$ is $q^2$-type $0,1,\nu$ if and only if $(M'\ N')$ is $q^2$-type $1,2,\nu$.
Then
$$\sum_Y\chi(\det N)=\cases 0&\text{if $\chi_q\not=1$,}\\
q^2(q-1)&\text{if $\chi_q=1$.}\endcases$$
Hence the contribution from this case is 0 if $\chi_q\not=1$, and
$q^{k-3}\cdot q^{-k}\cdot q^2(q-1)\E_{1,2,\nu}$ if $\chi_q=1$.

Next, suppose $\rank_qM'=0$.  To have $M$ $q^2$-type $0,1$, we need $M'$ to be $q^2$-type $0,1$ or $0,2$.
First suppose $M'$ is $q^2$-type $0,1$.  We have $q$ choices of $E$ and a unique choice of $G$
so that $M$ is $q^2$-type $0,1$ and $N$ is integral; for $(M'\ N')$ of $q^2$-type $0,1,\nu'$,
we have $(M\ N)$ of $q^2$-type $0,1,\nu'$ for all choices of $Y$.  So the
contribution from these terms is
$q^{k-3}q^{-k}qq^3\E_{0,1,\nu}=q\E_{0,1,\nu},$ for $\chi_q\not=1$ or 1.

Finally, suppose $M'$ is $q^2$-type $0,2$.  
For each $E$ there is a unique $G$ so that $N$ is integral.
We know $M$ is $q^2$-type $0,1$ if and only if $q^2\nmid m_2$
where $EM'G=\pmatrix m_1&m_2\\m_3&m_4\endpmatrix$.
Let $V=\F x_1\oplus\F x_2$ be equipped with the quadratic form given by $\frac{1}{q}M'\,^tN'$
relative to $(x_1\ x_2)$.  So relative to $(x_1'\ x_2')=(x_1\ x_2)\,^tE$, the quadratic form
is given by $\frac{1}{q}EM'\,^tN'\,^tE$, and $\F x_1'$ varies over all lines in $V$ as $E$
varies.  We have $q^2\nmid m_2$ if and only if $\F x_1'$ is anisotropic.  
Also, we have $\psi(\det N)=\psi(\det N')$, and $(M'\ N')$ is $q^2$-type $0,1,\nu$
if and only if $\nu\psi_1\left(\frac{1}{q}M,N\right)=\psi(\det N)$, or equivalently,
if and only if $\nu\psi(m_2n_2/q)=1$ where
$EN'\,^tG^{-1}=\pmatrix n_1&n_2\\n_3&n_4\endpmatrix.$  
We know that $\H$ contains $(q-1)/2$ lines representing only
squares, $(q-1)/2$ lines representing only (non-zero) non-squares.  Similarly, $\A$ contains
$(q+1)/2$ lines representing only squares, $(q+1)/2$ lines representing only (non-zero) non-squares.
Consequently, the contribution from this case is
$$(q-1)/2\cdot\E_{0,2,\epsilon}+(q+1)/2\cdot\E_{0,2,-\epsilon}$$
for $\chi_q\not=1$ or 1.
$\square$
\enddemo

\proclaim{Proposition 5.14} We have
$$\align
&\E_{0,2,\epsilon}|T_1(q^2)\\
&\quad = {\cases
\epsilon q^{-2}(q^2-1)\E_2
+\epsilon\chi_{\stufe/q}(q) q^{k-2}(q-1)/2\big(\E_{1,2,+}-\E_{1,2,-}\big) &\text{if $\chi_q\not=1$,}\\
q^{-2}(q^2-1)\E_2+\chi_{\stufe/q}(q)q^{k-2}(q-1)/2\cdot\E_{1,2}&\text{if $\chi_q=1$;}
\endcases}
\endalign$$
and
$$\E_{0,2,-\epsilon}|T_1(q^2)= 
\cases
-\epsilon \chi_{\stufe/q}(q)q^{k-2}(q-1)/2\big(\E_{1,2,+}-\E_{1,2,-}\big) &\text{if $\chi_q\not=1$,}\\
\chi_{\stufe/q}(q)q^{k-2}(q-1)/2\cdot\E_{1,2} &\text{if $\chi_q=1$.}
\endcases$$
\endproclaim

\demo{Proof}  Take $(M\ N)$ of $q^2$-type $0,2,\nu$ where $\nu=+$ or $-$.

Suppose $\rank_q\left(M\pmatrix 1/q\\&1\endpmatrix \ N\pmatrix q\\&1\endpmatrix
\right)=1$.  Then adjusting the representative for $SL_2(\Z)(M\ N)$, we have
$$(M'G\ N'\,^tG^{-1})=\pmatrix 1\\&1/q\endpmatrix\left(M\pmatrix 1/q\\&1\endpmatrix
\ \ N\pmatrix q\\&1\endpmatrix+M\pmatrix 1/q\\&1\endpmatrix Y\right)$$
integral pair with $\rank_qM'=2$.  
Also, $\chi_{\sigma}(M',N')=\chi_{\sigma}(M,N)$.
Reversing, take $(M'\ N')$ $q$-type 2.  For each $E$ there is a unique $G$ so that
$M$ is $q^2$-type $0,2$.  To have $N$ integral, we choose $y_1$ so that
$n_1\equiv m_1y_1\ (q)$ where $EM'G=\pmatrix m_1&m_2\\m_3&m_4\endpmatrix$
and $EN'\,^tG^{-1}=\pmatrix n_1&n_2\\n_3&n_4\endpmatrix$.  Then
$$\align
\psi(\det N)&=\psi(-(n_2-m_1y_1)(n_3-m_3y_1-m_4))\\
&=\psi(-m_1m_4)\psi((y_2-\overline m_1n_2)(y_2-\overline m_4n_3-\overline m_4m_3y_1))\\
&=\epsilon\psi(\det M')\psi^2(y_2-\overline m_1n_2).
\endalign$$
Hence for $q(q-1)$ choices of $Y$,
$(M\ N)$ is $q^2$-type $0,2,\epsilon$, and
there are no choices of $Y$ so that $(M\ N)$ is $q^2$-type $0,2,-\epsilon$.
So the contribution to $\E_{0,2,\epsilon}|T_1(q^2)$ is $\epsilon q^{-2}(q^2-1)\E_2$
if $\chi_q\not=1$, $q^{-2}(q^2-1)\E_2$ if $\chi_q=1$;
the contribution to $\E_{0,2,-\epsilon}|T_1(q^2)=0$.

Suppose $(M',N')=1$ where
$$(M'G\ N'\,^tG^{-1})=\left(M\pmatrix 1/q\\&1\endpmatrix
\ \ N\pmatrix q\\&1\endpmatrix+M\pmatrix q\\&1\endpmatrix Y\right).$$
So $M'$ is $q^2$-type $1,2$, and $\chi_{\sigma}(M',N')=\chi_{\stufe/q}(q)\chi_{\sigma}(M,N).$
Reversing, take $(M'\ N')$ $q^2$-type $1,2$.  There is a unique $G$ so that $q|M$, and
$q(q-1)/2$ choices for $Y$ so that $N$ is integral with $(M\ N)$ of
$q^2$-type $0,2,\nu$.  Then 
when $(M'\ N')$ is $q^2$-type $1,2,\nu'$, we have
$$\psi_1(M',N')=\nu'\left(\frac{1/q\cdot\det M'}{q}\right)
=\nu'\left(\frac{\det (M/q)}{q}\right)=\nu\nu'\left(\frac{\det N}{q}\right).$$
So when
$\chi_q\not=1$, the contribution is $\nu\chi_{\stufe/q}(q) q^{k-2}(q-1)/2\cdot(\E_{1,2,+}-\E_{1,2,-}).$
When $\chi_q=1$, the contribution is $\chi_{\stufe/q}(q)q^{k-2}(q-1)/2\cdot\E_{1,2}$.
$\square$
\enddemo

To prove our main theorem,
we will only be applying $A_0(q)$ in the case that $\chi_q=1$, and $A_1(q), A_2(q)$
in the case that $\chi_q\not=1$.  

\proclaim{Proposition 5.15} For $\chi_q=1$, we have
$$\E_2|A^*_0(q)=(q+1)\E_2+\chi_{\stufe/q}(q)q^k\E_{1,2}.$$
\endproclaim

\demo{Proof}  Note that since $\rank_qY\le 1$, we have $\rank_qMY=\rank_qY\le 1$.

In the case that $\rank_q MY=1$, we adjust the representative so that
$$(M'\,^tG^{-1}\ N'G)=\pmatrix 1\\&q\endpmatrix E'(M\ N+MY/q)$$
is integral; then
$(M'\ N')$ is $q^2$-type $1,2$, and $\chi_{\sigma}(M',N')
=\chi_{\stufe/q}(q)\chi_{\sigma}(M,N)$.
Reversing, 
there are unique $E,G,Y$ so that $M,N$ are integral; then $\rank_qM=2$,
and we get a contribution of
$$\chi_{\stufe/q}(q)q^k\E_{1,2}.$$

Now say $q|MY$; thus $Y$ must be 0, and we set
$(M'\,^tG^{-1}\ N'G)=(M\ N)$.  So the contribution from this case is $(q+1)\E_2$.
$\square$
\enddemo

\proclaim{Proposition 5.16} For $\chi_q=1$,
$$\E_{1,1}|A^*_0(q)=
2q\E_{1,1}+\chi_{\stufe/q}(q)q^{k+1}\E_{0,1}.$$
\endproclaim

\demo{Proof}  Take $(M\ N)$ of $q^2$-type $1,1$; we can assume $q^2$ divides row 2 of $M$.

Suppose $MY/q$ is integral.  So $q$ does not divide row 1 of $M$ or row 2
of $N+MY/q$; thus $(M'\,^tG^{-1}\ N'G)=(M\ N+MY/q)$ is of $q^2$-type $1,1$.
Also, $\chi_{\sigma}(M',N')=\chi_{\sigma}(M,N)$.
Reversing, take $(M'\ N')$ of $q^2$-type $1,1$; we can assume $q^2$ divides row 2 of $M'$.
To have $N$ integral, we need $q$ dividing column 1 of $M'\,^tG^{-1}Y$.
We have a unique choice of $G$ so that $q$ divides column 1
of $M'\,^tG^{-1}$, then we can choose any $y$.
For the other $q$ choices of $G$, we must take $y\equiv0\ (q)$.
Also, $\psi_1(M,N)=\psi_1(M',N')$.  So the contribution from this case is
$2q\E_{1,1}$.

Suppose $MY/q$ is not integral; then with
$$(M'\,^tG^{-1}\ N'G)=\pmatrix q\\&1\endpmatrix (M\ N-MY/q),$$
$(M'\ N')$ is $q^2$-type $0,1$, and $\chi_{\sigma}(M',N')=\chi_{\stufe/q}(q)\chi_{\sigma}(M,N)$.
Reversing, take $(M'\ N')$ of $q^2$-type $0,1$; assume $q^2$ divides row 2 of $M'$.
Then we have $q$ choices of $E$, and for each a unique choice of $G$ and $Y$ so that
$(M\ N)$ is integral of $q^2$-type $1,1$.
Thus the contribution from this case is
$\chi_{\stufe/q}(q)q^{k+1}(\E_{0,1,+}-\E_{0,1,-}).$
$\square$
\enddemo

\proclaim{Proposition 5.17}  For $\chi_q=1$, $\nu=+$ or $-$,
$$\align
\E_{1,2,\nu}|A_0^*(q)
&=\chi_{\stufe/q}(q)q^{-k}(q^2-1)/2\cdot \E_2
+(q+1)\E_{1,2,\nu}+(q-1)/2\cdot\E_{1,2}\\
&\quad+\chi_{\stufe/q}(q)q^{k+1}/2\cdot\E_{0,2}
-\chi_{\stufe/q}(q)q^k/2\cdot(\E_{0,2,\epsilon}-\E_{0,2,-\epsilon}).
\endalign$$
\endproclaim

\demo{Proof}  
When $MY/q$ is integral and $\rank_q(M\ N+MY/q)=1$, we set  
$$(M'\,^tG^{-1}\ \ N'G)=\pmatrix 1\\&1/q\endpmatrix(M\ N+MY/q).$$
Thus $\rank_qM'=2$, $\chi_{\sigma}(M',N')=\chi_{\stufe/q}(q)\chi_{\sigma}(M,N)$.
Reversing, take $(M'\ N')$ of $q$-type 2.  To have $(M\ N)$ integral of $q^2$-type
$1,2,\nu$, for each choice of $E$ we have a unique choice of $G$ and $q-1$
choices for $y$, giving us a contribution of $\chi_{\stufe/q}(q)q^{-k}(q^2-1)/2\cdot\E_2$.

When $MY/q$ is integral and $\rank_q(M\ N+MY/q)=2$, set
$$(M'\,^tG^{-1}\ \ N'G)=(M\ N+MY/q).$$
Then $\chi_{\sigma}(M',N')=\chi_{\sigma}(M,N)$.  Reversing, take $(M'\ N')$ of $q^2$-type $1,2$.
If $q|y$ then can choose any $G$, and $(M\ N)$ is $q^2$-type $1,2,\nu$ provided $(M'\ N')$ is,
giving a contribution of $(q+1)\E_{1,2,\nu}$.
There are $(q-1)/2$ choices for $y\not\equiv0\ (q)$ and for each a unique choice of $G$
so that $(M\ N)$ is integral of $q^2$-type $1,2,\nu$, giving us a contribution of
$(q-1)/2\cdot\E_{1,2}.$

Say $MY/q$ is not integral.  
Assuming $q$ divides row 2 of $M'$, we set
$$(M'\,^tG^{-1}\ \ N'G)=\pmatrix q\\&1\endpmatrix(M\ N+MY/q);$$
 so $M'$ is
$q^2$-type $0,2$, $\chi_{\sigma}(M',N')=\chi_{\stufe/q}(q)\chi_{\sigma}(M,N)$.
Reversing, take $(M'\ N')$ of $q^2$-type $0,2$.  For each choice of $E$ there
are unique choices for $G$ and $y$ so that $(M\ N)$ is an integral pair of $q^2$-type
$1,2$.  Looking at $V\simeq\frac{1}{q}M'\,^tN'$, to have $(M\ N)$ of $q^2$-type
$1,2,\nu$, we find there are $(q-1)/2$ choices for $E$ when $V\simeq\H$ (modulo $q$)
and $(q+1)/2$ choices for $E$ when $V\simeq\A$ (modulo $q$).  Thus
the contribution
from these terms is 
$$\chi_{\stufe/q}(q)\big(q^k(q-1)/2\cdot \E_{0,2,\epsilon}
+q^k(q+1)/2\cdot\E_{0,2,-\epsilon}\big).$$
Combining the terms yields the result.
$\square$
\enddemo

\proclaim{Proposition 5.18} For $\chi_q=1$,
$$\E_{0,0}|A_0^*(q)=q(q+1)\E_{0,0}.$$
\endproclaim

\demo{Proof} Set $(M'\,^tG^{-1}\ N'G)=(M\ N+MY/q).$  So
$N'G\equiv N\ (q)$, and $(M'\ N')$ is $q^2$-type $0,0$.  Reversing,
for $(M'\ N')$ is $q^2$-type $0,0$, $(M\ N)=(M'\,^tG^{-1}\ N'G-M'\,^tG^{-1}Y/q)$
is a coprime pair of $q^2$-type $0,0$ for each choice of $G,Y$. $\square$
\enddemo

\proclaim{Proposition 5.19}  For $\chi_q=1$, $\nu=+$ or $-$,
$$\E_{0,1,\nu}|A_0^*(q)=
 q\E_{0,1,\nu}+q(q-1)/2\cdot\E_{0,1}+\chi_{\stufe/q}(q)q^{1-k}(q-1)/2\cdot\E_{1,1}.$$
\endproclaim

\demo{Proof}
Take $(M\ N)$ of $q^2$-type $0,1,\nu$; we can assume $q^2$ divides row 2 of $M$.

Suppose $\rank_q(N+MY/q)=2$; set
$$(M'\,^tG^{-1}\ N'G)=(M\ N+MY/q).$$
So $(M'\ N')$ is $q^2$-type $0,1$, $\chi_{\sigma}(M',N')=\chi_{\sigma}(M,N)$.
Reversing, take $(M'\ N')$ of $q^2$-type $0,1$ and set
$$(M\ N)=(M'\,^tG^{-1}\ N'G-M'\,^tG^{-1}Y/q).$$
We can assume $q^2$ divides row 2 of $M'$.
Write $M'\,^tG^{-1}=q\pmatrix m_1&m_2\\qm_3&qm_4\endpmatrix$,
$N'G=\pmatrix n_1&n_2\\n_3&n_4\endpmatrix$.  First say $q|m_1$ (this is the case for
1 choice of $G$).  So $q\nmid m_2$ and by symmetry, $q|n_4$, so $q\nmid n_2n_3$.
Then $\det N\equiv -n_2n_3\not\equiv0\ (q)$, $\det N'\equiv\det N\ (q)$,
and $\psi_1(M/q,N)=\psi_1(M'/q,N')$.  So $(M\ N)$ is $q^2$-type $0,1,\nu$ if and only
if $(M'\ N')$ is.  We have $q$ choices for $Y$, so the contribution from these terms
is $q\E_{0,1,\nu}$.
Now suppose $q\nmid m_1$; there are $q$ choices for $G$ so that this is the case.
By symmetry, $q\nmid n_4$.  To have $(M\ N)$ $q^2$-type $0,1,\nu$, we need
$$\nu\left(\frac{\det N}{q}\right)=\nu\left(\frac{(n_1-m_1u)n_4-n_2n_3}{q}\right)
=\psi_1(M/q,N)=\left(\frac{m_1n_4}{q}\right);$$
so we have $(q-1)/2$ choices for $Y$.
We have $q$ choices for $G$, so the contribution from these terms is
$q(q-1)/2\cdot \E_{0,1}$.

Say $\rank_q(N+MY/q)=1$.  We know $q$ does not divide row 2 of 
$N+MY/q$, so adjusting our representative $(M\ N)$,
$$(M'\,^tG^{-1}\ N'G)=\pmatrix 1/q\\&1\endpmatrix (M\ N+MY/q)$$
is an integral pair of $q^2$-type $1,1$ and $\chi_{\sigma}(M',N')=
\chi_{\stufe/q}(q)\chi_{\sigma}(M,N)$.
Reversing, take $(M'\ N')$ $q^2$-type $1,1$; 
there is 1 choice of $E$, $q$ choices of $G$, and $(q-1)/2$ choices of $y$
so that $(M\ N)$ is $q^2$-type $0,1,\nu$.
So the contribution from these terms is
$q^{1-k}(q-1)/2\cdot\E_{1,1}$.
$\square$
\enddemo

\proclaim{Proposition 5.20} For $\chi_q=1$, we have
$$\E_{0,2,\epsilon}|A_0^*(q)=(q+1)^2/2\cdot\E_{0,2,\epsilon}+(q^2-1)/2\cdot\E_{0,2,-\epsilon}
+\chi_{\stufe/q}(q)q^{1-k}(q-1)/2\cdot\E_{1,2},$$
and
$$\E_{0,2,-\epsilon}|A_0^*(q)=(q-1)^2/2\cdot\E_{0,2,\epsilon}+(q^2-1)/2\cdot\E_{0,2,-\epsilon}
+\chi_{\stufe/q}(q)q^{1-k}(q-1)/2\cdot\E_{1,2}.$$
\endproclaim

\demo{Proof} Take $\nu=+$ or $-$.  Note that $\rank_q(N+MY/q)\ge 1$ since
$q$ does not divide column 2 of $N$.

When $\rank_q(N+MY/q)=2$.  Set $(M'\,^tG^{-1}\ N'G)=(M\ N+MY/q)$.
So $(M'\ N')$ is $q^2$-type $0,2$, and $\chi_{\sigma}(M',N')=\chi_{\sigma}(M,N)$.
Reversing, take $(M'\ N')$ of $q^2$-type $0,2,\nu'$, and
set $(M\ N)=(M'\,^tG^{-1}\ N'G-M'\,^tG^{-1}Y/q)$;
write $M'\,^tG^{-1}=q\pmatrix m_1&m_2\\m_3&m_4\endpmatrix$, $N'G=\pmatrix n_1&n_2\\n_3&n_4\endpmatrix.$
To have $(M\ N)$ of $q^2$-type $0,2,\nu$, we need $q\nmid\det N$ and 
$\nu\left(\frac{\det MN/q}{q}\right)=1.$  
Set $M_0=M'/q$, $d=\det M_0$.  
Let $V=\F x_1\oplus\F x_2\simeq\overline M_0N'$.  Then with $(x_1'\ x_2')=(x_1\ x_2)G$,
$\F x_2'$ varies over all lines in $V$.
Writing $M_0\,^tG^{-1}=\pmatrix m_1&m_2\\m_3&m_4\endpmatrix$,
$N'G=\pmatrix n_1&n_2\\n_3&n_4\endpmatrix$, we have $Q(x_2')=\overline d(m_1n_4-m_3n_2)$
and $$\det\overline M_0N\equiv u\overline d(m_3n_2-m_1n_4)+\det\overline M_0N'\ (q).$$
Suppose first that $\nu'=\epsilon$; then for 2 choices of $G$, $q|m_3n_2-m_1n_4$ and
hence $(M\ N)$ is $q^2$-type $0,2,\epsilon$ for all choices of $y$.  For the other
$q-1$ choices of $G$, $(M\ N)$ is $q^2$-type $0,2,\epsilon$ for $(q-1)/2$ choices
of $y$, and $q^2$-type $0,2,-\epsilon$ for $(q-1)/2$ choices of $y$
(and $(M,N)\not=1$ for 1 choice of $u$).   Suppose now that $\nu'=-\epsilon$; then
for each choice of $G$, we have $(M\ N)$ of $q^2$-type $0,2,\epsilon$ for
$(q-1)/2$ choices of $u$, and $q^2$-type $0,2,-\epsilon$ for $(q-1)/2$ choices of $y$.
So the contribution from these terms is $(q+1)^2/2\cdot\E_{0,2,\epsilon}
+(q^2-1)/2\cdot\E_{0,2,-\epsilon}$ when $\nu=\epsilon$, and
$(q-1)^2/2\cdot\E_{0,2,\epsilon}+(q^2-1)/2\cdot\E_{0,2,-\epsilon}$ when
$\nu=-\epsilon$.

When $\rank_q(N+MY/q)=1$, adjust the representative $(M\ N)$
to assume $q$ divides row 1 of $N+MY/q$, and set
$$(M'\,^tG^{-1}\ N'G)=\pmatrix 1/q\\&1\endpmatrix(M\ N+MY/q).$$
So is $q^2$-$q^2$-type $1,2$ and $\chi_{\sigma}(M',N')=\chi_{\stufe/q}(q)\chi_{\sigma}(M,N)$.
Reversing, take $(M'\ N')$ of $q^2$-type $1,2$; 
to have $(M\ N)$ of $q^2$-type $0,2,\nu$, we have 1 choice of $E$, $q$ choices of $G$,
and $(q-1)/2$ choices of $u$.  Then the contribution from these terms is
$\chi_{\stufe/q}(q)q^{1-k}(q-1)/2\cdot\E_{1,2}$. $\square$
\enddemo

\proclaim{Proposition 5.21} For $\chi_q\not=1$, $\E_2|A^*_1(q)=\chi_{\stufe/q}(q)q^{k+1}\,\E_{1,2}.$
\endproclaim

\demo{Proof}  Begin with $(M\ N)$ of $q^2$-type 2.  Note that if $\rank_qMY=0$
then $Y\equiv 0\ (q)$, so $\left(\frac{y_1}{q}\right)=0$; thus no contributions
will come from $\rank_qMY=0$.

Suppose $\rank_qMY=2$.  Then we take
$$(M'G\ N'\,^tG^{-1})=q(M\ N+MY/q);$$
note that $\chi_{\sigma}(M',N')=\chi_{\sigma}(M,N)$.
Reversing, take $(M'\ N')$ of $q^2$-type $0,2$. 
To ensure $N$ is integral, we need to take
$$Y\equiv G^{-1}M_0^{-1}N'\,^tG^{-1}\equiv\frac{1}{d}\pmatrix m_4&-m_2\\-m_3&m_1\endpmatrix
\pmatrix n_1&n_2\\n_3&n_4\endpmatrix\ (q)$$
where $M_0=M'/q$ and $d=\det M_0$.  
Let $V=\F x_1\oplus\F x_2\simeq M_0^{-1}N'$,
and let $(x_1'\ x_2')=(x_1\ x_2)\,^tG^{-1}$.  Then 
as $G$ (and thus $Y$) varies, $\F x_1'$ varies over all lines in $V$, and
since $V$ is regular,
$$\sum_Y\left(\frac{y_1}{q}\right)\left(\frac{\det M}{q}\right)
=\left(\frac{\det (M'/q)}{q}\right) \sum_{\F x'_1}\left(\frac{Q(x'_1)}{q}\right)=0.$$

Suppose $\rank_qMY=1$.  Since we only need to consider $Y$
where $q\nmid y_1$, and we know $\rank_q M=2$, we know that $q$ does not divide column 1
of $MY$.  Adjust the representative $(M\ N)$ so that $q$ divides row 1 of $MY$.  Then
$$(M'G\ N'\,^tG^{-1})=\pmatrix 1\\&q\endpmatrix (M\ N+\frac{1}{q}MY)$$
is $q^2$-type $1,2$ and
$\chi_{\sigma}(M',N')=\chi_{\stufe/q}(q)\chi_{\sigma}(M,N)$.
Reversing, take $(M'\ N')=1$ of $q^2$-type $1,2$.
To have $N$ integral, we need $Y\equiv\overline M_0\pmatrix q\\&1\endpmatrix N'\ (q)$
and for $q$ choices of $G$ we have $y_1\not\equiv0\ (q)$.  For each such $Y$, we have
$\left(\frac{y_1}{q}\right)\left(\frac{\det M}{q}\right)=\psi_1(M',N')$.
So  the contribution
from this case is $\chi_{\stufe/q}(q)q^{k+1}\E_{1,2}$.
$\square$
\enddemo

\proclaim{Proposition 5.22} We have $\E_{1,1}|A^*_1(q)=\chi_{\stufe/q}(q)q^{k+2}\E_{0,1}.$
\endproclaim

\demo{Proof} Take $(M\ N)$ of $q^2$-type $1,1$ with $q^2$ dividing row 2 of $M$.
So $\rank_qMY=0$ or 1.

Suppose $\rank_qMY=0$.
Set 
$$(M'G\ \ N'\,^tG^{-1})=\left(M\ \ N+MY/q\right);$$
 then $(M',N')=1$
and $\chi_{\sigma}(M',N')=\chi_{\sigma}(M,N)$.
Reversing, choose $(M'\ N')$ of $q^2$-type $1,1$.
We can assume $q^2$ divides row 2 of $M'$; then to have $N$ integral,
we need $q$ dividing row 1 of $M'GY$.  Then $\psi_1(M,N)=\psi_1(M',N')$ and
$\sum_Y\left(\frac{y_1}{q}\right)\psi_1(M,N)=0.$

Suppose $\rank_qMY=1$ and
$\rank_q\pmatrix 1\\&1/q\endpmatrix MY+\pmatrix q\\&1\endpmatrix N=2.$
Then
$$(M'G\ N'\,^tG^{-1})
=\pmatrix q\\&1\endpmatrix \left(M\  \ N+MY/q\right)$$
is a coprime pair; so $M'$ is $q^2$-type $0,1$, 
and $\chi_{\sigma}(M',N')=\chi_{\stufe/q}(q)\chi_{\sigma}(M,N)$.
Reversing, choose $(M'\ N')$ of $q^2$-type $0,1$.
To have $N$ integral, we need to choose $E$ so that $q$ divides row 2 of $EM'$.  Then
for 1 choice of $G$ there are unique choices for $y_1$ and $y_2$ so that $N$ is integral,
and in this case
$$\sum_Y\left(\frac{y_1}{q}\right)\psi_1(M,N)=q\left(\frac{\det N'}{q}\right).$$
For the other $q$ choices of $G$, $y_1$ can vary freely, and so the sum on corresponding $Y$
is 0.  So the contribution from these terms is $\chi_{\stufe/q}(q)q^{k+2}\,\E_{0,1}$.

Suppose $\rank_qMY=1$, 
$\rank_q \left(\pmatrix 1\\&1/q\endpmatrix MY+\pmatrix q\\&1\endpmatrix N\right)=1.$
We can adjust the equivalence class representative so that 
$$(M'G\ N'\,^tG^{-1})
=\pmatrix q\\&1/q\endpmatrix \left(M\ N+MY/q\right)$$
is a coprime pair; then $(M'\ N')$ is of $q^2$-type $1,1$, and
$\chi_{\sigma}(M',N')=\chi_{\sigma}(M,N)$.
Hence $q$ does not divide row 2 of $M'$, or row 1 of $N'$, and so $(M',N')=1$.
Also, $q^2$ does not divide row 1 of $M'$, but
$\rank_q\pmatrix 1/q\\&1\endpmatrix M'=1.$
Reversing, take $(M'\ N')$ of $q^2$-type $1,1$.
To have $M$ integral, we need to choose the unique $E$ so that $q$ divides
row 1 of $EM'$.  For 1 choice of $G$, we need to choose $y_1\equiv0\ (q)$ to have
$N$ integral, and then $\left(\frac{y_1}{q}\right)=0$.  For the other $q$ choices
of $G$, $y_1$ varies freely, and then 
$$\sum_Y\left(\frac{y_1}{q}\right)\psi_1(M,N)=\psi_1(M',N')\sum_Y\left(\frac{y_1}{q}\right)=0.$$
Thus there is no contribution from this final case.
$\square$
\enddemo

\proclaim{Proposition 5.23} For $\chi_q\not=1$, $\nu=+$ or $-$, we have
$$\align
\E_{1,2,\nu}|A^*_1(q)
&=-\nu q\E_{1,2,\epsilon\nu}+\epsilon \chi_{\stufe/q}(q)
q^{1-k}(q^2-1)/2\cdot\E_2
+\chi_{\stufe/q}(q)q^{k+2}/2\cdot\E_{0,2}\\
&\quad -\epsilon\chi_{\stufe/q}(q) q^{k+1}/2\cdot(\E_{0,2,+}-\E_{0,2,-}).
\endalign$$
\endproclaim

\demo{Proof}  Take $(M\ N)$ of $q^2$-type $1,2,\nu$; we can
assume that $q$ divides row 2 of $M$.

Suppose $q$ divides $MY$ and
$$(M'G\ N'\,^tG^{-1})= \left(M\ N+MY/q\right)$$
is a coprime pair.  Then $(M'\ N')$ is $q^2$-type $1,2$ and
$\chi_{\sigma}(M',N')=\chi_{\sigma}(M,N)$.
Reversing, take $(M'\ N')$ of $q^2$-type $1,2,\nu'$.
To have $N$ integral with $q\nmid y_1$, we have $q$ choices for $G$ and
$q-1$ choices for $Y$ with $y_1$ free.  Then 
$\psi_1(M,N)=\left(\frac{-m_2n_3-dy_1}{q}\right)$ where
$M'G=\pmatrix m_1&m_2\\m_3&m_4\endpmatrix$, $N'\,^tG^{-1}=\pmatrix n_1&n_2\\n_3&n_4\endpmatrix$
and $d=\frac{1}{q}\det M'$.  To have $(M\ N)$ of $q^2$-type $1,2,\nu$,
we need to choose $y_1$ so that for some
$u\not\equiv0\ (q)$,
$$y_1\equiv\cases -u^2-\overline d m_2n_3&\text{if $\nu=+$,}\\
-\omega u^2-\overline d m_2n_3&\text{if $\nu=-$.}\endcases$$
Then summing over the corresponding $Y$ we have
$$\align
\sum_Y\left(\frac{y_1}{q}\right)\psi_1(M,N)
&=\nu\left(\frac{d}{q}\right)\sum_Y\left(\frac{y_1}{q}\right)\\
&=\nu\nu'\psi_1(M',N')\sum_Y\left(\frac{y_1}{q}\right)\\
&=\nu\nu'\psi_1(M,N)\sum_Y\left(\frac{y_1}{q}\right),
\endalign$$
and it is a standard exercise in number theory to verify that
$$\align
\sum_Y\left(\frac{y_1}{q}\right)
&=-\epsilon\nu\left(1+\nu\left(\frac{dm_2n_3}{q}\right)\right)/2\\
&={\cases -\epsilon\nu&\text{if $\nu'=\epsilon\nu$,}\\
0&\text{otherwise.}\endcases}
\endalign$$
So the contribution from this case is $-\nu q\E_{1,2,\epsilon\nu}$.

Suppose $q|MY$ and $\rank_q(M\ N+MY/q)=1.$ Then
$$(M'G\ N'\,^tG^{-1})=\pmatrix 1\\&1/q\endpmatrix (M\ N+MY/q)$$
is a coprime pair with $M'$ of
$q^2$-type 2 and $\chi_{\sigma}(M',N')=\chi_{\stufe/q}(q)\chi_{\sigma}(M,N)$.  Reversing,
choose $(M'\ N')$ of $q$-type 2.
For each $E$ there is a unique $G$ so that $q|y_1$ when $Y$ is chosen so that
$N$ is integral.  So for each $E$ we have $q$ other choices for $G$, and choosing
$Y$ so that $(M\ N)$ is an integral pair of $q^2$-type $1,2,\nu$, we get
$$\sum_Y\left(\frac{y_1}{q}\right)\psi_1(M,N)=\frac{\epsilon(q-1)}{2}\left(\frac{\det M'}{q}\right).$$
This gives us a contribution of
$\epsilon\chi_{\stufe/q}(q)q^{1-k}(q^2-1)/2\cdot\E_2.$

Suppose $q$ does not divide row 1 of $MY$ and
$$(M'G\ N'\,^tG^{-1})=\pmatrix q\\&1\endpmatrix (M\ N+MY/q)$$  
is a coprime pair.  So $(M'\ N')$
is $q^2$-type $0,2$ and $\chi_{\sigma}(M',N')=\chi_{\stufe/q}(q)\chi_{\sigma}(M,N)$.
Reversing, take $(M'\ N')$ of $q^2$-type $0,2,\nu'$.
Take $d$ so that $\det M'=q^2d$.  Write 
$$EM'G=q\pmatrix m_1&m_2\\m_3&m_4\endpmatrix,\ EN'\,^tG^{-1}=\pmatrix n_1&n_2\\n_3&n_4\endpmatrix.$$

Suppose first that $q|m_2$; for each choice of $E$, there is one choice of $G$ so that this
is the case.  
To have $N$ integral we need $(n_1\ n_2)\equiv m_1(y_1\ y_2)\ (q)$.  Then
$$\psi_1(M,N)=\left(\frac{m_1(n_4-m_3y_2-m_4y_4)}{q}\right)
=\left(\frac{m_1n_4-m_3n_2-dy_4}{q}\right),$$
so we have $(q-1)/2$ choices for $y_4$ so that $(M\ N)$ is $q^2$-type $1,2,\nu$.
Summing over corresponding $Y$, we have
$$\sum_Y\left(\frac{y_1}{q}\right)\psi_1(M,N)=\nu\left(\frac{d}{q}\right)\frac{(q-1)}{2}
\left(\frac{m_1n_1}{q}\right).$$
We know $V\simeq\frac{1}{q}M'\,^tN'$ is regular, thus as we let $E$ vary we get
$$\sum_E\left(\frac{m_1n_1}{q}\right)=0.$$

Now suppose $q\nmid m_2$; for each $G$, we have $q$ choices of $E$ so that this is the case.
Choose $Y$ so that $N$ is integral; hence $y_1$ is free and
$$\psi_1(M,N)=\left(\frac{-m_2(n_3-m_3y_1-m_4y_2)}{q}\right)
=\left(\frac{m_4n_1-m_2n_3-dy_1}{q}\right).$$
Restricting our choice of $y_1$ so that $(M\ N)$ is $q^2$-type $1,2,\nu$, and
then summing over corresponding $Y$, we have
$$\sum_Y\left(\frac{y_1}{q}\right)
=\cases \epsilon\nu(q-1)/2&\text{if $q|m_4n_1-m_2n_3$,}\\
-\epsilon\nu\left(1+\epsilon\nu\left(\frac{d(m_4n_1-m_2n_3)}{q}\right)\right)/2
&\text{if $q\nmid m_4n_1-m_2n_3$.}\endcases$$
With $M_0=\frac{1}{q}M'$, take $V\simeq\overline M_0N'$.  Thus with
$(x_1'\ x_2')=(x_1\ x_2)\,^tG^{-1}$, $Q(x_1')=\overline d(m_4n_1-m_2n_3)$.
When $\epsilon=\nu'$, we have $V\simeq\H$ and then
$$\sum_{G,E,Y}\left(\frac{y_1}{q}\right)
=q\cdot 2\epsilon\nu(q-1)/2 - q\cdot\epsilon\nu(q-1)/2
=\nu\nu' q(q-1)/2.$$
When $\epsilon\not=\nu'$, we have $V\simeq\A$ and then
$$\sum_{G,E,Y}\left(\frac{y_1}{q}\right)
=-q\cdot \epsilon\nu(q+1)/2
=\nu\nu' q(q-\epsilon\nu')/2.$$
Thus the contribution is
$$\chi_{\stufe/q}(q)q^{k+1}(q-\epsilon)/2\cdot\E_{0,2,+}+q^{k+1}(q+\epsilon)/2\cdot\E_{0,2,-}.$$

Suppose $q$ does not divide row 1 of $MY$ and
$\rank_q\pmatrix q\\&1\endpmatrix(N+MY/q)=1$.  Then 
we can adjust the equivalence class representative to assume $q$ divides
row 2 of $N+MY/q$;  set
$$(M'G\ N'\,^tG^{-1})=\pmatrix q\\&1/q\endpmatrix(M\ N+MY/q).$$
So $M'$ is $q^2$-type $1,2$  and
$\chi_{\sigma}(M',N')=\chi_{\sigma}(M,N)$.
Reversing, take $(M'\ N')$ of $q^2$-type $1,2,\nu'$.  Then
$$(M\ N)
=\pmatrix 1/q\\&q\endpmatrix E(M'G\ \ N'\,^tG^{-1}-M'GY/q)$$
where $E=\pmatrix 1&q\alpha\\0&1\endpmatrix$
with $\alpha$ varying modulo $q$.

Let $M_0=\pmatrix 1/q\\&1\endpmatrix M'$, $d=\det M_0$.  
Take 
$$Y\equiv G^{-1}\overline M_0\pmatrix 1\\&0\endpmatrix N'\,^tG^{-1}
-G^{-1}\overline M_0\pmatrix 0&0\\0&v\endpmatrix\,^tG^{-1}\ (q).$$
For $(q-1)/2$ choices of $v$ we have
$(M\ N)$ of $q^2$-type $1,2,\nu$.  With $V=\F x_1\oplus\F x_2\simeq GY\,^tG$, $V$ is regular, and
with $(x_1'\ x_2')=(x_1\ x_2)\,^tG^{-1}$, we have $Q(x_1')=y_1$.  So letting
$Y$ vary with $G$, we have
$$\sum_G\psi_1(M,N)\left(\frac{y_1}{q}\right)
=\nu\left(\frac{d}{q}\right)\sum_{\F x_1'}\left(\frac{Q(x_1')}{q}\right)=0$$
since $\F x_1'$ varies over all lines in $V$ as $G$ varies.
Thus there is no contribution from this case. $\square$
\enddemo

\proclaim{Proposition 5.24} We have $\E_{0,0}|A^*_1(q)=0.$ 
\endproclaim

\demo{Proof} With $(M\ N)$ of $q^2$-type $0,0$
$$(M'G\ N'\,^tG^{-1})=\left(M\ \ N+MY/q\right)$$
is a coprime pair of $q^2$-type $0,0$ and $\chi_{\sigma}(M',N')=\chi_{\sigma}(M,N)$.
Reversing, we have $\det N'\equiv\det N\ (q)$ for all $Y$, so
So $\det N'\equiv\det N\ (q),$ and 
$$\sum_Y\left(\frac{y_1}{q}\right)\left(\frac{\det N}{q}\right)=0,$$
proving the proposition. $\square$
\enddemo

\proclaim{Proposition 5.25} For $\chi_1\not=1$, $\nu=+$ or $-$, 
$$\E_{0,1,\nu}|A^*_1(q)
=\epsilon\chi_{\stufe/q}(q) q^{2-k}(q-1)/2\cdot\E_{1,1}-\nu q^2\E_{0,1,\epsilon\nu}.$$
\endproclaim

\demo{Proof}  Take $(M\ N)$ of $q^2$-type $0,1$; we can assume
$q^2$ divides row 2 of $M$.

Say $(M'G\ N'\,^tG^{-1})=(M\ N+MY/q)$ is a coprime pair; then
$\chi_{\sigma}(M',N')=\chi_{\sigma}(M,N)$.
Reversing, assume $(M'\ N')$ is $q^2$-type $0,1,\nu'$ with $q^2$ dividing row 2 of $M'$.
$$(M\ N)=(M'G\ \ N'\,^tG^{-1}-M'GY/q).$$
Write $M'G=q\pmatrix m_1&m_2\\qm_3&qm_4\endpmatrix$,
$N'\,^tG^{-1}=\pmatrix n_1&n_2\\n_3&n_4\endpmatrix$; to have $(M,N)=1$, we need
$q\nmid \det N$.
First suppose $q|n_3$; this is the case for 1 choice of $G$.  Then $q\nmid n_1n_4$, and
by symmetry, $q|m_2$.
To ensure $(M\ N)$ is $q^2$-type $0,1,\nu$, we need that
for some $u\not\equiv0\ (q)$,
$$y_1\equiv\cases \overline m_1n_1-u^2&\text{if $\nu=+$,}\\
\overline m_1n_1-\omega u^2&\text{if $\nu=-$.}\endcases$$
Summing over the corresponding $Y$, we have
$$\sum_Y\left(\frac{y_1}{q}\right)\left(\frac{\det N}{q}\right)
=\nu\left(\frac{m_1n_4}{q}\right)\sum_Y\left(\frac{y_1}{q}\right)$$
and
$$\sum_Y\left(\frac{y_1}{q}\right)=
\cases -\epsilon\nu q^2&\text{if $\epsilon\nu'=\nu$,}\\
0&\text{otherwise.}\endcases$$
So the contribution from this case is
$-\nu q^2\E_{0,1,\epsilon\nu}.$
Now say $q\nmid n_3$; there are $q$ choices of $G$ so that this is the case,
and by symmetry, $q\nmid m_2$.  To have $q\nmid \det N$ with $(M\ N)$ $q^2$-type $0,1,\nu$,
we can choose $y_1,y_2$ freely, leaving us with $(q-1)/2$ choices for $y_4$.
Note that $\psi_1(M/q,N)=\psi_1(M'/q,N').$ Thus,
summing over corresponding $Y$, we have
$$\sum_Y\left(\frac{y_1}{q}\right)\left(\frac{\det N}{q}\right)
=\nu\psi_1(M',N')\sum_Y\left(\frac{y_1}{q}\right),
\ \sum_Y\left(\frac{y_1}{q}\right)=0.$$
Thus the contribution from these terms is 0.

Suppose $\rank_q(N+MY/q)=1$.  
We can adjust the equivalence class representative so that $q$ divides row 1 of
$N+MY/q$, and then
$$(M'G\ N'\,^tG^{-1})=\pmatrix 1/q\\&1\endpmatrix
(M\ N+MY/q)$$
is an integral pair of $q^2$-type $1,1$.
To have $M$ of $q^2$-type $0,1$, we need $q$ to divide row 2 of $EM'$,
thus we only need to consider $E=I$.  
For $q$ choices of $G$, $y_1$ can be chosen freely; consequently 
$\psi_1(M/q,N)=\psi_1(M',N')$ and
$$\sum_Y\left(\frac{y_1}{q}\right)\left(\frac{\det N}{q}\right)
=\nu\psi_1(M',N')\sum_Y\left(\frac{y_1}{q}\right)=0.$$
For 1 choice of $G$, to have $(M\ N)$ of $q^2$-type $0,1,\nu$ we need
$\epsilon\nu\left(\frac{y_1}{q}\right)=1$; then
$$\sum_Y\left(\frac{y_1}{q}\right)\left(\frac{\det N}{q}\right)=
\epsilon q(q-1)/2\cdot\psi_1(M',N').$$
So the contribution from these terms is
$\epsilon\chi_{\stufe/q}(q) q^{2-k}(q-1)/2\cdot\E_{1,1}.$
$\square$
\enddemo

For the next proposition, we need the following lemma.

\proclaim{Lemma 5.26} The $GL_2(\F)$-orbit of a matrix for $\H$ has
$q(q^2-1)/2$ elements; for any $a\in\F^{\times}$, there are $q(q-1)/2$
elements in this orbit of the form $\pmatrix a&*\\*&*\endpmatrix$, and
$q(q-1)$ elements of the form $\pmatrix 0&*\\*&*\endpmatrix$.  The
$GL_2(\F)$-orbit of a matrix for $\A$ has $q(q-1)^2/2$ elements; for
any $a\in\F^{\times}$, there are $q(q-1)/2$ elements in this orbit of
the form $\pmatrix a&*\\*&*\endpmatrix$, and none of the form
$\pmatrix 0&*\\*&*\endpmatrix.$
\endproclaim

\demo{Proof}  The number of elements in the orbit of $\H$ is
$$\frac{|GL_2(\F)|}{o(\H)}=\frac{(q^2-1)(q^2-q)}{2(q-1)}=\frac{q(q^2-1)}{2}.$$
Choose $a\in\F^{\times}$, and let $V$ be a hyperbolic plane.  There
are $q-1$ vectors $x_1$ in $V$ so that $Q(x_1)=a$.  Thus the number of
ways to choose a basis $(x_1\ x_2)$ for $V$ so that $Q$ is given by
$\pmatrix a&*\\*&*\endpmatrix$ relative to this basis is
$(q-1)(q^2-q)$.
Thus the number of matrices in the orbit of $\H$ of the form $\pmatrix
a&*\\*&*\endpmatrix$ is
$$\frac{(q-1)(q^2-q)}{o(\H)}=q(q-1)/2.$$
Hence the number of matrices in the orbit of $\H$ of the form
$\pmatrix 0&*\\*&*\endpmatrix$ is
$$\frac{q(q^2-1)}{2}-(q-1)\frac{q(q-1)}{2}=q(q-1).$$

The argument for the orbit of $\A$ is virtually identical, but here
$o(\A)=2(q+1)$, for for $V$ an anisotropic plane and
$a\in\F^{\times}$, there are $q+1$ vectors $x_1$ in $V$ so that
$Q(x_1)=a$.
$\square$
\enddemo

\proclaim{Proposition 5.27} For $\chi_q\not=1$, $\nu=+$ or $-$,
$$\E_{0,2,\nu}|A^*_1(q)=
\epsilon\chi_{\stufe/q}(q) q^{2-k}(q-1)/2\cdot\E_{1,2}.$$
\endproclaim

\demo{Proof} 
Suppose $\rank_q(M+MY/q)=0$.  Set
$$(M'G\ N'\,^tG^{-1})=\frac{1}{q}(M\ \ N+MY/q);$$
so $M'$ is $q$-type $2$ and $\chi_{\sigma}(M',N')=\chi_{\sigma}(M,N)$.
Reversing, take $(M'\ N')$ of $q$-type 2.
To have $(M\ N)$ $q^2$-type $0,2,\nu$, we need $\nu\left(\frac{\det Y}{q}\right)=1.$
Thus $Y$ varies over the orbit of $I$ if $\nu=+$ and over the orbit of
$\pmatrix 1\\&\omega\endpmatrix$ if $\nu=-$.  Summing over these $Y$ and using Lemma 3.13, we have
$$\sum_Y\left(\frac{y_1}{q}\right)\left(\frac{\det N}{q}\right)
=\nu\left(\frac{\det M'}{q}\right)\frac{q(q-1)}{2}\sum_{y_1\in\F^{\times}}\left(\frac{y_1}{q}\right)=0.$$

Suppose $\rank_q(N+MY/q)=1$.  So adjusting $(M\ N)$, we can assume $q$ divides
row 1 of $N+MY/q$.  Then
$$(M'G\ N'\,^tG^{-1})=\pmatrix 1/q\\&1\endpmatrix (M\ \ N+MY/q)$$
is $q^2$-type $1,2$, and $\chi_{\sigma}(M',N')=\chi_{\stufe/q}(q)\chi_{\sigma}(M,N)$.
Reversing, take $(M'\ N')$ of $q$-type 2.
We can assume $q$ divides row 2 of $M'$.  To have $M$ of $q^2$-type $0,2$ we need 
only consider $E=I$.
Let $M_0=\pmatrix 1\\&1/q\endpmatrix M'$ and set $d=\det M_0$; so $M_0$ is $q$-type 2, and 
$$G^{-1}\overline M_0N\equiv G^{-1}\overline M_0\pmatrix 0\\&1\endpmatrix N'\,^tG^{-1}-Y\ (q)$$
is symmetric modulo $q$.
To have $(M\ N)$ of $q^2$-type $0,2,\nu$,
we need $\nu\left(\frac{\det\overline M_0N}{q}\right)=1.$  So we vary $Y$ so that $G^{-1}\overline M_0N$
varies over the orbit of $I$ if $\nu=+$, and over the orbit of $\pmatrix 1\\&\omega\endpmatrix$ if
$\nu=-$.   Thus
$$\sum_Y\left(\frac{y_1}{q}\right)\psi_1(M,N)=\nu\left(\frac{d}{q}\right)\sum_Y\left(\frac{y_1}{q}\right).$$
Using Lemma 3.13, we find that for $q$ choices of $G$ we have
$$
\sum_Y\left(\frac{y_1}{q}\right)=
\epsilon\nu q(q-1)/2\cdot\left(\frac{d}{q}\right)\psi_1(M',N'),$$
and for 1 choice of $G$ the sum on $Y$ is 0.
So the contribution from these terms is $\epsilon\chi_{\stufe/q}(q) q^{2-k}(q-1)/2\cdot\E_{1,2}.$

Say $\rank_q(N+MY/q)=2$; then
$$(M'G\ N'\,^tG^{-1})=(M\ N+MY/q)$$ is $q^2$-type $0,2$
and $\chi_{\sigma}(M',N')=\chi_{\sigma}(M,N)$.
Reversing, take $(M'\ N')$ of
$$(M\ N)=(M'G\ \ N'\,^tG^{-1}-M'GY/q)$$
$q^2$-type $0,2$.  To have $(M\ N)$ $q^2$-type $0,2,\nu$, we need
$\nu\left(\frac{\det NM_0}{q}\right)=1$ where $M_0=\frac{1}{q}M'$.  
We choose those $Y$ so that $G^{-1}\overline M_0N$ is in the $GL_2(\F)$-orbit of $I$ if $\nu=+$,
and in the orbit of $\pmatrix 1\\&\omega\endpmatrix$ if $\nu=-$.  Then with $d=\det M_0$,
and again using Lemma 3.13,
$$\align
&\sum_Y\left(\frac{y_1}{q}\right)\left(\frac{\det N}{q}\right)\\
&\quad = \nu\left(\frac{d}{q}\right)\frac{q(q-1)}{2}\sum_{a\in\F^{\times}}\left(\frac{a}{q}\right)
+\epsilon\left(\frac{d}{q}\right)\frac{q(q-1)}{2}\left(\frac{Q(x_1')}{q}\right)
\endalign$$
where $V=\F x_1\oplus\F x_2\simeq\overline M_0N'$, and $(x_1'\ x_2')=(x_1\ x_2)\,^tG^{-1}$.
As $G$ varies, $\F x_1'$ varies over all lines in $V$ and since $V$ is regular,
$$\sum_{x_1'}\left(\frac{Q(x_1')}{q}\right)=0.$$
Thus the contribution from these terms is 0.  $\square$
\enddemo

\proclaim{Proposition 5.28} For $\chi_1\not=1$, we have
$\E_2|A^*_2(q)=
\epsilon q^{2k}\E_{0,2}.$
\endproclaim

\demo{Proof}  Take $(M\ N)$ of $q^2$-type 2, $Y\in\Z^{2,2}_{\sym}$
with $\rank_qY=2$.  Then $\rank_qMY=2$, so setting
$$(M'\ N')=q(M\ N+MY/q),$$
we have $(M'\ N')$ of $q^2$-type $0,2$
and $\chi_{\sigma}(M',N')=\chi_{\sigma}(M,N)$.
Reversing, take $(M'\ N')$ of $q^2$-type $0,2$.  
To have $N$ integral, we need
to choose $\frac{1}{q}M'Y\equiv N'\ (q)$; then
$$\left(\frac{\det Y}{q}\right)\left(\frac{\det
M}{q}\right)=\left(\frac{\det N'}{q}\right),$$
and $\E_2|A^*_2(q)=q^{2k}\E_{0,2}.$ 
$\square$
\enddemo

\proclaim{Proposition 5.29} We have $\E_{1,1}|A^*_2(q)=q(q-1)\E_{1,1}.$
\endproclaim

\demo{Proof}  
Take $(M\ N)$ of $q^2$-type $1,1$ with $q^2$ dividing row 2 of $M$.
Hence $\rank_qMY=1$.

Suppose first
$$(M'\ N')=\pmatrix q\\&1\endpmatrix \left(M\
N+\frac{1}{q}MY\right)$$
is a coprime pair.
Then $(M'\ N')$ is of $q^2$-type $0,1$,
$\chi_{\sigma}(M',N')=\chi_{\stufe/q}(q)\chi_{\sigma}(M,N)$.
Reversing, take $(M'\ N')$ of $q^2$-type $0,1$.  For 1 choice of $G$,
we have $q$ choices of $Y$ so that $N$ is integral; here $y_4$ is free,
and consequently
$$\sum_Y\left(\frac{\det Y}{q}\right)\psi_1(M,N)=0.$$
For each of the other $q$ choices of $G$, there are $q$ choices of $Y$ so that
$N$ is integral; here $y_1$ is free and again
$\sum_Y\left(\frac{\det Y}{q}\right)\psi_1(M,N)=0.$

Say $\rank_q\pmatrix q\\&1\endpmatrix \left(N+\frac{1}{q}MY\right)=1;$
we can adjust the equivalence class representative so that 
$$(M'\ N')=\pmatrix q\\&1/q\endpmatrix (M\ N+MY/q)$$
is integral of $q^2$-type $1,1$ and $\chi_{\sigma}(M',N')=\chi_{\sigma}(M,N)$.
Reversing, take $(M'\ N')$ of $q^2$-type $1,1$; we can assume $q$ does not divide row 2 of $M'$
and $q$ exactly divides row 1 of $M'$.  
For $E\in SL_2(\Z)$, we have $\pmatrix 1/q\\&q\endpmatrix
E\pmatrix q\\&1/q\endpmatrix\in SL_2(\Z)$ if and only if
$E\equiv\pmatrix *&0\\*&*\endpmatrix\ (q^2),$  
and to have $(M\ N)$ of $q^2$-type $1,1$ we need $q$ exactly dividing row 1 of $EM'$.
Thus we need to
consider $E(M'\ N')$ where $E=\pmatrix 1&qb\\&1\endpmatrix$, $b$
varying through the $q-1$ congruence classes modulo $q$ so that $q^2$ 
does not divide row 1 of $\pmatrix 1&qb\\&1\endpmatrix M'$.
For each such $E$, there are $q$ choices for $Y$ so that $N$ is integral, and then
$$\sum_Y\left(\frac{\det Y}{q}\right)\psi_1(M,N)=\epsilon q\psi_1(M',N'),
\ \sum_Y\left(\frac{\det Y}{q}\right)=\epsilon q.$$
So the contribution from these terms is
$q(q-1)\E_{1,1}.$ $\square$
\enddemo

\proclaim{Proposition 5.30} For $\chi_1\not=1$, $\nu=+$ or $-$,
$$\align
\E_{1,2,\nu}|A_2^*(q)&=
\nu\chi_{\stufe/q}(q) q^k(q-1)/2\cdot\E_{0,2,\epsilon}\\
&\quad -\nu\chi_{\stufe/q}(q) q^k(q+1)/2\cdot\E_{0,2,-\epsilon}\\
&\quad +\epsilon q(q-1)/2\cdot\E_{1,2}.
\endalign$$
\endproclaim

\demo{Proof}  Choose $(M\ N)$ of $q^2$-type $1,2,\nu$ with $q$ dividing
row 2 of $M$.

Suppose $$(M'\ N')=\pmatrix q\\&1\endpmatrix(M\
N+\frac{1}{q}MY)$$ 
is a coprime pair.  So $(M'\ N')$ has $q^2$-type $0,2$
and $\chi_{\sigma}(M',N')=\chi_{\stufe/q}(q)\chi_{\sigma}(M,N)$.
Reversing, 
take $(M'\ N')$ of $q^2$-type $0,2,\nu'$; so with $M_0=\frac{1}{q}M'$,
$d=\det M_0$, we have $\nu'\left(\frac{d\cdot \det N'}{q}\right)=1$.
To have $N$ integral with $(M,N)$ of $q^2$-type $1,2,\nu$, we need
$$EM_0Y\,^tM_0\,^tE\equiv EN'\,^tM_0\,^tE-\pmatrix 0&0\\0&t\endpmatrix\ (q)$$
where  $\nu\left(\frac{t}{q}\right)=1$.
Then for fixed $E$,
$$\sum_Y\left(\frac{\det Y}{q}\right)\psi_1(M,N)
=\nu\left(\frac{d}{q}\right)\sum_Y\left(\frac{\det Y}{q}\right)$$
and writing $EN'\,^tM_0\,^tE=\pmatrix a&b\\b&c\endpmatrix$,
$$\align
\sum_Y\left(\frac{\det Y}{q}\right)
&=\sum_t\left(\frac{-at+\det N'M_0}{q}\right)\\
&={\cases (q-1)/2\cdot\left(\frac{d\cdot\det N'}{q}\right)&\text{if $q|a$,}\\
-\left(\left(\frac{d\cdot\det N'}{q}\right)+\nu\left(\frac{-a}{q}\right)\right)/2
&\text{if $q\nmid a$.}\endcases}
\endalign$$
With $V=\F x_1\oplus\F x_2\simeq N'\,^tM_0$ and
$(x_1'\ x_2')=(x_1\ x_2)\,^tE$, $Q(x_1')=a$.  
Then letting $E$ vary,
$$\align
\sum_{E,Y}\left(\frac{\det Y}{q}\right)
&={\cases (q-1)/2\cdot\left(\frac{d\cdot\det N'}{q}\right)&\text{if $\nu'=\epsilon$,}\\
-(q+1)/2\cdot\left(\frac{d\cdot \det N'}{q}\right)&\text{if $\nu'\not=\epsilon$}
\endcases}\\
&\epsilon(q-1)/2.
\endalign$$
So the contribution from these terms is
$$\nu\chi_{\stufe/q}(q)q^k(q-1)/2\cdot\E_{0,2,\epsilon}-\nu\chi_{\stufe/q}(q)q^k(q+1)/2\cdot\E_{0,2,-\epsilon}.$$

Now suppose $\rank_q\pmatrix q\\&1\endpmatrix
(N+\frac{1}{q}MY)=1.$  
Adjust the equivalence class representative 
so that $q$ divides row 2 of $N+MY/q$;
$$(M'\ N')=\pmatrix q\\&1/q\endpmatrix(M\ N+\frac{1}{q}MY)$$
is coprime with $q^2$-type $1,2$, and
$\chi_{\sigma}(M',N')=\chi_{\sigma}(M,N)$.
Reversing, take $(M'\ N')$ $q^2$-type $1,2,\nu'$ with $q$ dividing row 1 of $M'$; set
$$(M\ N)=\pmatrix 1/q\\&q\endpmatrix E(M'\ \ N'-M'Y/q).$$
To have $M$ of $q^2$-type $1,2$ we need $q$ dividing row 1 of $EM'$; so we take
$E=\pmatrix 1&qb\\&1\endpmatrix$ with $b$ varying
modulo $q$.
Let $M_0=\pmatrix 1/q\\&1\endpmatrix EM'$; so $M_0$ is integral with $\rank_qM_0=2.$
To have $N$ integral with $(M\ N)$ of $q^2$-type $1,2,\nu$, we need
$$M_0Y\,^tM_0\equiv
M_0\,^tN'\,^tE\pmatrix 1\\&q\endpmatrix-\pmatrix 0&0\\0&t\endpmatrix\ (q)$$
where  $\nu\left(\frac{t}{q}\right)=1$.  Then
$$\sum_Y\left(\frac{\det Y}{q}\right)\psi_1(M,N)
=\nu\left(\frac{d}{q}\right)\sum_Y\left(\frac{\det Y}{q}\right),$$
and with $EN'=\pmatrix n_1&n_2\\n_3&n_4\endpmatrix$, $EM'=\pmatrix qm_1&qm_2\\m_3&m_4\endpmatrix$,
we have $M_0Y\equiv\pmatrix n_1&n_2\\0&0\endpmatrix-\overline d\pmatrix0&0\\-m_2t&m_1t\endpmatrix\ (q).$
Thus $\det M_0Y\,^tM_0\equiv-t(m_1n_1+m_2n_2)\ (q).$  Since
$$
\psi_1(M',N')
=\left(\frac{d}{q}\right)\psi_1(M'\overline M_0,N'\,^tM_0)
=\left(\frac{d}{q}\right)\left(\frac{m_1n_1+m_2n_2}{q}\right),
$$
we have
$$
\sum_Y\left(\frac{\det Y}{q}\right)
=\epsilon\nu\left(\frac{d}{q}\right)\psi_1(M',N')(q-1)/2.
$$
Thus the contribution from these terms is
$\epsilon q(q-1)/2\cdot\E_{1,2}.$
$\square$
\enddemo

\proclaim{Proposition 5.31} We have
$E_{0,0}|A^*_2(q)=\epsilon q(q-1)\E_{0,0}$.
\endproclaim

\demo{Proof}  With $(M\ N)$ of $q^2$-type $0,0$, set $(M'\ N')=(M\
N+MY/q).$  Reversing, we can choose any $Y$, and we have
$N$ integral with $N'\equiv N\ (q)$.  Thus
$$\sum_Y\left(\frac{\det Y}{q}\right)\left(\frac{\det N}{q}\right)
=\left(\frac{\det N'}{q}\right)\epsilon q(q-1),$$
hence the claim follows. $\square$
\enddemo

\proclaim{Proposition 5.32}  For $\nu=+$ or $-$,
$\E_{0,1,\nu}|A_2^*(q)=\epsilon q(q-1)\E_{0,1,\nu}.$
\endproclaim

\demo{Proof}  Take $(M\ N)$ of $q^2$-type $0,1,\nu$ with $q^2$ dividing
row 2 of $M$.  Then $q$ does not divide row 2 of $N+\frac{1}{q}MY$, so
$\rank_q(N+\frac{1}{q}MY)\ge 1$.

Say $(M'\ N')=(M\ N+\frac{1}{q}MY)$ is a coprime pair. 
Reversing, take $(M'\ N')$ of $q^2$-type $0,1,\nu'$.  Then there are $q^2(q-1)/2$
choices of $Y$ so that $(M\ N)$ is $q^2$-type $0,1,\nu$, and
$$\align
\sum_Y\left(\frac{\det Y}{q}\right)\left(\frac{\det N}{q}\right)
&=\nu\psi_1(M'/q,N')\sum_Y\left(\frac{\det Y}{q}\right)\\
&={\cases \nu\psi_1(M'/q,N')\cdot\epsilon q(q-1)&\text{if $\nu'=\nu$,}\\
0&\text{otherwise.}\endcases}
\endalign$$
So the contribution from these terms is
$\epsilon q(q-1)\E_{0,1,\nu}$.

Now say $\rank_q(N+\frac{1}{q}MY)=1.$  Thus we can adjust the equivalence
class representative so that 
$$(M'\ N')=\pmatrix 1/q\\&1\endpmatrix (M\ N+MY/q)$$
is a coprime pair of  $q^2$-type $1,1$, and
$\chi_{\sigma}(M',N')=\chi_{\stufe/q}(q)\chi_{\sigma}(M,N)$.
Reversing, take $(M'\ N')$ of $q^2$-type $1,1$.  For any $Y$ with $q\nmid\det Y$,
$(M\ N)$ is integral and
coprime.  To have $(M\ N)$ of $q^2$-type $0,1$, we need 
$E\equiv \pmatrix 1&\alpha\\0&1\endpmatrix\ (q).$
Take $G\in SL_2(\Z)$ so that $M'G\equiv\pmatrix
m_1&0\\0&0\endpmatrix\ (q).$  Thus $N'\,^tG^{-1}\equiv
\pmatrix n_1&n_2\\0&n_4\endpmatrix\ (q),$ $q\nmid n_4$, and hence with 
$U=G^{-1}Y\,^tG^{-1}=\pmatrix u_1&u_2\\u_2&u_4\endpmatrix,$
$$N\,^tG^{-1}=\pmatrix q\\&1\endpmatrix(N'\,^tG^{-1}-\frac{1}{q}M'GU)
\equiv\pmatrix -m_1u_1&-m_1u_2\\0&n_4\endpmatrix\ (q).$$
To have $(M\ N)$ of $q^2$-type $0,1,\nu$, we need to choose $u_1$ so that
$\nu\left(\frac{-u_1}{q}\right)=1$.  So $u_1$ is fixed with $q\nmid u_1$ and
$u_2,u_4$ are unconstrained; thus
$$\sum_Y\left(\frac{\det Y}{q}\right)\left(\frac{\det N}{q}\right)
=\nu\psi_1(M/q,N)\sum_U\left(\frac{\det U}{q}\right)=0.$$
The proposition now follows.
$\square$
\enddemo

\proclaim{Proposition 5.33} For $\chi_1\not=1$, $\nu=+$ or $-$, we have 
$$\align
\E_{0,2,\nu}|A^*_2(q)&=
\epsilon q(q-2)\E_{0,2,\nu}
+\epsilon\nu\chi_{\stufe/q}(q) q^{1-k}(q-1)/2\cdot(\E_{1,2,+}-\E_{1,2,-})\\
&\quad +q^{1-2k}(q-1)(q+\epsilon\nu)/2\cdot\E_2.
\endalign$$
\endproclaim

\demo{Proof}  Begin with $(M\ N)$ of $q^2$-type $0,2,\nu$.

Say  $(M'\ N')=(M\ N+MY/q)$ is coprime; then $(M'\ N')$
 has $q^2$-type $0,2$ and $\chi_{\sigma}(M',N')=\chi_{\sigma}(M,N)$.
Reversing,
take $(M'\ N')$ of $q^2$-type $0,2,\nu'$.  Choose $Y\in\F^{2,2}_{\sym}$ so that
$$(M\ N)=(M'\ N'-M'Y/q)$$
is $q^2$-type $0,2,\nu$.  Set $M_0=\frac{1}{q}M'$; thus we need $Y$ so that
$\nu\left(\frac{\det\overline M_0N'-Y}{q}\right)=1.$  
We know there is some $G\in SL_2(\Z)$ so that $^tG\overline M_0N'G\equiv\pmatrix 1&0\\0&c\endpmatrix\ (q)$
where $\nu'\left(\frac{c}{q}\right)=1$, so choosing $Y$ is equivalent to choosing $U=\,^tGYG$ so that
$\det\pmatrix 1&0\\0&c\endpmatrix-U\equiv t\ (q)$ where $\nu\left(\frac{t}{q}\right)=1$.
So we need either $u_1\not\equiv1\ (q),$ $u_4\equiv(t+u_2^2)\overline {(u_1-1)}+c\ (q)$,
or $u_1\equiv 1\ (q)$, $-u_2^2\equiv t\ (q)$.  Summing over all such $U$, we have
$$\align
\sum_U\left(\frac{\det U}{q}\right)
&=\sum_{{u_1\not\equiv1}\atop{u_2,\,t}}\left(\frac{u_1-1}{q}\right)\left(\frac{u_2^2+u_1(t+cu_1-c)}{q}\right)\\
&\quad+
{\cases 0&\text{if $\epsilon\not=\nu$,}\\
\sum_{u_2\not\equiv0,\,u_4}\left(\frac{u_4-u_2^2}{q}\right)&\text{if $\epsilon=\nu$,}\endcases}
\endalign$$
where $u_2,u_4$ vary modulo $q$, and $t$ varies over all quadratic residues modulo $q$ if
$\nu=+$, and over all quadratic non-residues if $\nu=-$.  Note that for any $u_2$,
$\sum_{u_4}\left(\frac{u_4-u_2^2}{q}\right)=0$.  
For fixed $t$,
$$\align
\sum_{u_1,u_2}\left(\frac{\det U}{q}\right)&=\sum_{u_1\not\equiv 1}\left(\frac{u_1-1}{q}\right)\left(\frac{u_1(t+u_1c-c)}{q}\right)\\
&\quad +(q-1)\sum_{u_1\equiv0,1-\overline ct}\left(\frac{u_1-1}{q}\right)\\
&\quad - \sum_{u_1\not\equiv0,1,1-\overline ct}\left(\frac{u_1-1}{q}\right)\left(1+\left(\frac{u_1(t+u_1c-c)}{q}\right)\right).
\endalign$$
Since $\left(\frac{u_1(t+u_1c-c)}{q}\right)=0$ when $u_1\equiv 0$ or $1-\overline ct\ (q)$, and since
$\sum_{u_1}\left(\frac{u_1-1}{q}\right)=0$, we have
$$\align
\sum_{u_1,u_2}\left(\frac{\det U}{q}\right)
&=(q-1)\sum_{u_1\equiv0,1-\overline ct}\left(\frac{u_1-1}{q}\right)\\
&\quad -\sum_{u_1\not\equiv0,1-\overline ct}\left(\frac{u_1-1}{q}\right)
+\sum_{u_1}\left(\frac{u_1-1}{q}\right)\\
&=q\sum_{u_1\equiv0,1-\overline ct}\left(\frac{u_1-1}{q}\right)\\
&={\cases \epsilon q&\text{if $c\equiv t\ (q)$,}\\
\epsilon q\left(1+\left(\frac{ct}{q}\right)\right)&\text{if $c\not\equiv t\ (q)$.}\endcases}
\endalign$$
We know $\nu\nu'\left(\frac{ct}{q}\right)=1$, so
$$\sum_{t,u_1,u_2}\left(\frac{\det U}{q}\right)=\cases 0&\text{if $\nu\not=\nu'$,}\\
\epsilon q+2\epsilon q(q-3)/2&\text{if $\nu=\nu'$.}
\endcases$$
Hence when $\nu=\nu'$,
$$\sum_Y\left(\frac{\det Y}{q}\right)\left(\frac{\det N}{q}\right)
=\epsilon \nu\left(\frac{d}{q}\right)q(q-2)=\epsilon\left(\frac{\det N'}{q}\right)q(q-2).$$
Thus the contribution from these terms is $\epsilon q(q-2)\E_{0,2,\nu}$.

Now say $\rank_q(N+\frac{1}{q}MY)=1.$  Adjust the
equivalence class representative for $(M\ N)$ so that $q$ divides row
1 of $N+\frac{1}{q}MY$, and set 
$$(M'\ N')=\pmatrix 1/q\\&1\endpmatrix(M\ N+\frac{1}{q}MY).$$
So $(M'\ N')$ has $q^2$-type $1,2$ with $q$ dividing row 2 of $M'$
and $\chi_{\sigma}(M',N')=\chi_{\stufe/q}(q)\chi_{\sigma}(M,N)$.
Reversing, 
take $(M'\ N')$ of $q^2$-type $1,2$ with $q$ dividing row 2 of $M'$, and set
$$(M\ N)=\pmatrix q\\&1\endpmatrix E(M'\ N'-M'Y/q),$$
$E\in\,^t\G_1$.
To have $M$ of $q^2$-type $0,2$, we need $q$ to divide row 2 of $EM'$; then
$$\pmatrix q\\&1\endpmatrix E\pmatrix 1/q\\&1\endpmatrix\in GL_2(\Z).$$
So we only need to consider $E=I$.
Set $M_0=\pmatrix 1\\&1/q\endpmatrix M'$.
Choose $G\in SL_2(\Z)$ so that $M_0G\equiv \pmatrix d\\&1\endpmatrix \ (q)$; 
by the symmetry of $M'\,^tN'$, we have
$N'\,^tG^{-1}\equiv\pmatrix *&*\\0&c\endpmatrix\ (q)$, and $q\nmid c$ since
$(M',N')=1$.  (So $\psi_1(M',N')=\left(\frac{dc}{q}\right)$.)  With
$U=M_0Y\,^tM_0=\pmatrix u_1&u_2\\u_2&u_4\endpmatrix$, we have
$\left(\frac{\det U}{q}\right)=\left(\frac{\det Y}{q}\right)$ and
$$N\,^tM_0=\pmatrix q\\&1\endpmatrix N'\,^tM_0-U.$$  To have $(M\ N)$ of
$q^2$-type $0,2,\nu$, we need $u_1(u_4-c)-u_2^2\equiv t\ (q)$ for some $t$ so that
$\nu\left(\frac{t}{q}\right)=1$.  For each such $t$ and $u_1\not\equiv 0\ (q)$, we
need $u_4\equiv\overline u_1(t+u_2^2)+c\ (q)$; if $\nu=\epsilon$, we can allow $q|u_1$
(and then $u_4$ is unconstrained).  Letting $t$ vary so that
$\nu\left(\frac{t}{q}\right)=1$,
$$\align
\sum_{{t,u_2}\atop{u_1\not\equiv0}}\left(\frac{t+u_1c}{q}\right)
&={\cases -q\sum_{u_1\not\equiv0}\left(1+\left(\frac{u_1c}{q}\right)\right)/2
&\text{if $\nu=+$,}\\
q\sum_{u_1\not\equiv0}\left(1-\left(\frac{u_1c}{q}\right)\right)/2&\text{if $\nu=-$}\endcases}\\
&=-\nu q(q-1)/2.
\endalign$$
If $\nu=\epsilon$ then we also have a contribution to $\sum_Y\left(\frac{\det Y}{q}\right)$ of
$$\sum_{t,u_4}\sum_{-u_2^2\equiv t}\left(\frac{-u_2^2}{q}\right)=\epsilon q(q-1).$$
Thus the sum on $Y$ is $\epsilon q(q-1)/2$.  Since
$$\sum_Y\left(\frac{\det Y}{q}\right)\left(\frac{\det N}{q}\right)
=\nu\nu'\psi_1(M',N')\sum_Y\left(\frac{\det Y}{q}\right),$$
the contribution from this case is
$\nu\chi_{\stufe/q}(q) q^{1-k}(q-1)/2\cdot(\E_{1,2,\epsilon}-\E_{1,2,-\epsilon}).$

Say $\rank_q(N+\frac{1}{q}MY)=0$.  With $(M'\
N')=\frac{1}{q}(M\ N+MY/q)$, we have $\rank_qM'=2$
and $\chi_{\sigma}(M',N')=\chi_{\sigma}(M,N)$.
Reversing, 
with $(M'\ N')$ of $q$-type 2, set
$$(M\ N)=q(M'\ N'-M'Y/q).$$
So $(M\ N)$ is $q^2$-type $0,2,\nu$ if and only if $\nu\left(\frac{\det Y}{q}\right)=1.$
Summing over such $Y$, we have
$$\sum_Y\left(\frac{\det Y}{q}\right)\left(\frac{\det N}{q}\right)
=\nu\left(\frac{\det M'}{q}\right)
\#\left\{\sym\ Y\ (q):\ \nu\left(\frac{\det Y}{q}\right)=1\ \right\}$$
and
$$\#\left\{\sym\ Y\ (q):\ \nu\left(\frac{\det Y}{q}\right)=1\ \right\}
=\nu q(q-1)(q+\epsilon\nu)/2.$$
Thus the contribution from this case is $\nu q^{1-2k}(q-1)(q+\epsilon\nu)/2\cdot\E_2$
if $\chi_q=1$, and $q^{1-2k}(q-1)(q+\epsilon\nu)/2\cdot\E_2$ if $\chi_q\not=1$. $\square$
\enddemo

\bigskip
\head{\bf \S6. The Main Theorem and its proof}\endhead
\smallskip

\proclaim{Theorem 6.1}  Let
$K$ be a maximal even integral lattice with $\rank K=2k\ge 8$,
odd level $\stufe$, and $\chi(*)=\left(\frac{*}{\disc K}\right)$.  
For $q$ a prime dividing $\stufe$, set
$$c_1(q)=\cases 
\chi_{\stufe/q}(q)\left(\frac{\eta_q}{q}\right)q^{-1}\g(q)&\text{if $\chi_q\not=1$,}\\
-q^{-1}&\text{if $\chi_q=1$,}
\endcases$$
where $\Z_qK\simeq2\big<1,-1,\ldots,1,-1,\eta_q\big>\perp2q\big<\eta'_q\big>$ when $\chi_q\not=1$,
and set
$$c_2(q)=\cases 
\left(\frac{-1}{q}\right)q^{-1}&\text{if $\chi_q\not=1$,}\\
q^{-2}&\text{if $\chi_q=1$.}
\endcases$$
Extending $c_1,c_2$ multiplicatively, we have
$$\theta(\gen K)=\sum_{\stufe_0\stufe_1\stufe_2=\stufe} c_1(\stufe_1)c_2(\stufe_2) \E_{(\stufe_0,\stufe_1,\stufe_2)}$$
where $\E_{(\stufe_0,\stufe_1,\stufe_2)}$ is the Eisenstein series of level $\stufe$, character
$\chi$  defined in \S2.
\endproclaim

\demo{Proof} Using [3], we know that $\theta(\gen K)\in \Eis_k(\stufe,\chi)$, 
the space of degree 2 Siegel Eisenstein series of weight $k$, level $\stufe$, character $\chi$;
we normalised $\theta(\gen K)$ to have 0th Fourier coefficient equal to 1, and we showed
in Theorem 4.11 that for each prime $q|\stufe$, $\theta(\gen K)$ is an eigenform
for $T_K(q)$.  We will show there is a unique $\E\in\Eis_k(\stufe,\chi)$ with 0th Fourier
coefficient 1 so that for each prime $q|\stufe$, $\E|T_K(q)\in\Eis_k(\stufe,\chi)$.

We know $\{\E_{(\stufe_0,\stufe_1,\stufe_2)}:\ \stufe_0\stufe_1\stufe_2=\stufe\ \}$ is a basis for
$\Eis_k(\stufe,\chi)$, and (using [3]) the 0th Fourier coefficient of
$\E_{(\stufe_0,\stufe_1,\stufe_2)}$ is 0 unless $\stufe_0=\stufe$, in which case it is 1.
So set
$$\E=\sum_{\stufe_0\stufe_1\stufe_2=\stufe} c(\stufe_1,\stufe_2)\E_{(\stufe_0,\stufe_1,\stufe_2)}$$
where $c(\stufe_1,\stufe_2)\in\C$ with $c(1,1)=1$, and suppose that for each prime
$q|\stufe$, we have $\E|T_K(q)\in\Eis_k(\stufe,\chi)$.  Fix a prime $q|\stufe$;
rewrite
$$\E=\sum_{Q_1Q_2|\stufe/q}\E(Q_1,Q_2)$$
where
$$\align
\E(Q_1,Q_2)&=c(Q_1,Q_2)\E_{(\stufe/(Q_1Q_2),Q_1,Q_2)}
+ c(qQ_1,Q_2)\E_{(\stufe/(qQ_1Q_2),qQ_1,Q_2)}\\
&\quad + c(Q_1,qQ_2)\E_{(\stufe/(qQ_1Q_2),Q_1,qQ_2)}.
\endalign$$
For $Q_1Q_2|\stufe/q$ and $\sigma=(\stufe/(qQ_1Q_2),Q_1,Q_2)$, we supplement the notation
of \S5 (where $\sigma$ was fixed) by setting
$$\E_0^{\sigma}=\E_{(\stufe/(Q_1Q_2),Q_1,Q_2)},\ 
\E_1^{\sigma}=\E_{(\stufe/(qQ_1Q_2),qQ_1,Q_2)},\ 
\E_2^{\sigma}=\E_{(\stufe/(qQ_1Q_2),Q_1,qQ_2)}.$$
So again supplementing the notation, the results of \S5 tell us how to write
$\E_i^{\sigma}|T_K(q)$ in terms of the elements of
$$\B^{\sigma}=\{\E_{0,0}^{\sigma},\E_{0,1,\nu}^{\sigma},\E_{0,2,\alpha}^{\sigma},
\E_{1,1}^{\sigma},\E_{1,2,\nu}^{\sigma},\E_2^{\sigma}:\ \nu=+\text{ or }-,\ \alpha\in\F_q^{\times}\ \}.$$
So we have $\E_i^{\sigma},\E_i^{\sigma}|T_K(q)\in\spn\B^{\sigma}$, and
with $\sigma$ varying over all partitions of $\stufe/q$,
$\cup_{\sigma}\B^{\sigma}$ is a linearly independent set.
Hence to have $\E|T_K(q)\in\Eis_k(\stufe,\chi)$, for each $\sigma=(\stufe/(qQ_1Q_2),Q_1,Q_2)$
we must have
$$(c(Q_1,Q_2)\E_0^{\sigma}+c(qQ_1,Q_2)\E_1^{\sigma}+c(Q_1,qQ_2)\E_2^{\sigma})|T_K(q)
\in\spn\{\E_0^{\sigma},\E_1^{\sigma},\E_2^{\sigma}\}.$$
(To satisfy this, there are many conditions that must hold, and some of these
are easier to work with than others;
three of these conditions tell us that in the expression for
$\E(Q_1,Q_2)|T_K(q)$
in terms of the elements of $\B^{\sigma}$, the
coefficients on $(\E_{0,1,+}^{\sigma}-\E_{0,1,-}^{\sigma})$ and $(\E_{0,2,+}^{\sigma}-\E_{0,2,-}^{\sigma})$
must be 0, and the coefficients on $\E_{0,0}^{\sigma}$ and $\E_{0,1}^{\sigma}$ must be equal.
After using this criteria to solve for $c(qQ_1,Q_2),c(Q_1,qQ_2)$, it is not difficult
to verify that all the required conditions are met.)
We find that we must have
$$c(qQ_1,Q_2)=c_1(q)c(Q_1,Q_2),\ c(Q_1,qQ_2)=c_2(q)c(Q_1,Q_2).$$
These relations must hold for all primes dividing $\stufe$, so we must have
$$c(Q_1,Q_2)=c_1(Q_1)c_2(Q_2)c(1,1)=c_1(Q_1)c_2(Q_2),$$
proving the theorem.
$\square$
\enddemo

\enddocument